\documentclass[11pt] {article} 
  \usepackage{amsmath}
    \usepackage{amssymb}
    
    \usepackage{pdfsync} 

\newtheorem{example}{Example}[section]
\newtheorem{theorem}{Theorem}[section]
\newtheorem{lemma}{Lemma}[section]
\newtheorem{corollary}{Corollary}[section]

\newtheorem{remark}{Remark}[section]
\newcommand{\eqnsection}{
   \renewcommand{\theequation}{\thesection.\arabic{equation}}
   \makeatletter
   \csname @addtoreset\endcsname{equation}{section}
   \makeatother}
 
   \def\njc{{\bf note Jay's change !!! }}
   
\def \ov{\overline}

\def \be{\begin{equation}}
\def \ee{\end{equation}}
\def \bt{\begin{theorem}} 
\def \et{\end{theorem}}
\def \bl{\begin{lemma}} 
\def \el{\end{lemma}}
\def \bea{\begin{eqnarray}}
\def \eea{\end{eqnarray}}
\def \bas{\begin{eqnarray*}}
\def \eas{\end{eqnarray*}}

\def \al{\alpha}
\def \bb{\beta}
\def \ga{\gamma}
\def \Ga{\Gamma}
\def \de{\delta}
\def \De{\Delta}
\def \ep{\epsilon}
\def \vep{\varepsilon}
\def \la{\lambda} 
\def \La{\Lambda}
\def \ka{\kappa}
\def \om{\omega}
\def \Om{\Omega}

\def \vf{\varphi}
\def \si{\sigma}

\def \th{\theta}

\def \Th{\Theta}

\def \ze{\zeta}

\def \ff{\infty}
\def \wh{\widehat}
\def \wt{\widetilde}

\def \cd{\,\cdot\,}

\def\stl{\stackrel{law}{=}}

\def \FF{{\cal F}}

\def \II{{\cal I}}

\def \RR{{\cal R}}

\def \TT{{\cal T}}

\def \({\left(}
\def \){\right)}

\def \nn{\nonumber}
\def \Proof{\noindent{\bf Proof $\,$ }}

\def \bc{\begin{center} }
\def \ec{\end{center} }
\def \bs{\begin{slide} }
\def \es{\end{slide} }

\def\square{{\vcenter{\vbox{\hrule height.3pt
        \hbox{\vrule width.3pt height5pt \kern5pt
           \vrule width.3pt}
        \hrule height.3pt}}}}
\def\qed{{\hfill $\square$ \bigskip}}

\eqnsection
\begin{document}

\title{Sample path properties of permanental processes}

     \author{  Michael B. Marcus\,\, \,\, Jay Rosen \thanks{Research of     Jay Rosen was partially supported by  grants from the National Science Foundation and the Simons Foundation.   }}
\maketitle

\begin{abstract}    
  Let $X_{\al}=\{X_{\al}(t),t\in \TT\}$,  $\al>0$, be an $\al$-permanental process with kernel $u(s,t)$. We show that     $X^{1/2}_{\al}$ is a subgaussian process with respect  to the metric
\begin{equation}
   \si(s,t)= (u(s,s)+u(t,t)-2(u(s,t)u(t,s))^{1/2})^{1/2} .\nn
   \end{equation}
This allows us to use the vast literature on sample path properties of subgaussian processes to extend these properties to $\al$-permanental processes. Local and uniform moduli of continuity are obtained as well as the behavior of the processes at infinity. Examples are given of permanental processes with kernels that are the potential density of transient L\'evy processes that are not necessarily symmetric, or with kernels of the form  
\begin{equation}
   \wt    u(x,y)= u(x,y)+f(y),\nn
   \end{equation}
where $u$ is the potential density of a  symmetric transient Borel right process  and $f$ is an  excessive function for the process.

  \end{abstract}

\bibliographystyle{amsplain}

  \section{Introduction}

 An $R^{n}$ valued   $\al$-permanental random variable $ (X_{\al }(1),\ldots, X_{\al}(n))$ is  a non-negative random variable with Laplace transform  
\begin{equation}
   E\(e^{-\sum_{i=1}^{n}s_{i}X_{\al}(i)}\) 
 = \frac{1}{ |I+US|^{ \al}},   \label{int.1} 
 \end{equation}
for some  $n\times n$ matrix $U$  and diagonal matrix $S$ with positive entries $ s_{1},\ldots,s_{n} $, and $\al>0$.     We refer to the matrix $U$ as the kernel of $(X_{\al }(1),\ldots, X_{\al}(n))$.     

An $\al$-permanental process $X_{\al}=\{X_{\al}(t),t\in \TT\}$   is a stochastic process that has finite dimensional distributions that are $\al$-permanental random variables.  In this paper $\TT$ is   usually $R^{+}$,  a    subset of $R^{+}$ or  $\mathbb N$, the set of integers.

 An $\al$-permanental process $X_{\al}$ is determined by a kernel $\{u(s,t),s,t 	\in \TT\}$ with the property that for all $t_{1},\ldots,t_{n}$ in $\TT$,  $\{u(t_{i},t_{j}),i,j\in [1,n]\}$ is the kernel of the $\al$-permanental random variable $(X_\al(t_{1}),\ldots,X_\al(t_{n}))$.  To avoid trivialities we restrict our attention to kernels with the property that for all $\de>0$, $\sup_{t\le \de }u(t,t)>0$.
 
 \  An extensive class of examples of kernels of permanental processes is given by  $  U=\{u(s,t),s,t\in \TT\}$ when $U$ is the potential density of a transient Markov process with state space $\TT$,    with respect to some $\si$-finite measure $m$ on $\TT$  and $u(s,t)$ is finite for all $s,t\in\TT$.
 In this case $  U$ is the kernel of an $\al$-permanental process $X_{\al}$ for all $\al>0$; see  \cite[Theorem 3.1]{EK}, and   for all   distinct $(t_{1},\ldots,t _{n})$ in $\TT$, the $n\times n$ matrix    $\{u(t_{i},t_{j})\}_{i,j=1}^{n}$ is invertible;  (see \cite[Lemma A.1]{MRnec}).  We refer to these   permanental processes  as  associated  $\al$-permanental processes because they are associated with the  transient Markov process.   
 
 \medskip	The following observation about bivariate $\al$-permanental random variables is the key to the results   in this paper: Suppose that  $(X_{\al}(s),X_{\al}(t))$ is an $\al$-permanental random variable with  kernel 
   \begin{equation} U_{s,t}=
\left(
\begin{array}{ cc}
u(   s,s)   &  u(   s,t)  \\
u(  t,s)   &   u(   t,t)
\end{array}
\right).  \label{1.8nint}
   \end{equation}
 It follows from \cite[p. 135]{VJ} that $u(s,s)$, $u(t,t)$ and $u(s,t)u(t,s)$ are all greater than or equal to 0.    Furthermore, one can see from (\ref{int.1}) that $(X_{\al}(s),X_{\al}(t))$  also has the  symmetric  kernel
 \begin{equation} \wt U_{s,t}=
\left( 
\begin{array}{ cc}
u(   s,s)   &  (u(s,t)u(t,s))^{1/2}  \\
(u(s,t)u(t,s))^{1/2}  &   u(   t,t)
\end{array}
\right).\label{1.8pint}
   \end{equation}
and 
 \begin{equation}
|\wt U_{s,t}|=| U_{s,t}|=u(s,s)u(t,t)-u(s,t)u(t,s)\geq 0, \label{jss.2}
 \end{equation}
  (see (\ref{3.50xx})), so that the symmetric  kernel $\wt U_{s,t}$ is   positive definite. 
 
 Let $X_{\al}$  be an $\al$-permanental process with kernel $\{u(s,t),s,t 	\in \TT\}$. It follows from (\ref{jss.2}) and the inequality between arithmetic and geometric means that  
 \begin{equation}
u(s,s)+u(t,t)-2(u(s,t)u(t,s))^{1/2}\geq 0. \label{arith-geom}
 \end{equation}
 We then define   the function $\{\si(s,t),s,t\in \TT\}$  by  
\begin{equation}
   \si(s,t) = \(u(s,s)+u(t,t)-2(u(s,t)u(t,s))^{1/2}\)^{1/2}\label{1.2}\ge 0. 
   \end{equation}
 We refer to  $\{\si(s,t),s,t\in \TT\}$  as the sigma function of $X_{\al}$. Note that although we don't require that $u(s,t)$ is symmetric, $\si(s,t)$ obviously is symmetric.  When $X_{\al}$ is an associated  $\al$-permanental process then$\{\si(s,t),s,t\in \TT\}$  is a metric,  \cite[Lemma 4.2]{KMR}.  \label{p3}
      
\medskip	When $u(s,t)$ is symmetric and is a kernel that  determines a $1/2$-permanental process,  $Y_{1/2}=\{Y_{1/2}(t),t\in \TT\}$,   then $Y_{1/2}\stl \{G^{2}(t)/2,t\in \TT\}$ where $\{G  (t) ,t\in \TT\}$ is a mean zero Gaussian process with covariance $u(s,t)$. In this case
\begin{equation}
   \si (s,t)=\(E\(G (t)-G (s)\)^{2}\)^{1/2}=\|G (t)-G (s)\|_{2}.\label{1.3}
   \end{equation}
   
    In \cite{KMR} and  \cite{MRsuf} we use these observations   to show that many properties of Gaussian processes which are obtained  using their bivariate   distributions  also hold for $1/2$-permanental processes. In this paper we show how such properties   also hold for $\al$-permanental processes, for all $\al>0$. Theorem \ref{theo-1.3} is the critical step in this work.
  
 	   \medskip	Let $\psi_{2}(x)=\exp(x^{2})-1$ and $\|\cdot\|_{\psi_{2}}$ be the Orlicz norm in the corresponding Banach space. A fundamental relationship in the study of Gaussian processes is   
\begin{equation}
   \| G(s)-G(t)\|_{\psi_{2}} =2 \|G (t)-G (s)\|_{2}\label{1.53qq}.
   \end{equation}
    We extend this result to $\al$-permanental processes.

\begin{theorem}\label{theo-1.3} Let   $X_{\al}=\{X_{\al}(t),t\in \TT\}$ be an $\al$-permanental process with kernel  $  U=\{u(s,t),s,t\in \TT\}$.    Then 
\begin{equation}
   \| X^{1/2}_{ \al}(s)-X^{1/2}_{ \al}(t)\|_{\psi_{2}} \leq  C_{\al}  \si(s,t):= d_{C _{\al} ,\si }(s,t)\label{1.53},
   \end{equation}
for some finite constant  $C_{\al} $. 
 \end{theorem}
 
 As we just stated, when $u(s,t)$ is symmetric and is a kernel that  determines a $1/2$-permanental process,  the process is  $ \{G^{2}(t)/2,t\in\TT\}$, where $\{G  (t) ,t\in \TT\}$ is a mean zero Gaussian process with covariance $u(s,t)$. In this case (\ref{1.53}) gives
 \begin{equation}  
   \| |G(s)|  - |G(t)| \|_{\psi_{2}} \leq  C'  \|G (t)-G (s)\|_{2}.   \end{equation}    
for some absolute constant $C'$,   which is just a bit weaker than (\ref{1.53qq}) because $||G(s)|  - |G(t)||\le |G(s)   - G(t)|$.

\medskip  Theorem \ref{theo-1.3} shows that  $\{d_ {C ,\si} (s,t),  s,t\in R^{+}\}$   dominates $  \| X^{1/2}_{ \al}(s)-X^{1/2}_{ \al}(t)\newline \|_{\psi_{2}}$.    Therefore the proof used in \cite[Theorem 3.1]{KMR} to obtain results  for the uniform modulus of continuity of 1/2-permanental processes extends immediately to  $\al$-permanental processes for which (\ref{1.53}) holds.  Similarly,    the results in   \cite[Theorems 4.1 and 4.2]{MRsuf} on the local modulus of continuity of  $1/2$-permanental processes, extend to all $\al$-permanental processes for which (\ref{1.53}) holds.  These    are given in Theorems \ref{theo-1.1a} and \ref{theo-1.1au} in the Appendix. 

The reason we present these theorems in an appendix   is because
the results for the uniform modulus of continuity and local modulus of continuity in (\ref{2.1sv}) and (\ref{2.1ww}) are quite abstract. However, as we show in \cite[Example 4.1]{MRsuf}, when $T=R^{+} $ and   there exists an increasing  function $\vf$ such that for   all $0\le s,t<\ff$, 
   \begin{equation}
   \si(s ,t )\le \vf ( |t-s|),\label{8.14dw}
   \end{equation}
 where  
   \begin{equation}
  \int_{0}^{1/2}\frac{\vf( u)}{u(\log 1/u)^{1/2}}\,du<\ff,\label{3.21w}
   \end{equation} then
(\ref{2.1sv}) and (\ref{2.1ww}) give  results  for 1/2-permanental processes, that are  the same as familiar results for Gaussian processes, although we don't get the best constants. 

In order to get precise results with the best constants, we examine more closely permanental processes with kernels that allow (\ref{8.14dw}) and (\ref{3.21w}) to be satisfied and obtain stronger versions of (\ref{2.1sv}) and (\ref{2.1ww})  by  generalizing a classical inequality of Fernique. Here are some examples of    results we obtain.  
  
    \bt \label{theo-1.4}  	   Let  $X_{\al}=\{X_{\al}(t ), t\in [0,1]\}$  be an $\al$-permanental process with   kernel $u(s,t)$ and sigma   function $\si (s,t)$ for which   (\ref{8.14dw})  and (\ref{3.21w}) hold  for some function  $\vf(t)$ that  is regularly varying at zero with positive index.
    Then,   
	            \begin{equation}
   \limsup_{ h\to 0}\sup_{0\le t\le h}\frac{|X_{\al} (t)-X_{\al} (0)|}{    \vf(h)( \log\log 1/h)^{1/2}}\le  2 X_{\al}^{1/2}(0)\qquad a.s.\label{1.23e}
   \end{equation} 
   
   If   
 $u(0,0)=0$,   
    and if  in addition to the conditions on $\vf$ above,   $\vf^{ 2}(h)= O(u (h,h))$,   then
         \begin{equation}
   \limsup_{ t\to 0} \frac{ X_{\al} (t) }{    u(t,t) \log\log 1/t }\le  1\qquad a.s.\label{1.24e}
   \end{equation} 
and 
  \begin{equation}
   \limsup_{ h\to 0}\sup_{0\le t\le h}\frac{ X_{\al} (t) }{   u^{*}(h,h) \log\log 1/h }\le  1\qquad a.s.\label{1.25e}
\ee
where
   \begin{equation}
   u^{*}(h,h)=\sup_{t\le h}u(t,t).\label{ustar}
   \end{equation}

    \et

Note that   when $X_{1/2}$ is 1/2 times the square of    Brownian motion, its kernel $u$  satisfies $u(h,h)=h$.  In this case the upper bound  in   (\ref{1.25e}) is well known to be best possible.  

\medskip	The next theorem gives uniform moduli of continuity.

  \bt \label{theo-1.4unif}  	   Let  $X_{\al}=\{X_{\al}(t ), t\in [0,1]\}$  be an $\al$-permanental process with    kernel $u(s,t)$ and sigma   function $\si (s,t)$ for which   (\ref{8.14dw})  and (\ref{3.21w}) hold.    Assume furthermore,  
	that    $\vf(t)$  is regularly varying at zero with positive index.        Then  
         \begin{equation}
   \limsup_{ h\to 0}\sup_{\stackrel{|s-t|\le h}{s,t\in [0,1]}}\frac{|X _{\al} (s)-X _{\al} (t)|}{    \vf(h)( \log  1/h)^{1/2}}\le  2\sup_{t\in [0,1]}X^{1/2} _{\al} (t)\qquad a.s.\label{8.35rwvq}
   \end{equation}
    \et

   \medskip		 We also investigate the behavior of permanental processes at infinity.    
    	The following simple limit theorem is best possible for the squares of some  Gaussian processes with stationary increments.
  
    \bt \label{theo-3.1w}Let $X_{\al}=\{X_{\al}(t),t\in R^{+}\}$ be an $\al$-permanental process with kernel $u(s,t)$ and  sigma  function $\si (s,t)$    that satisfies (\ref{8.14dw})  and (\ref{3.21w}).  Then
    \begin{equation}
   \limsup _{t\to\ff}  \frac{X_{\al}(t )}{u^{*}(t,t)\log t}  \le 1\qquad a.s.\label{9.10wr}
   \end{equation}
 \et

Under additional hypotheses we get the familiar iterated logarithm in the denominator.
  
   \bt \label{theo-3.2q} Under the hypotheses of Theorem \ref{theo-3.1w} assume furthermore that  $u(t,t)$ is regularly varying at infinity with positive index
and     $ \vf^{2}(t) = O(u(t,t))$ as $t\to \ff$.       Then
     \begin{equation}
   \limsup _{t\to\ff}  \frac{X_{\al}(t )}{u(t,t)\log \log t}  \le 1\qquad a.s.\label{3.67}
   \end{equation}
\et

   We   require  (\ref{3.21w}) because we are considering the processes on $R^{+}$ so in order for them to behave well for all $0<t<\ff$ they must be continuous.  We also want to study the behavior of  permanental sequences   $X_{\al}=\{X_{\al}(t_{n}), n\in \mathbb N\}$ in which $\{t_{n}\}$ has no limit points,   or in which $\lim_{n\to\ff}t_{n}=t_{0}$ but  $X_{\al}$ is not continuous at $t_{0}$. Processes of this sort as treated in the next theorem. 
 	  
  \begin{theorem} \label{lem-3.9a} Let $X_{\al}=\{X_{\al}(t_{n}),n\in \mathbb N\}$ be an $\al$-permanental sequence with kernel $u(t_{j},t_{k})$   and  sigma  function $\si (t_{j},t_{k})$. 
  
Assume also that   $\bb^{2}:= \lim_{n\to\ff}\si^{2}(t_{n} ,0)\log n$ exists. 
  If $\bb=\ff$,
  \begin{equation}
      \limsup_{n\to\ff}\frac{ X_{\al} (t_{n} ) }{ \si^{2}(t_{n} ,0)\log n }\le 1\qquad a.s.\label{3.73q}
   \end{equation}
  If $0<\bb<\ff$ then
 \begin{equation}
   \limsup_{n\to\ff} {|X_{\al}(t_{n} )-X_{\al}(0)|} \le \bb^{2}+2\bb X_{\al}^{1/2}(0)\qquad a.s.\label{3.71}
   \end{equation}
 If $ \bb=0$,
 \begin{equation}
   \limsup_{n\to\ff}\frac{|X_{\al}(t_{n} )-X_{\al}(0)|}{ \si (t_{n} ,0)(\log n)^{1/2}}\le  2  X_{\al}^{1/2}\label{3.72}(0)\qquad a.s.
   \end{equation}
     \end{theorem}

The reader may wonder why (\ref{1.53}) is given in terms of  $X_{\al}^{1/2} $ when we are really concerned with the permanental processes    $X_{\al}$.
The reason is that  for $X_{\al}$ we only have
 \be
  \| X _{ \al}(s)-X _{ \al}(t)\|_{\psi_{1}} \leq  C'_{\al} \si(s,t)\label{1.6},
   \end{equation}
for some finite constant  $C'_{\al}$, where 
  $\psi_{ 1}(x)=\exp(x )-1$ and $\|\cdot\|_{\psi_{1}}$ is the Orlicz norm in the corresponding Banach space. Sufficient conditions for the continuity and boundedness of processes satisfying (\ref{1.6}) are weaker than those for processes satisfying (\ref{1.53}). See, e.g., \cite[Theorem 11.4]{LT}.    (It is much easier to obtain (\ref{1.6}) than   (\ref{1.53}), we do this in \cite[Lemma 3.2 and Corollary 3.1]{KMR}.)

\medskip	 We apply the results above in the following examples:

 	 \begin{example} \label{Levy}{\rm {\bf FBMQ$^{\ga,\bb}$}   Let $Z_{\ga,\bb}=\{Z_{ t},t\in R_{+}\}$ be a L\'evy process with characteristic function
\begin{equation}
 Ee^{i\la Z_{\ga,\bb}}=e^{-\psi_{\ga,\bb}(\la)t},\label{10.1}
 \end{equation}
 where  
 \begin{equation}
\psi_{\ga,\bb} (\la) = |\la|^{\ga+1}(1-i\bb \mbox{ sign }(\la)\,\tan (\ga\pi/2))\label{9.1},
\end{equation}
 for   $0<\ga<1$ and $|\bb|\le 1$. 
 Consider the transient Markov process    that is $Z_{\ga,\bb}$ killed at the first time  it hits zero and  let $u_{T_{0};\ga,\bb}(x,y)$ denote its zero potential.  Therefore, $u_{T_{0};\ga,\bb}(x,y)$ is also     the kernel of $\al$-permanental processes for all $\al>0.$ These permanental processes    belong to the class of processes   $\mbox{FBMQ}^{\ga ,\bb}$. (See  page \pageref{page 1} and \cite[Section 5.1]{KMR} for the explanation of this notation.)
 
 \begin{theorem}\label{theo-1.8new} Let $Y_{\al;\ga,\bb}=\{Y_{\al;\ga,\bb}(x),x\in R^{+}\}$ be an $\al$-permanental process with kernel $u_{T_{0};\ga,\bb} $. Then 
 for all $T>0$
   \begin{equation}
   \limsup_{ h\to 0}\sup_{\stackrel{|s-t|\le h}{s,t\in [0,T]}}\frac{|Y_{\al:\ga,\bb} (s)-Y_{\al:\ga,\bb}(t)|}{  (  h^{\ga}  \log  1/h)^{1/2}}\le  2\sqrt 2 ( (1+|\bb|) C_{\ga,\bb})^{1/2}\(\sup_{t\in T}Y_{\al:\ga,\bb} (t)\)^{1/2} \label{8.35rw5}
   \end{equation}   
 almost surely,  
       \begin{equation}
  2 C_{\ga,\bb} \le \limsup_{ h\to 0}\sup_{0\le t\le h}\frac{Y_{\al:\ga,\bb} (t) }{    h^{\ga} \log\log 1/h }\le   2 (1+|\bb|) C_{\ga,\bb}    \qquad a.s.\label{8.35sw5}
   \end{equation} 
and
     \begin{equation}
   \limsup _{t\to\ff}  \frac{Y_{\al;\ga,\bb}(t )}{ t^{\ga} \log \log t} = 2C_{\ga,\bb}\qquad a.s.\label{3.67ss},
   \end{equation}
where
\begin{equation}
 C_{\ga,\bb}={-\sin \left( (\ga+1)\frac{\pi}{2} \right) 
\Gamma( -\ga ) \over \pi (1+\bb^{2}  \,\tan^{2} ( (\ga+1)\pi/2))}>0.
 \end{equation}
  \end{theorem}
  
     }\end{example}
     
  \begin{example}\label{Levy-2} {\rm   We also consider the transient Markov process    that is $Z_{\ga,\bb}$ killed at the end of an independent exponential time with mean $\rho$.  Let $u_{\rho;\ga,\bb}(x,y)$ denote its zero potential.  Therefore, $u_{\rho;\ga,\bb}(x,y)$ is also     the kernel of $\al$-permanental processes for all $\al>0$.   The next theorem is similar to   Theorem   \ref{theo-1.8new} except   that $Y_{\al; \ga,\bb}\equiv 0$, whereas $Y_{\al;\rho,\ga,\bb}>0$ almost surely.
  
  \begin{theorem}\label{theo-1.8exp} Let $ Y_{\al;\rho,\ga,\bb}=\{Y_{\al;\rho,\ga,\bb}(x),x\in R^{+}\}$ be an $\al$-permanental process with kernel $u_{\rho;\ga,\bb} $. Then 
 for all $T>0$
   \bea
 &&  \limsup_{ h\to 0}\sup_{\stackrel{|s-t|\le h}{s,t\in [0,T]}}\frac{|Y_{\al:\rho,\ga,\bb} (s)-Y_{\al:\rho,\ga,\bb}(t)|}{  (  h^{\ga}  \log  1/h)^{1/2}}\\
 &&\qquad\qquad\le  2\sqrt 2 ( (1+|\bb|) C_{\ga,\bb})^{1/2}\(\sup_{t\in T}Y_{\al:\rho,\ga,\bb} (t)\)^{1/2},  \nn
   \eea   
 almost surely,  
       \begin{equation}
  \limsup_{ h\to 0}\sup_{0\le t\le h}\frac{|Y_{\al:\rho,\ga,\bb} (t) - Y_{\al:\rho,\ga,\bb} (0)| }{  (  h^{\ga} \log\log 1/h )^{1/2}}\le 2\sqrt 2 ( (1+|\bb|) C_{\ga,\bb})^{1/2}  Y_{\al:\rho,\ga,\bb} (0)  \quad a.s. 
   \end{equation} 
and  
  \begin{equation}
   \limsup _{t\to\ff}  \frac{Y_{\al;\rho,\ga,\bb}(t )}{\log t } = D_{\rho,\ga,\bb} \qquad a.s.\label{3.67ss0},
   \end{equation}

   \begin{equation}
  D_{\rho,\ga,\bb}= \frac{1}{\pi}\int_{0}^{\ff}\frac{ \RR e(\rho+\psi_{\ga,\bb}(\la))}{|\rho+\psi_{\ga,\bb}(\la)|^{2}}\,d\la.\label{1.44mm}
   \end{equation}

  \end{theorem}

   }\end{example}   
   
  The potentials in Examples \ref{Levy} and  \ref{Levy-2} are not symmetric when $\bb\ne 0$. Prior to this paper there have not been examples of permanental processes that do not have symmetric kernels other than this case, i.e., when the kernel is the potential of a L\'evy process. The next theorem shows how we can modify a very large class of symmetric potentials  so that they are no longer symmetric but are still kernels of permanental processes.

    \bt
   \label{theo-borelN}
        Let $S$ a be locally compact set with a countable base. 
Let $X\!=\!
(\Om,  \FF_{t}, X_t,\th_{t},P^x
)$ be a transient symmetric Borel right process with state space $S$  and continuous strictly positive  potential   densities  $u(x,y)$ with respect to some $\si$-finite measure $m$ on $S$.  
 Then for any excessive function   $f$ of $X$ and $\al>0$,
\begin{equation}
   \wt u(x,y)= u(x,y) +f(y),\qquad x,y\in S,\label{1.10}
   \end{equation}
is the kernel of an $\al$-permanental  process.
\et

   A function $f$ is said to be  excessive for $X$ if $  E^{x}\(f(X_{t})\)\uparrow  f(x)$ as $t\to 0$ for all  $x\in S$.     
  It is easy to check that for any positive measurable function $h$, 
\begin{equation}
   f(x)=\int u(x,y) h(y) \,dm(y)=E^{ x} \(   \int_{0}^{ \ff}h\(   X_{t}\)\,dt\)\label{potdef}
   \end{equation}
  is excessive for $X$. 
Such a function $f$   is called  a potential function for $X$.  All the potential functions  considered  in this paper are continuous. This is discussed   at the beginning of Section \ref{sec-Borel}.  	 We describe other  excessive functions,   some of which are not potentials, in the next two examples.

   \begin{example}\label{ex-1.2} {\rm    Let   $ \ov B=\{ \ov B_{t},t\in R^{+} \}$ be Brownian motion  killed after an independent exponential time with parameter $\la^{2}/2$, or,  equivalently,   with mean $2/\la^{2}$. The process $\ov B$ has potential densities,  
\begin{equation}
\ov u(x,y) ={e^{-\la |y-x|}\over \la},\qquad x,y\in  R^{ 1}  \label{expke.13},
\end{equation} 
and $f$ is excessive for $\ov  B$  if and only if  $f$ is positive and   is a  $\la^{2}/2$--superharmonic function, \cite[p. 659]{Doob}. In particular,  $f\in C^{2}$  is excessive for $\ov  B$  if and only if it is positive and $f''(x)\leq \la^{2} f(x)$   for all $x\in R^{1}$.    Examples of such functions are   $e^{r x}$ for $|r|\leq \la$ and $q+|x|^{ \bb}$ for $\bb\geq 2$ and $q\ge q_{0}$, where $q_{0}$ depends on $\bb$ and $\la$. 
 It follows   Theorem \ref{theo-borelN} that for   functions  $f$ that are excessive for $\ov  B$,   
\be
  \wh u_{f}(s,t)=e^{-\la |s-t|}+f(t),\qquad s,t\in [0,1],\label{1.33ss}
  \ee  
  is the kernel of an $\al$-permanental process on $ [0,1]$ for all $\al>0$. 
 
\medskip	  The function  $e^{-\la |s-t|}$, $s,t\in R^{+}$, is also the covariance of a time changed   Ornstein-Uhlenbeck   process.

\begin{theorem}\label{theo-1.9mm} Let $\ov X_{\al,f}= \{\ov X_{\al,f}(t), t\in  [0,1]\}$ be an $\al$-permanental process with kernel $\wh u_{f}(s,t)$  where  $f\in C^{2}$   and   is excessive for $\ov  B$. Then     
 \begin{equation}
  \limsup_{ h\to 0} \sup_{ 0\le t\le h}\frac{|\ov X_{\al,f} (t)-\ov X_{\al,f}(0)|}{  (  h  \log \log 1/h)^{1/2}}\le 2\sqrt {2 \la   } \,  \ov X^{1/2}_{\al,f} (0) \quad a.s.\label{1.43rwx}
  \end{equation}   
and,  for all $T>0$,  
  \begin{equation}
  \limsup_{ h\to 0}\sup_{\stackrel{|s-t|\le h}{s,t\in [0,T]}}\frac{|\ov X_{\al,f} (s)-\ov X_{\al,f}(t)|}{  (  h  \log  1/h)^{1/2}}\le 2\sqrt   {2 \la   } \(\sup_{t\in T}\ov X_{\al,f} (t)\)^{1/2},  \label{8.35rww}
  \end{equation}   
  almost surely.
   \end{theorem} 
 
 }\end{example}
  
 	  We now obtain an upper bound for the behavior at infinity.

 \begin{theorem}\label{theo-1.9inf} Let $\wh X_{\al,f}= \{\wh X_{\al,f}(t), t\in  [0,\ff)\}$ be an $\al$-permanental process with kernel $\wh u_{f}(s,t)$ in (\ref{1.33ss}) but for $s,t\in [0,\ff)$, where  $f\in C^{2}$   and   is excessive for $\ov  B$. Then    
 \begin{equation}
  \limsup_{ t\to \ff} \frac{ \wh  X_{\al,f} (t)  }{  (  1+f(t))\log t }\le  1   \qquad a.s.,\label{1.43rw}
  \end{equation}    
  and when $\lim_{t\to \ff}f(t)=0$,  
  \be
      \limsup_{ t\to \ff} \frac{ \wh  X_{\al,f} (t)  }{  \log t }= 1 \qquad a.s.\label{1.43ro}
  \end{equation}  
  In particular this holds if $f$ is a potential for $\ov  B$, with  $h\in L_{+}^{1}\(0,\ff\)$.   (See   (\ref{potdef}).)
    \end{theorem}

 \begin{example} \label{ex-1.3}{\rm  
Let $\wt B=\{\wt B_{t},t\in R^{+}\}$ be Brownian motion killed the first time it hits $0$. $\wt B$ has state space $D=(0,\ff)$   and potential densities  
\begin{equation}
\bar u(x,y) =2(x\wedge y),\hspace{.2 in}x,y>0.\label{expke.14}
\end{equation}

 \begin{theorem} \label{theo-1.10mm}
For any positive  concave function $f$ on $\(   0,\ff\)$, and $\al>0$,
   \begin{equation}
 \wt u_{f}(s,t)=s  \wedge t+f(t),\qquad s,t>0.\label{1.39mm}
  \end{equation}
is the kernel of an $\al$-permanental processes, $  \wt Z_{\al,f}  = \{\wt Z_{\al ,f}(t),t>0\} $,   
and 
  \be  1\le  \limsup_{t\to \ff}\frac{\wt Z_{\al,f  } (t) }{ t\, \log\log t}\leq 1+C_{0}\qquad a.s.\label{1.25ss},
  \end{equation}
  where  $C_{0}= \lim_{t\to \ff} f(t)/t$,   which is  necessarily   finite .  
  
   If $f$ is a potential for   $\wt B$, with  $h\in L_{+}^{1}\(0,\ff\)$, (see   (\ref{potdef})), then  $f(t)=o(t)$ and
    \be   \limsup_{t\to \ff}\frac{\wt Z_{\al,f  } (t) }{ t\, \log\log t}=1\qquad a.s.\label{1.25ssj}
  \end{equation}
\end{theorem}
        }\end{example}
        
  The next theorem  describes the behavior of  $  \{\wt Z_{\al ,f}(t),t>0\}$ as $t\to 0$   when $\lim_{t\to 0}f(t)/t=\ff$.     
      
        \bt\label{theo-cup}
Let    $  \wt Z_{\al,f}$ be the    $\al$-permanental process defined in   Theorem \ref{theo-1.10mm}, 
for     $f$   a positive concave function on   $[0,\ff)$,  with the additional property that for some $\de>0$ 
\begin{equation}
 \int_{0}^{\de}\frac{1}{f(u)(\log 1/u)^{1/2}}\,du<\ff\label{ddd}.
   \end{equation} 
Then there exists a coupling of $ \wt Z_{\al,f}$ with a   gamma random variable $\xi_{\al,1}$, with shape $\al$ and scale $1$, such that
\begin{equation}
\label{cup.0}
\lim_{t\to 0}\frac{\wt Z_{\al ,f}(t)}{f(t)}=\xi_{\al,1}\hspace{.2 in}a.s.\end{equation}
\et

With an additional condition on $f$ we can describe the behavior of $ \wt Z_{\al,f}$ near $\xi_{\al,1}$ more precisely.

        \bt\label{theo-cupx}
Let    $  \wt Z_{\al,f}$ be the    $\al$-permanental process defined in   Theorem \ref{theo-1.10mm}, 
for    $f$  a positive concave function on $(0,\ff)$ that is regularly varying at  zero  with index less than 1. Then
	     \begin{equation}
   \limsup_{ t\to 0} \frac{|\wt Z^{1/2}_{\al ,f}(t)-(f(t)\xi_{\al,1})^{1/2}|}{    (t\log\log 1/t)^{1/2}}\le 1\qquad a.s.\label{cup.4}, 
   \end{equation} 
   where $\xi_{\al,1}$ is defined in Theorem \ref{theo-cup}.
\et

\medskip	  Other examples are given in the body of this paper.

Theorem \ref{theo-1.3} is the main result in this paper. It   deals, simply,  with pairs of  permanental random variables. Let $Z_{\al}=(Z_{\al}(1),Z_{\al}(2))$ be an $\al$-permanental random variable with  kernel
 \begin{equation} K=
\left(
\begin{array}{ cc}
b   &  \ga  \\
 \ga   &   a 
\end{array}
\right).  \label{1.8}
   \end{equation}
     We point out in the paragraph  containing (\ref{1.8nint}) that $a$, $b$ and $|K|\geq 0$. In addition 
  we can take $\ga\geq 0. \footnote{ An $\al$-permanental random can have more than one  kernel. In particular  if it has kernel (\ref{1.8}) then it also has the kernel with $\ga$ replaced by $-\ga$.}$ We point out in the paragraph  containing (\ref{1.2})   that   $\si=(a+b-2\ga)^{1/2}\ge 0$.

\medskip	 The next theorem is the main ingredient in the proof of Theorem \ref{theo-1.3}. 
        
   \bt \label{theo-3.2} Let   $X_{\al}=(X_{\al}(1),X_{\al}(2))$ be an $\al$-permanental random variable with  kernel $K$ in (\ref{1.8}).
  Set 
  \begin{equation}
 \si^{2}= {a+b-2\ga} .\label{1.3ax}
   \end{equation}
Then     for all $\la\ge 1$,  
    \be   
    P\(\frac{|X_{\al}^{1/2}(1) -X_{\al}^{1/2}(2)|}{\si}\geq \la\)   \leq     C_{\al}\,\la^{(4\al-2) \vee 0 }\,e^{-\la^{2}}, \label{1.9fa} 
   \ee 
  for some constant $C_{\al}$, depending only on $\al.$   
    \et  
  
 Note that  by inequality between arithmetic and geometric means and the fact that $|K|\ge 0$,  
  \begin{equation}
 \frac{a+b}{2}\ge (ab)^{1/2}\ge \ga.\label{sigzero}
   \end{equation}
 Hence if $\si=0$ we have equality throughout (\ref{sigzero}) which implies that $a=b$ and $|K|=0$.  
  We point out prior to the statement of Lemma \ref{lem-2.1n} that this implies that   $X_{\al}^{1/2}(1) =X_{\al}^{1/2}(2)$ almost surely.  We take the quotient $0/0$ in (\ref{1.9fa}) to be 0. (See Lemma \ref{lem-2.1n}.)
      
 \medskip	We are interested in  local moduli of continuity and rate of growth of $X_{\al}^{1/2}$ and $X_{\al}$. It follows from Theorem \ref{theo-3.2} that when   $\al\le 1/2$,  
  \begin{equation}
  P\(\frac{|X_{\al}^{1/2}(s) -X_{\al}^{1/2}(t)|}{\si(s,t)}\geq \la\)\leq    C_{\al}  e^{-\la^{2}},  \label{1.9d}
   \end{equation}
  for some constant $C_{\al}$, depending only on $\al.$ 
 It is well known; (see, e.g.,   \cite[Lemma 5.1.3]{book}), that
 for
$\wt G=\{\wt G(t), t\in \TT\}$,   a mean zero Gaussian process with covariance $u(s,t)/ 2$, 
   \begin{equation}
   P\(\frac{| G(s) - G(t)|}{\si (s,t)}\geq \la\)\leq      e^{-\la^{2} } \label{1.9f}.
   \end{equation}
    Comparing (\ref{1.9d})  and (\ref{1.9f}) it is clear that upper bounds for the rates of growth of Gaussian processes that are proved solely by using the Borel-Cantelli Lemma applied to increments of the process, should also hold for the square roots of permanental processes when    $\al\le 1/2$.   This is the case even when $\al>1/2$  and the exponential in (\ref{1.9d}) is multiplied by a power of $\la$.
 
 \medskip	An extensive treatment of moduli of continuity of Gaussian processes is given in \cite[Chapter 7]{book}.	However the  proofs use properties of Gaussian processes  more sophisticated than (\ref{1.9f}). Earlier treatments give many of the same results and are  obtained solely from  (\ref{1.9f}).  In Lemma \ref{lem-8.2a} we   modify   
an inequality of Fernique, (see, e.g., 
\cite[Chapter IV.1, Lemma 1.1]{JM}), and   use it to obtain the results in Theorems \ref{theo-1.4}--\ref{lem-3.9a}. 

\medskip The proofs of   Theorems \ref{theo-1.3} and \ref{theo-3.2} are given in Section \ref{sec-2}. The proofs of  Theorems \ref{theo-1.4},   \ref{theo-1.4unif} and   \ref{theo-1.9mm}  are given in Section \ref{sec-8}.   Theorem  \ref{theo-1.8exp} and the proof of the upper bounds in  Theorem \ref{theo-1.8new}  are proved in Section \ref{newest}. Theorems \ref{theo-3.1w}--\ref{lem-3.9a}   and  the proof of upper bounds in   Theorems   \ref{theo-1.9inf}--\ref{theo-cupx} are proved in Section \ref{sec-4}.   The proof  of   Theorem  \ref{theo-borelN} is given in Section \ref{sec-Borel}.  The proof of the lower bounds in  Theorems \ref{theo-1.8new},    \ref{theo-1.9inf} and   \ref{theo-1.10mm} are given in Section \ref{sec-7}.   The extension of 
\cite[Theorem 3.1]{KMR} and  \cite[Theorems 4.1 and 4.2]{MRsuf}   to all $\al$-permanental processes for which (\ref{1.53}) holds is given  in  Section \ref{Appendix}.

     \section{Proofs of  Theorems \ref{theo-1.3} and \ref{theo-3.2}}\label{sec-2}

\medskip	Let   $\xi_{\al,v}$  denote a gamma random variable with 
 probability density function 
 \be
f(\al,v;x) = \frac{v^{\al} x^{\al-1} e^{-v x}}{\Gamma(\al)}, \qquad    x \geq 0 \text{ and } \al,v > 0,\label{1.2q}
 \ee and equal to 0 for $x\le0$, where
$\Gamma(\al) =\int_{0}^{\ff}  x^{\al-1} e^{-  x}\,dx 
$ is the gamma function.  

  It is easy to see that $\xi_{\al,v}\stackrel{(law)}{=}v\xi_{\al,1}$ and 
\begin{equation}
E\(e^{-s \xi_{\al,1}}\)={1 \over (1+s)^{\al} },\label{int.1m}
\end{equation}
so that $\xi_{\al,1}$ is an $\al$-permanental random variable with  kernel $1$. We also note that
\begin{equation}
    P\( \xi_{\al,1}^{1/2}  \geq \la\)=P\( \xi_{\al,1}   \geq \la^{2}\)={1 \over \Gamma(\al)}\int_{\la^{2}}^{\ff}x^{\al-1} e^{- x}\,dx.\label{int.1ma}
\end{equation}
Let $\la\geq 1$. For $\al\leq 1$, (\ref{int.1ma})  is bounded by
\begin{equation}
 {\la^{2(\al-1)} \over \Gamma(\al)}\int_{\la^{2}}^{\ff} e^{- x}\,dx={\la^{2(\al-1)} e^{- \la^{2}}\over \Gamma(\al)}.\label{2.4m}
\end{equation}
For $\al>1$, we see by integration by parts that (\ref{int.1ma}) is bounded by $C_{\al}\la^{2(\al-1)} e^{- \la^{2}}$. Thus, for $\la\geq 1$,
\begin{equation}
    P\( \xi_{\al,1}^{1/2}  \geq \la\)\leq C_{\al}\la^{2(\al-1)} e^{- \la^{2}}.\label{int.1mb}
\end{equation}

\medskip	Let $X_{\al}=(X_{\al}(1),X_{\al}(2))$ be as in Theorem  \ref{theo-3.2}. Then 
\begin{equation}
|I+KS|=1+s_{1}b+s_{2}a+s_{1} s_{2}|K|.\label{}
\end{equation}
If $|K|=0$
it follows   from (\ref{int.1}) that
\begin{equation}
   E\(e^{-(s_{1}X_{\al}(1)+s_{2}  X_{\al}(2))} \) 
 = \frac{1}{ (1+s_{1}b+s_{2}a )^{ \al}}.  \label{int.1easy} 
 \end{equation}
 Let $ (X_{\al}(1),X_{\al}(2))\stackrel{(law)}{=}(b,a)\xi_{\al,1}$. Then
 \bea
   E\(e^{-(s_{1}X_{\al}(1)+s_{2}  X_{\al}(2))} \) 
 &=&E\(e^{-(s_{1}b+s_{2}a)\xi_{\al,1}} \)  \\
  &=&\frac{1}{\Ga(\al) }\int_{0}^{\ff}e^{-(s_{1}b+s_{2}a)x} x^{\al-1}e^{-x}\,dx\nn\\
  &=&\frac{1}{\Ga(\al) }\int_{0}^{\ff}e^{-(1+s_{1}b+s_{2}a)x} x^{\al-1} \,dx\nn\\
  &=&\frac{1}{ (1+s_{1}b+s_{2}a )^{ \al}}\nn
  \eea
 Therefore, when $|K|=0$,  $ (X_{\al}(1),X_{\al}(2))\stackrel{(law)}{=}(b,a)\xi_{\al,1}$, and since $|K|=0$ implies that $ab=\ga^{2}$ we see that $\si^{2}=(\sqrt{b}-\sqrt{a})^{2}$. 
 
 When $a\ne b$ \begin{equation}
\frac{|X_{\al}^{1/2}(1) -X_{\al}^{1/2}(2)|}{\si}\stackrel{(law)}{=}\frac{|\sqrt{b}-\sqrt{a}|\xi_{\al,1}}{|\sqrt{b}-\sqrt{a}|}=\xi_{\al,1}.\label{int.1mj}
\end{equation}
Thus, when $|K|=0$ and $a\ne b$, (\ref{1.9fa}) follows from (\ref{int.1mb}).

When $|K|=0$ and $a= b$, $X_{\al}(1)=X_{\al}(2)$ so the numerator on the left-hand side of  (\ref{int.1mj}) is equal to zero. Of course, $\si$ is also equal to zero. We take $0/0=0$ and get $(\ref{1.9fa})$ in this case also.

For general  $|K|$, assume that $\ga \le (a\wedge b)/2$. It follows that $\si^{2}\geq a\vee b$.
 Consequently\bea
  &&  P\(\frac{|X_{\al}^{1/2}(1) -X_{\al}^{1/2}(2)|}{\si}\geq \la\)   \leq     P\(\frac{X_{\al}^{1/2}(1)  }{\si}\geq \la\)+ P\(\frac{X_{\al}^{1/2}(2)  }{\si}\geq \la\)\nn\\
   &&\qquad  \leq     P\(\frac{b^{1/2}\xi_{\al,1}^{1/2}    }{\si}\geq \la\)+P\(\frac{a^{1/2}\xi_{\al,1}^{1/2}    }{\si}\geq \la\)\le 2C_{\al}\la^{2(\al-1)} e^{- \la^{2}}, \label{2.10mm}
\eea
where we use the facts that  $ X_{\al}(1) \stackrel{(law)}{=}b\xi_{\al,1}$, $ X_{\al}(2) \stackrel{(law)}{=}a\xi_{\al,1}$  and (\ref{int.1mb}).
Thus we get (\ref{1.9fa}) when  $\ga\le (a\wedge b)/2$.

\medskip	We collect these observations into the following lemma:

\bl\label{lem-2.1n}
   If $|K|=0$ or if $\ga\le (a\wedge b)/2$, Theorem  \ref{theo-3.2} holds. In particular Theorem  \ref{theo-3.2} holds when $\ga=0$.
   \el

    We continue with the proof of Theorem \ref{theo-3.2}. From now on we assume that $\ga>0$ and $|K|>0$ which implies that  $\si>0$.   We use the probability distribution of $ (X_{1},X_{2} )$ that is given in \cite[Theorem 1.1]{Marcus}       in terms of  
     \be 
     K^{-1}={1 \over |K|}
     \begin{pmatrix}a&-\ga\\
     -\ga&b \end{pmatrix}.
     \ee

      \begin{theorem}\label{theo-1.5} Let $X=(X_{1},X_{2} )$ be an $\al$-permanental random variable with  kernel $K$.   The    probability density function of $X$ is  
      \be  \wt  g(\al ;{(x_{1},x_{2})}) \label{1.11b} =  \frac{ \ga^{1-\al}}{  \Ga(\al)\de}   \frac{
 \II_{\al-1}\(2\ga\sqrt{x y}/\de\,\)} {(xy)^{(1-\al)/2}}e^{-( a x /\de+b y/\de)} 
   \ee 
on $R^{2}_{+}$,  and zero elsewhere, where   $\de=|K|$ and 
   \begin{equation}
   \II_{\nu}(z)= \sum_{n=0}^{\ff}\frac{1}{\Ga(n+\nu+1)\,n!}\(\frac{z}{2}\)^{2n+\nu}\label{3.2}
   \end{equation}
   is the modified Bessel function of the first kind.
 \end{theorem}

  We use the notation   $E(\xi;A)=E(\xi I_{A})$ for sets $A$.

 	  \bl \label{lem-3.1} Let $X=(X_{1},X_{2} )$ be an $\al$-permanental random variable with  kernel $K$ as in  Theorem  \ref{theo-1.5}.   Then when   $\ga>0$ and $\rho\ge \sqrt 2$,
	     \be
E\(\exp\(\frac{\rho}{\si}\( {X^{1/2}_{1}-X^{1/2}_{2}}\)\); X_{1} X _{2}\ge (\de/2\ga)^{ 2}\,\)\le  C_{\al}\(\frac{\rho}{\sqrt 2}\)^{(4\al-2)\vee 0}e^{\rho^{2} /4},\label{1.9lo}
   \ee
   for some constant $C_{\al}$ depending only on $\al$.  
 \el

\medskip	\Proof   By Theorem \ref{theo-1.5}  
  \bea 
\lefteqn{    E\(\exp\(\frac{\rho}{\si}\(X^{1/2}_{1}-X^{1/2}_{2}\)\); X_{1} X _{2}\ge (\de/2\ga)^{ 2}\,\)\label{1.11}=\frac{ \ga^{1-\al}}{  \Ga(\al)\de}}\\
     && \int_{0}^{\ff}\!\!\int_{0}^{\ff}e^{ \rho (x^{1/2} -y^{1/2} )/\si}  \frac{
 \II_{\al-1}\(2\ga\sqrt{x y}/\de\,\)} {(xy)^{(1-\al)/2}}e^{-( a x /\de+b y/\de)}1_{\{   x y \ge (   \de/2\ga)^{ 2}\}}\,d x \,d y  \nn.
  \eea 
     We make the change of variables  $x=u^{2} /2$, $y=v^{2} /2$ and see that (\ref{1.11})
       \begin{equation}
 =\frac{ (2\ga)^{1-\al} }{ \Ga (\al)\de} \int_{0}^{\ff} \!\!\int _{0}^{\ff}  (uv)^ \al e^{ \rho'  (u-v )}  
 \II_{\al-1}\(\ga uv/\de\) e^{-( au^{2}+bv^{2})/2 \de}1_{\{   \ga uv/\de\ge 1 \}}\,d u\,dv  \label{1.27j}
   \end{equation}
where $\rho'=\rho/(\sqrt2 \si)$. 
  Over the range $ \ga uv/\de\ge 1$, we bound 
     \begin{equation}
   \II_{\al-1}(\ga uv/\de)\le C' _{\al}\frac{\de^{1/2}e^{\ga uv/\de}}{\sqrt{2  \ga uv}} \label{2.8t}
      \end{equation}
 for some  constant $C'_{\al}$, depending on $\al$;  see \cite[8.451.5]{GR}. Therefore    (\ref{1.27j})  
               \begin{equation}
\le \frac{ (2\ga)^{(1/2)-\al}C'_{\al}}{  \Ga (\al)\de^{1/2}}   \int _{0}^{\ff}\!\!\int _{0}^{\ff}(uv) ^{\al-1/2} e^{ \rho'  (u-v )}  
 e^{-( au^{2}-2\ga uv+bv^{2})/2 \de}1_{\{   \ga uv/\de\ge 1 \}}\,d u\,dv  \label{3.36ss}.
   \end{equation} 
   \begin{equation}
   \le \frac{ (2\ga)^{(1/2)-\al}C'_{\al}}{  \Ga (\al)\de^{1/2}}   \int _{-\ff}^{\ff}\!\!\int _{-\ff}^{\ff}|uv| ^{\al-1/2} e^{ \rho'  (u-v )}  
 e^{-( au^{2}-2\ga uv+bv^{2})/2 \de} \,d u\,dv  \label{3.36sss}
   \end{equation}
 	Consider  
\begin{equation}
 \rho'  (u-v )-\frac{  au^{2}-2\ga uv+bv^{2}}{2 \de}
   \end{equation}	
Let   $w=(u,v)$ and  $z=\rho' (1,-1)$, and write this as 	
\begin{equation}
   (z,w)-\frac{(w,K^{-1} w)}{ 2}.\label{2.17n}
   \end{equation}
Let 	$w=s+Kz$. With this substitution (\ref{2.17n}) is equal to
\begin{equation}
     \frac{ (z,Kz)}{2}-\frac{(s,K^{-1} s)}{ 2}.\label{2.18n}
   \end{equation}
It follows from this that if we make the change of variables $u=s_{1}+ K(z)_{1}$ and $v=s_{2}+ K(z)_{2}$, and $\al\ge 1/2$,   (\ref{3.36sss}) is  less than or equal to  
   \begin{equation}
 \frac{  C''_{\al}}{ 2\pi |K|^{1/2}}   \int _{-\ff}^{\ff}\!\!\int _{-\ff}^{\ff}\(|s_{1}+ K(z)_{1}||s_{2}+ K(z)_{2}| \)^{\al-1/2}    
 e^{-(   s,K^{-1} s )/2  } \,d s_{1}\,ds_{2}\,e^{(z, Kz)/2}  \label{3.36n},
\ee
in which
\begin{equation}
C_{\al}'' =  \frac{ (2\pi)(2\ga)^{(1/2)-\al}C'_{\al}}{  \Ga (\al) }. 
   \end{equation}
(Recall that $\de=|K|$.)

Note that (\ref{3.36n})
\begin{equation}
  \le C_{\al}'' \(E [\(|\xi_{1}|+ |K(z)_{1}|\)(|\xi_{2}|+ |K(z)_{2}|) ]^{\al-1/2}\)e^{(z, Kz)/2},\label{2.21n}
   \end{equation}
where $(\xi_{1},\xi_{2})$ is a mean zero 2-dimensional Gaussian random variable with covariance matrix $K$. Furthermore, by Lemma \ref{lem-2.3},  
\begin{equation}
   |K(z)_{1}|=\frac{\rho}{\sqrt {2}}\frac{|a-\ga|}{\si}\le \frac{\rho\, a^{1/2}}{\sqrt {2}}\quad\mbox{and}\quad   |K(z)_{2}|=\frac{\rho}{\sqrt {2}}\frac{|b-\ga|}{\si}\le \frac{\rho\, b^{1/2}}{\sqrt {2}}.
   \end{equation}
Therefore, 
\begin{equation}
     \(|\xi_{1}|+ |K(z)_{1}|\)(|\xi_{2}|+ |K(z)_{2}|)  \le   (ab)^{1/2}   \(\Big|\frac{\xi_{1}}{a^{1/2}}\Big|+ \frac{\rho}{\sqrt {2}}\) \(\Big|\frac{\xi_{2}}{b^{1/2}}\Big|+ \frac{\rho}{\sqrt {2}}\). \label{2.23n}
   \end{equation}
Note that $\xi_{1}/a^{1/2}$ and $\xi_{ 2}/b^{1/2}$ have variance 1. Using this, (\ref{2.23n}) and the Cauchy Schwartz Inequality, and the fact that
\begin{equation}
   (z,Kz)=(\rho'\si)^{2}=\frac{\rho^{2}}{2},\label{2.24nx}
   \end{equation}
we see that (\ref{2.21n})
\begin{equation}
  \le C_{\al}''(ab) ^{(2\al-1)/4} E \(|\eta|+\frac{\rho}{\sqrt 2}\)^{2\al-1 } e^{\rho^{2}/4}.\label{2.21nx}
   \end{equation}
where $\eta$ is a standard normal random.

Since $\al\ge 1/2$, this is  
\bea
 && \le C_{\al}''2^{2\al-1}(ab) ^{(2\al-1)/4} \(E \(|\eta|^{2\al-1} \)  +  \( \frac{\rho}{\sqrt 2}\)^{2\al-1 }\) e^{\rho^{2}/4} \label{2.24n}\\
  && \le C_{\al}''2^{2\al-1}(ab) ^{(2\al-1)/4}E \(|\eta|^{2\al-1} \) \(1  +  \( \frac{\rho}{\sqrt 2}\)^{2\al-1 }\) e^{\rho^{2}/4}.\nn
   \eea
   
   Note that 
   \begin{equation}
   C_{\al}''2^{2\al-1}(ab) ^{(2\al-1)/4}=\frac{    2^{\al+1/2}\,\pi \,C'_{\al}}{\Ga(\al)}\(\frac{ (ab)^{1/2} }{\ga}\)^{( \al-1/2)}. 
   \end{equation} 
 Therefore, if $\ga\ge (ab)^{1/2}/\sqrt 2$, the left-hand side of (\ref{1.11})  
   \be  
   \le \frac{ 2 ^{(5\al+1)/4}\,\pi\,C_{\al}'}{ \Ga(\al)}  E \(|\eta|^{2\al-1} \) \(1  +  \( \frac{\rho}{\sqrt 2}\)^{2\al-1 }\) e^{\rho^{2}/4}.\label{2.28n}
   \ee

 \medskip	We now obtain an upper bound for the left-hand side of (\ref{1.11})  when $\al\ge 1/2$ and $\ga<(ab)^{1/2}/\sqrt 2$. 
   To begin note that left-hand side of (\ref{1.11})  is bounded by the integral in (\ref{3.36ss}), in which the integrand is restricted to the region, $\ga uv/\de\ge 1$. On this region
\begin{equation}
   \frac{1}{\ga^{\al-1/2}}\le \(\frac{uv}{\de}\) ^{\al-1/2}.
   \end{equation}
 Using this inequality we can bound the integral in (\ref{3.36ss}) by 
  \begin{equation}
  \frac{ D_{\al} }{    \de^{\al-1/2}2\pi|K|^{1/2}}   \int _{0}^{\ff}\!\!\int _{0}^{\ff}(uv) ^{2\al-1} e^{ \rho'  (u-v )}  
 e^{-( au^{2}-2\ga uv+bv^{2})/2 \de}1_{\{   \ga uv/\de\ge 1 \}}\,d u\,dv,  \label{2.25mm}
   \end{equation}
   where
  \begin{equation}
 D_{\al}=   \frac{ 2\pi C'_{\al} }{  2^{\al-1/2}\Ga (\al) }.  
   \end{equation} 
  Following the argument from (\ref{3.36ss})--(\ref{2.24n}), and in particular focusing on (\ref{3.36n}) we see that
   \bea
\lefteqn{ \frac{ 1 }{   2\pi|K|^{1/2}}   \int _{0}^{\ff}\!\!\int _{0}^{\ff}(uv) ^{2\al-1} e^{ \rho'  (u-v )}  
 e^{-( au^{2}-2\ga uv+bv^{2})/2 \de}1_{\{   \ga uv/\de\ge 1 \}}\,d u\,dv  \nn}\\
 &&  \le  4^{2\al-1}(ab) ^{ \al-1/2}E \(|\eta|^{2(2\al-1)} \) \(1  +  \( \frac{\rho}{\sqrt 2}\)^{2(2\al-1) }\) e^{\rho^{2}/4}.
   \eea
Therefore (\ref{2.25mm})
\begin{equation}
   \le   \frac{ 4^{2\al-1}\pi C'_{\al} }{   2^{\al-1/2} \Ga (\al) } \(\frac{ ab }{\de}\)^{( \al-1)/2}E \(|\eta|^{2(2\al-1)} \) \(1  +  \( \frac{\rho}{\sqrt 2}\)^{2(2\al-1) }\) e^{\rho^{2}/4}.\label{2.33n}
   \end{equation}
Now, note that $\ga<(ab)^{1/2}/\sqrt 2$ implies that $\de>(ab/2)$. Consequently, (\ref{2.33n})
\be
  \le  \frac{ 4^{2\al-1}\pi C'_{\al} }{     \Ga (\al) }  E \(|\eta|^{2(2\al-1)} \) \(1  +  \( \frac{\rho}{\sqrt 2}\)^{2(2\al-1) }\) e^{\rho^{2}/4}.\label{2.33nn}
   \end{equation}
 Using (\ref{2.28n}) and (\ref{2.33nn}) we see that  when $\al\ge 1/2$,   $\ga>0$ and $\rho\ge \sqrt 2$,
	     \be
E\(\exp\(\frac{\rho}{\si}\( {X^{1/2}_{1}-X^{1/2}_{2}}\)\); X_{1} X _{2}\ge (\de/2\ga)^{ 2}\,\)\le  C_{\al}\(\frac{\rho}{\sqrt 2}\)^{4\al-2}e^{\rho^{2} /4},\label{1.9lox}
   \ee
   for some constant $C_{\al}$ depending only on $\al$.

\medskip	 
To  obtain an upper bound for the left-hand side of (\ref{1.11})  when  $\al< 1/2$ we follow the argument from   (\ref{3.36ss})--(\ref{3.36n}) and use (\ref{2.24nx}) to see that it is 
\bea
 && \le C_{\al}'' \(E \Big |\frac{1}{ ( \xi_{1} +  K(z)_{1} ) (\xi_{2} +  K(z)_{2})}\Big | ^{ 1/2-\al}\)e^{\rho^{2} /4}\label{2.21nn}\\
 && \le C_{\al}'' \(E \Big |\frac{1}{ ( \xi_{1} +  K(z)_{1} ) (\xi_{2} +  K(z)_{2})}\Big | ^{ 1/2-\al}\)e^{\rho^{2} /4}\nn\\
 &&\le C_{\al}''\left[ E \( \(  \frac{1}{ |\xi_{1} | }\)^{ 1-2 \al } \)E \( \(  \frac{1}{ |\xi_{2} |}\)^{ 1-2 \al } \)\right] ^{1/2}e^{\rho^{2} /4}\nn,
   \eea
     where, as we point out above, $\xi_{1}$ is mean zero normal random variable  with variance  $a$ and $\xi_{2}$ is mean zero normal random variable  with variance  $b$. Let   $\eta$ be a standard normal random variable. The last line of (\ref{2.21nn}) is equal to   
  \bea
  && \frac{  2  ^{(3/2)-\al} \ga ^{(1/2)-\al}\pi C_{\al}' }{ \Ga (\al) }\frac{1}{(ab)^{(1/4)-(\al/2)}}E \(\(  \frac{1}{ |\eta |}\) ^{ 1-2 \al }\)e^{\rho^{2} /4}\label{2.34n}\\
  &&\qquad\le \frac{  2  ^{(3/2)-\al}  \pi C_{\al}' }{ \Ga (\al) } E \(\(  \frac{1}{|\eta |}\) ^{ 1-2 \al }\)e^{\rho^{2} /4},\nn
   \eea
  where we use the fact that $\ga<\sqrt{ab}.$
    Using   (\ref{1.9lox}) and (\ref{2.34n})  we get (\ref{1.9lo}). \qed
    
 We use the next lemma in the proof of Lemma \ref{lem-3.2}.

   \begin{lemma} \label{lem-3.3}  For $K$   and $\si^{2}$ in (\ref{1.8}) and (\ref{1.3ax}),   
   \begin{equation}
|K| \le (a\wedge b){ \si^{2}} .\label{2.2ww}
   \end{equation}
 \end{lemma}
 
\Proof    To prove (\ref{2.2ww})
 we first  show that $ |K| \le a\si^{2} $. This is 
 \begin{equation}
ab-\ga^{2}\le a^{2}+ab-2a\ga,
   \end{equation}
which is equivalent to $(a-\ga)^{2}\geq 0$.
The same argument works with $a$ replaced by $b$.\qed

 	   \begin{lemma}\label{lem-3.2}  For   $\ga>0$ and $\la\ge 1$,   
   \begin{equation}
 P\(\frac{ |X^{1/2}_{\ 1}-X^{1/2}_{ 2}|}{\si}\ge \la;\,\, X_{1} X _{2}\le (\de/2\ga)^{ 2}\,\)\leq \wh C_{\al}    \la^{(2\al-2)  } \,\, e^{-  \la^{2} },\label{3.27a}
   \end{equation}
where $\wh C_{\al}$ is a constant depending only on $\al$.
 \end{lemma}

\Proof 	Let  $(   U,V):=2^{1/2}(X^{1/2}_{ 1},X^{1/2}_{ 2})$ so that the left-hand side of (\ref{3.27a}), without the absolute value signs,  can be written as 
\begin{equation}
   P\(U-V\ge \sqrt2 \la\,\si;\,\, UV\leq \de/\ga \)\leq   P\(U\ge \sqrt2 \la\,\si;\,\, UV\leq \de/\ga \).\label{3.27m}
   \end{equation} 
Using (\ref{1.27j}) we see that  
  \be 
 \frac{ (2\ga)^{1-\al}}{\Ga(\al)\de } \label{3.31}\\
(uv)^{\al} \II_{\al-1}\( \ga {uv}/\de\,\) e^{-( au^{2}  +bv^{2})/2\de} 
   \ee  
  is the joint probability density function of $(   U,V)$.  
 To find an upper bound for (\ref{3.27m}) we note that by (\ref{3.2})
    \begin{equation}
   \II_{\al-1}(w )\le C_{\al}w^{\al-1},     \label{2.30}
   \end{equation}
 for $w\leq 1 $,  where $C_{\al}$ is a constant depending on $\al$. 
     With this substitution (\ref{3.27m}) is less than or equal to 
      \bea
  &&  \frac{C_{\al}2^{1-\al}}{\Ga(\al)\de^{\al}}\int_{ \sqrt2 \la\,\si}^{\ff}\int _{0}^{\ff}(uv)^{ 2\al -1} e^{ -( au^{2}  +bv^{2})/2\de} \,dv\,du.\label{3.33}
   \eea
 
  We have 
 \be \int _{0}^ {\ff}   v^{ 2\al -1}e^{- bv ^{2}  /2\de   } \,dv=\sqrt\frac{\pi}{2}\(\frac{\de}{b}\)^{\al}E \(| Z |^{2\al -1 }\),\label{1.37qa} 
   \ee 
 where $Z$ is a normal random variable with mean zero and variance 1. In addition  
  \bea
 \int_{\sqrt2\la\si}^{\ff} u^{ 2\al -1}e^{- a u^{2}/2\de }  \,du&=&   \(\frac{\de}{a}\)^{\al}  \int_{(a/\de)^{1/2} \sqrt2\la\si}^{\ff} s^{ 2\al -1}\,e^{- s^{2}/2}  \,ds\nn\\
 &\le &   \(\frac{\de}{a}\)^{\al}  \int_{  \sqrt2\la }^{\ff}  s^{ 2\al-1}   \,e^{- s^{2}/2}  \,ds,\label{1.37q}
  \eea
 since $a\si^{2}/\de\ge 1$ by Lemma \ref{lem-3.3}.   Setting $x=s^{2}/2$ in (\ref{int.1ma}), it follows from (\ref{int.1mb}) that for all 
$\la\geq 1$ and $\al>0$,
 \begin{equation}
    \int_{  \sqrt2\la }^{\ff}  s^{ 2\al-1}   \,e^{- s^{2}/2}  \,ds\le C_{\al}\la^{2(\al-1)} e^{- \la^{2}}.\label{2.35j}
  \end{equation} 

  Using (\ref{1.37qa})--(\ref{2.35j}) in (\ref{3.33}) 
  we see that (\ref{3.33})  is bounded by 
   \begin{equation}
    D''_{\al} \sqrt\frac{\pi}{2}\(\frac{\de}{ab}\)^{\al}E \(| Z |^{2\al -1 }\)\la^{2\al-2}e^{-\la^{2}},
   \end{equation}
     where $D''_{\al}$ is a constant depending only on $\al$.  
Since   $\de/(ab)<1$, we get (\ref{3.27a}). (We can include the absolute value sign by multiplying the  probability on the right by 2.)\qed

\noindent{\bf Proof of Theorem \ref{theo-3.2} } When $\ga=0$   or $|K|=0$ this follows from Lemma \ref{lem-2.1n}. 

 Now suppose that $\ga>0$ and $||K>0$. We write
 \bea
   P\(X^{1/2}_{1}-X^{1/2}_{2}\geq\si \la\)& \le &  P\( X^{1/2}_{1}-X^{1/2}_{2}\geq\si \la; \,\, X_{1} X _{2}\le (\de/2\ga)^{ 2}\)  \label{3.37}\\
  & &\quad +   P\( X^{1/2}_{1}-X^{1/2}_{2}\geq \si\la; \,\, X_{1} X _{2}\ge (\de/2\ga)^{ 2}\).  \nn
   \eea
Using Lemma \ref{lem-3.2} we see that we have the upper bound in (\ref{1.9fa}) for the first probability on the right-hand side of  the inequality in (\ref{3.37}).

By Lemma \ref{lem-3.1}  the second term on the right-hand side of the inequality  
\bea
  && \le e^{-\rho\la}    E\(\exp\(\frac{\rho\(X^{1/2}_{1}-X^{1/2}_{2}\)}{\si}\) ;\,\, X_{1} X _{2}\ge (\de/2\ga)^{ 2}\)\\
  &&  \le C_{\al}\(\frac{\rho}{\sqrt 2}\)^{(4\al-2)\vee 0} e^{-\rho\la} e^{\rho^{2} /4}\nn.
   \eea
 Taking $\rho=2\la $  the upper bound in (\ref{1.9fa}) for the second probability on the right-hand side of the inequality in (\ref{3.37}). \qed

\medskip	We need the next two lemmas in the proof of Theorem \ref{theo-1.3}.
\begin{lemma}\label{lem-1.4} If
\begin{equation}
     P\(|Z|\ge \la\)  \le     K e^{-\la^{2}} \label{1.50}
   \end{equation}
then 
\begin{equation}
   \|Z\|_{\psi_{2}}\le c^{*},
   \end{equation}
   where $c^{*}$ is the value of $c>1$ such that
   \begin{equation}
   K^{1/c^{2}} \frac{c^{2}}{c^{2}-1 }=2\label{3.45}.
   \end{equation} 
  
 \end{lemma}

\Proof For  $c>1$  and   $   K=e^{y^{2}_{0}}$, 
\bea
 &&  E\(\exp\(\frac{Z^{2}}{c^{2}} \)-1\)\label{1.51}\\&&\qquad=-\int_{0}^{\ff}\(e^{\la^{2}/c^{2}}-1\)\,d  P\(|Z|\ge \la\) = \int_{0}^{\ff}P\(|Z|\ge \la\)\,d e^{\la^{2}/c^{2}}  \nn\\
      && \qquad\leq \int_{0}^{y_{0}}\,d e^{\la^{2}/c^{2}}+K \int_{y_{0}}^{\ff}e^{-\la^{2}}\,d e^{\la^{2}/c^{2}}\nn\\
      &&\qquad= e^{y^{2}_{0}/c^{2}}-1+K \int_{y_{0}}^{\ff}e^{-\la^{2}}\,d e^{\la^{2}/c^{2}}.\nn
   \eea 
 We have
 \bea
  \int_{y_{0}}^{\ff}e^{-\la^{2}}\,d e^{\la^{2}/c^{2}}  &=&\frac{1}{c^{2}(   1-1/c^{2})} \int_{y_{0}}^{\ff}e^{-\la^{2}(   1-1/c^{2})}2(   1-1/c^{2})\la \,d\la\nn\\
      &= &\frac{e^{-y_{0}^{2}(   1-1/c^{2})}}{c^{2}(   1-1/c^{2})},  
 \eea 
 so that (\ref{1.51})
 \begin{eqnarray}
 &&= e^{y^{2}_{0}/c^{2}}-1+ e^{y^{2}_{0}}\frac{e^{-y_{0}^{2}(   1-1/c^{2})}}{c^{2}(   1-1/c^{2})} \\
 &&= e^{y^{2}_{0}/c^{2}}-1+  \frac{ e^{y^{2}_{0}/c^{2}}}{c^{2}(   1-1/c^{2})} = e^{y^{2}_{0}/c^{2}}\frac{c^{2}}{c^{2}-1 }-1 \nn,
 \end{eqnarray}
    which gives (\ref{3.45}). 
   \qed

  	\begin{lemma}\label{lem-1.5} If 
	 \begin{equation}
  P\(|Z|\ge \la\) \le    K (  \la ^{n}+1)e^{-\la^{2}},\qquad \la\ge 0,
   \end{equation}
for $n>0$, then 
\begin{equation}
   \|Z\|_{\psi_{2}}\le c^{*},
   \end{equation}
   where $c^{*}$ is the value of $c>1$ such that
   \begin{equation}
 e^{y^{2}_{0}/c^{2}} +\frac{ K }{c^{2}} \int_{y_{0}}^{\ff}\la (  \la ^{n}+1)e^{-\la^{2}(1-1/c^{2})}\,d \la =2\label{3.45q},
   \end{equation} 
 for $y_{0}$,  the solution of 
 \begin{equation}
   K   y_{0} ^{n}e^{-y_{0}^{2}}=1.
   \end{equation}
 \end{lemma}

\Proof For $c>1$  and $   K  y_{0} ^{n}e^{-y^{2}_{0}}=1$,  
\bea
 &&  E\(\exp\(\frac{Z^{2}}{c^{2}} \)-1\)\\&&\qquad=-\int_{0}^{\ff}\(e^{\la^{2}/c^{2}}-1\)\,d  P\(|Z|\ge \la\) = \int_{0}^{\ff}P\(|Z|\ge \la\)\,d e^{\la^{2}/c^{2}}  \nn\\
      &&\qquad= e^{y^{2}_{0}/c^{2}}-1+\frac{ 2K}{c^{2}} \int_{y_{0}}^{\ff} \la (  \la ^{n}+1)e^{-\la^{2}(1-1/c^{2})}\,d \la  \nn,
   \eea which gives (\ref{3.45q}).
\qed

 	\noindent{\bf   Proof of Theorem \ref{theo-1.3} }
By hypothesis $(X_{\al}(s),X_{\al}(t))$ is an $\al$-permanental random variable with  kernel
 \begin{equation} K_{s,t}=
\left(
\begin{array}{ cc}
u(   s,s)   &  u(   s,t)  \\
u(  t,s)   &   u(   t,t)
\end{array}
\right)  \label{1.8n}.
   \end{equation}
We point out in (\ref{1.8pint}) that $(X_{\al}(s),X_{\al}(t))$  also has  the  symmetric   kernel
 \begin{equation} U_{s,t}=
\left(
\begin{array}{ cc}
u(   s,s)   &  (u(s,t)u(t,s))^{1/2}  \\
(u(s,t)u(t,s))^{1/2}  &   u(   t,t)
\end{array}
\right)  \label{1.8p},
   \end{equation}
    and that the function
     $\si(s,t)$ corresponding to this is as given in (\ref{1.2}). Therefore,
it follows from   Theorem \ref{theo-3.2} that  for $\la\geq 1$ 
  \be
    P\(\frac{|X^{1/2}_{\al}(s)- X^{1/2}_{\al}(t) |}{\si(s,t)}\geq \la\) \leq    
       C_{\al}\,\la^{(4\al-2) \vee 0 }\,e^{-\la^{2}}.\label{2.58q}
   \ee
   for some absolute constant $ C_{\al}$ that depends only on $\al$. 
  Theorem \ref{theo-1.3} follows from  (\ref{2.58q})  and Lemmas \ref{lem-1.4} and \ref{lem-1.5}.\qed
  
  The next lemma is used in the proof of  Lemma \ref{lem-3.1}.

   \begin{lemma}\label{lem-2.3}  When $|K|>0$, 
   \be
\frac{ |a-\ga| }{\si}\le a^{1/2}\qquad\mbox{and}\qquad \frac{ |b-\ga| }{\si}\le b^{1/2}.\label{plac.1}
\ee
\end{lemma}

\Proof  Since $a$ and $b$ are interchangeable, it suffices to prove the first inequality in (\ref{plac.1}). Suppose that $a-\ga\ge 0$. Then if, in addition,  $b>\ga$, 
\begin{equation}
  \frac{  a-\ga  }{\si}=    \frac{  a-\ga  }{(a +b-2\ga)^{1/2}}\le (a-\ga)^{1/2}.\label{2.30r}
   \end{equation}
(Note that since $\si>0$, this holds when $a=\ga$.)
 
 Next, suppose that $a> \ga$ and $b<\ga$. Then, since  
   
 \begin{equation}
 \frac{d}{d\ga}\( \frac{  (a-\ga)^{2}  }{ a +b-2\ga } \)=2\frac{\(\ga-b\)\(a-\ga\)}{(a +b-2\ga)^{2} }>0\label{2.31r},
   \end{equation}
 we see that because $\ga\le (ab)^{1/2}$,  
  \bea
      \frac{  a-\ga  }{\si}&=&  \frac{  a-\ga  }{(a +b-2\ga)^{1/2}} \le     \frac{  a-(ab)^{1/2}  }{(a +b-2(ab)^{1/2})^{1/2}} \label{2.32r}\\
   &=&  \frac{  a-(ab)^{1/2}  }{a^{1/2}-b^{1/2}}= a^{1/2}.\nn   
   \eea
If $a<\ga$, then $b> \ga$, and
  \begin{equation}
      \frac{   \ga-a  }{\si}=   \frac{   \ga -a }{(a +b-2\ga)^{1/2}} .\label{2.32t}
   \end{equation}
   Then, since   by (\ref{2.31r}),
  \begin{equation}
 \frac{d}{d\ga}\( \frac{  ( \ga-a)^{2}  }{ a +b-2\ga } \)>0.\label{2.33r} 
   \end{equation}
  It follows,  as in (\ref{2.32r}), that  
 \begin{equation}
      \frac{   \ga-a  }{\si}\le  \frac{  (ab)^{1/2} -a }{(a +b-2(ab)^{1/2})^{1/2}}=\frac{  (ab)^{1/2} -a }{b^{1/2}-a^{1/2}}=a^{1/2} .\label{2.32tt}
   \end{equation}
\qed

\section{  Upper bounds for the local moduli of continuity and rate of growth of  permanental processes, I}\label{sec-8}
 
All the results in this section follow from the next lemma which is a modification of an   inequality of Fernique as presented in
\cite[Chapter IV.1, Lemma 1.1]{JM}.

   \begin{lemma}  \label{lem-8.2a} Let    $Y=\{Y(t), t\in R^{+}\}$ be a stochastic process.  For $s,t\in [0,S]$ and $a\ge 0$,   let
   \begin{equation}
  F(a)=\sup_{s,t\in [0,S]}P\(\frac{|Y(s)- Y(t) |}{\si (s,t)}\geq a\) \label{8.15z}  
     \end{equation} 
 for some positive function $\si(s,t)$ on $[0,S]\times [0,S]$.
 Assume furthermore, that there exist an increasing  function $\vf$ such that for $S>0$ and all $0\le s,t\le S$,
   \begin{equation}
   \si(s ,t )\le \vf ( |t-s|).\label{8.14}
   \end{equation}
Let  $n$ be an integer greater than 1. Then
  \bea
 &&P\(   \sup_{t\in [0,S]}|Y(t)-Y(0)|> a \vf(S)+\sum_{p=1}^{\ff }\ka (p)\vf(S/ n(p))\)\label{8.15}\\
 &&\qquad
   \le n^{2}F(a)+\sum_{p=1}^{\ff }n^{2}(p) F(\ka(p))  ,\nn
   \eea
   where $n(p)  = n^{2^{p}}$ and $\ka$ is a positive function with $\ka \ge1$.  
    \end{lemma}
    
    To give the reader some idea of where we are heading we mention that in using (\ref{8.15}) we take both $a$ and $S$ to depend on $n$ in such a way that the right-hand side of  (\ref{8.15}) is a converging sequence in $n$.  This enables us to use the Borel--Cantelli Lemma to get upper bounds for the limiting behavior of $|Y(t)-Y(0)|$.    
       
\medskip	\noindent{\bf Proof of Lemma \ref{lem-8.2a} }  Consider \cite[Chapter IV.2, Lemma 1.1]{JM}. This lemma is proved for  $S=1$. (It also assumes that $Y$ is a Gaussian process, but only uses $F(a)$ as defined  in (\ref{8.15z}), in which case the right-hand side of (\ref{8.15z}) is independent of $s$ and $t$.) 

 Let us first assume that $S=1$. 
   The only other thing in \cite[Chapter IV.1, Lemma 1.1]{JM} that might be confusing is the term $\al$ which is
  \begin{equation}
  \sup_{t\in [0,1]}\si(0,t)\le \vf(1).
   \end{equation}
   Thus we get (\ref{8.15}) with $S=1$. In particular we require that        \begin{equation}
   \si(s ,t )\le \vf ( |t-s|),\qquad \forall  \,0\le s,t\le 1\label{8.14a}.
   \end{equation}
Now replace $\vf(\cd)$ by $\vf(S\cd)$. We have for all $0\le u,v\le S$
 \begin{equation}
   \si(u,v )\le \vf (S |(u/S)-(v/S)|)  ,\qquad \forall  \,0\le    (u/S),(v/S) \le 1\label{8.14aa}.
   \end{equation}
If we replace $\vf(\cd)$ with $\vf(S\cd)$ we get (\ref{8.15}) for arbitrary $S$.    Here we also use the fact that the  right-hand side of (\ref{8.15z}) is defined for all $s,t\in [0,S]$.\qed

 Lemma \ref{lem-8.2a} can used to find upper bounds for the local and uniform moduli of continuity of $Y$ and the behavior of $Y(t)$ as $t\to\ff$.   This is done in  \cite[Theorem 1.3, Chapter IV.2]{JM} for the local and uniform moduli of continuity of Gaussian processes  $G=\{G(t),t\in R^{+}\}$, but the same proof also  gives the behavior of $G(t)$ as $t\to\ff$.    We generalize  \cite[Theorem 1.3, Chapter IV.2]{JM} in the case of local modulus of continuity and behavior at infinity by applying  Lemma \ref{lem-8.2a}  to stochastic processes with the property that  
  \begin{equation}
   F(\la) \le  C_{m} \(\la^{m} +1\)e^{-\la^{2}},\label{8.28}
   \end{equation}
where $m\ge 0$ and $C_{m}$ is a   constant.   

     \begin{lemma} \label{theo-8.1} Let $Y=\{Y(t), t\in R^{+}\}$ be a stochastic process for which (\ref{8.28})  holds.  Assume that there exist an increasing  function $\vf$ such that for   all $0\le s,t<\ff$,
   \begin{equation}
   \si(s ,t )\le \vf ( |t-s|),\label{8.14d}
   \end{equation}
 where  
   \begin{equation}
  \int_{0}^{1/2}\frac{\vf( u)}{u(\log 1/u)^{1/2}}\,du<\ff.\label{3.21}
   \end{equation}
   Consider the condition    
   \begin{equation}
   \quad \lim_{\ga\to 1}\frac{\vf(\ga V)}{\vf( V)}=1.\label{8.33}
   \end{equation}

Assume that  (\ref{8.33})  holds  uniformly in  $0<V\leq V_{0}$ for some $V_{0}\le 1$. For $h$ near zero, set 
   \begin{equation}
   \tau(h)=  \vf(h)(\log\log 1/h)^{1/2}+\frac{1}{\log 2}\int_{0}^{1/2}\frac{\vf(hu)}{u(\log 1/u)^{1/2}}\,du .\label{3.23}
   \end{equation}
   Then  
   \begin{equation}
   \limsup_{ h\to 0}\sup_{0\le t\le h}\frac{|Y(t)-Y(0)|}{\tau(h)}\le \sqrt 3\qquad a.s.\label{8.35}
   \end{equation} 
   Similarly, assume that   (\ref{8.33})   holds  uniformly in  $V\geq V_{0}$ for some $V_{0}\ge 1$.   
        For $T$ large, set  
      \begin{equation}
 \wt  \tau(T)=  \vf(T)(\log\log T)^{1/2}+\frac{1}{\log 2}\int_{0}^{1/2}\frac{\vf(Tu)}{u(\log 1/u)^{1/2}}\,du .\label{3.25}
   \end{equation} Then
   \begin{equation}
   \limsup_{ T\to\ff}\sup_{0\le t\le T}\frac{|Y(t)-Y(0)|}{\wt \tau(T)}\le \sqrt3\qquad a.s.\label{8.37}
   \end{equation}   
 \end{lemma}

 \Proof    We first prove (\ref{8.37}).  Fix  $1<\th<2$  and set   $V=\th^{n}$. Choose $n_{0}$   such that  $\th^{n_{0}}\geq V_{0}$.    Consider Lemma \ref{lem-8.2a}  with $S=\th^{n}$ for $n\geq n_{0}$.  Let $\ep>0$ and take
 \be
 a =((3+\ep)\log\log \th^{n})^{1/2}. 
  \ee  
    Note that by   (\ref{8.28}) there exists a constant  $C $ such that 
   \begin{equation}
F(a )\le 1\wedge    C \frac{(\log n)^{m/2} +1  }{n^{3+\ep }}.\label{3.21q}
   \end{equation}
   Therefore 
    \begin{equation}
n^{2}F(a )\le      \frac{1 }{n^{1+\ep/2}},\label{8.33w}
   \end{equation}
   for all  $n\geq n_{0}$ sufficiently large.

\medskip	Consider the last term in (\ref{8.15}). Take  $\ka(p)=( 3\log n(p))^{1/2}$. We have that 
\begin{equation}
   \sum_{p=1}^{\ff }n^{2}(p) F(\ka(p))  \le  \wt C  \sum_{p=1}^{\ff }\frac{(2^{p}\log n )^{m/2}}{n( p )}\le \frac{1}{n^{3/2}} ,\label{8.39}
   \end{equation}
   for some constant $\wt C$ and all $n$ sufficiently large. 
   Considering (\ref{8.33w}) and (\ref{8.39}) we see   that the right-hand side of (\ref{8.15}) is a converging sequence in $n$.
   
   \medskip	Now consider the sum in the first line of (\ref{8.15}). We have  
   \bea 
  && \sum_{p=1}^{\ff} ( 3\log n(p ))^{1/2}\vf(\th^n/n(p))\label{3.35} \\
    &&\qquad\le \sqrt3  \sum_{p=1}^{\ff}  (  \log n(p ))^{1/2}\vf(\th^n/n(p))\nn\\
      &&\qquad\le \frac{\sqrt3}{\log 2}  \int_{0}^{1/n }  \frac{\vf(\th^nu)}{u(\log 1/u)^{1/2}}\,du \nn.
 \eea
 (This is obtained by replacing the sum by the integral with respect to $p$, from zero to infinity, and making the change of variables $n(p )=1/u$).
 We can now use the Borel-Cantelli Lemma to get 
    \begin{equation}
   \limsup_{ n\to\ff}\sup_{0\le t\le \th^{n}}\frac{|Y(t)-Y(0)|}{ \wt\tau(\th^{n})}\le \sqrt3\qquad a.s.\label{8.47a}
   \end{equation}

 To prove (\ref{8.37})   it suffices to show that
      \begin{equation}
   \limsup_{ T\to\ff}\sup_{0\le t\le T}\frac{|Y(t)-Y(0)|}{\wt\tau(T)}\le\sqrt3 \(1+\ep\)^{2}\qquad a.s.,\label{8.37qm}
   \end{equation}
   for all $\ep>0$.
   
   Fix  $\ep>0$. We show below that we can find    $1<\th_{1}<2$ and  $n_{1}<\ff$ such that 
 \begin{equation}
\frac{ \wt\tau(\th_{1}^{n+1})}{ \wt\tau(\th_{1}^{n })}\le 1+\ep \label{8.37qm1}
 \end{equation}
   for all $n\geq n_{1}$. Furthermore, by (\ref{8.47a})   for all $\om\in \Omega $ for some set $\Omega$ with probability one, we can   find  $n_{2}=  n_{2}(\om)<\ff$  such that  
      \begin{equation}
 \sup_{n\geq n_{2}}\,\,\sup_{0\le t\le \th_{1}^{n}}\frac{|Y(t)-Y(0)|}{ \wt\tau(\th_{1}^{n})}\le \sqrt3(1+\ep)\qquad a.s.\label{8.47aklm}
   \end{equation}  
   Then,    for any   $\th_{1}^{n}<\rho <\th_{1}^{n+1}$ with $n\geq n_{1}, n_{2}$ 
   \bea
   \sup_{0\le t\le \rho }\frac{|Y(t)-Y(0)|}{ \wt\tau(\rho )}\label{8.48}&\le   & \sup_{0\le t\le \th_{1}^{n+1}}\frac{|Y(t)-Y(0)|}{ \wt\tau(\th_{1}^{n})}\\
 & \le   & \sup_{0\le t\le \th_{1}^{n+1}}\frac{|Y(t)-Y(0)|}{ \wt\tau(\th_{1}^{n+1})}\frac{ \wt\tau(\th_{1}^{n+1})}{ \wt\tau(\th_{1}^{n })}\le \sqrt3\(1+\ep\)^{2}\nn, 
   \eea
   where the final inequality uses (\ref{8.37qm1}) and (\ref{8.47aklm}).  This gives (\ref{8.37qm}).
   
   To obtain  (\ref{8.37qm1})  
let
   \begin{equation}
   \II(T)=\int_{0}^{1/2}\frac{\vf( Tu)}{u(\log 1/u)^{1/2}}\,du.
   \end{equation}
 	   We have
         \begin{equation}
 \wt  \tau( \th_{1}^{n+1})= \vf( \th_{1}^{n+1})(\log\log  \th_{1}^{n+1})^{1/2}+\frac{1}{\log 2}   \II(\th_{1}^{n+1}) .\label{3.26}
   \end{equation}
Note that 
 \begin{equation}
 \log\log  \th_{1}^{n+1}=   \log (n+1)+ \log\log  \th_{1} =       \log\log  \th_{1}^{n } +O\(\frac{1}{n}\).\label{3.27}
   \end{equation}

 	Also  
	\bea
&&  \II(\th_{1}^{n+1}) =  \int_{0}^{\th_{1}/2}\frac{\vf(\th_{1}^{n} u)}{u(\log 1/ u+ \log \th_{1})^{1/2}}\,du \label{3.28} \\
  & &\le  \int_{0}^{1/2}\frac{\vf(\th_{1}^{n}  u)}{u(\log 1/ u+\log \th_{1})^{1/2}}\,du\nn+ \int_{1/2}^{\th_{1}/2}\frac{\vf(\th_{1}^{n}  u)}{u(\log 1/ u+\log \th_{1})^{1/2}}\,du.
 \eea
The last integral above  
\begin{equation}
   \le {\vf( \th_{1}^{n+1}/2)}  \int_{1/2}^{\th_{1}/2}\frac{1}{u(\log 1/u)^{1/2} }\,du\le 2 \vf( \th_{1}^{n+1}) (\log \th_{1})^{1/2}   .\label{3.29}
   \end{equation}
     Combining (\ref{3.28}) and (\ref{3.29}) we have  
\begin{equation}
    \II(\th_{1}^{n+1}) \le  \II(\th_{1}^{n})+  2 \vf( \th_{1}^{n+1}) (\log \th_{1})^{1/2}.
   \end{equation}
Using this and (\ref{3.26}) and (\ref{3.27})  and the fact that $  { \vf( \th_{1}^{n+1})}/{\vf( \th_{1}^{n })}\geq 1$ we get  
\bea
  \lefteqn{  \wt  \tau( \th_{1}^{n+1})}\\
    &\le  & \vf( \th_{1}^{n+1}) \( \log\log  \th_{1}^{n } +O\(\frac{1}{n}\)\)^{1/2}+ \frac{\II(\th_{1}^{n})}{\log 2}+ \frac{2\vf( \th_{1}^{n+1})( \log \th_{1})^{1/2} }{  \log 2 }\nn\\
    &\le  &  \frac{ \vf( \th_{1}^{n+1})}{\vf( \th_{1}^{n })}\(\vf( \th_{1}^{n })\( \log\log  \th_{1}^{n }\)^{1/2}+ \frac{\II(\th_{1}^{n})}{\log 2}\right.\\
    &&\hspace{1.5in}+\left. \vf( \th_{1}^{n })\( \frac{2 ( \log \th_{1})^{1/2} }{  \log 2 }+O\(\frac{1}{n^{1/2}}\)\)\)\nn\\
   & =&\frac{ \vf( \th_{1}^{n+1})}{\vf( \th_{1}^{n })}\( \wt  \tau( \th_{1}^{n})+o( \wt  \tau( \th_{1}^{n }))\).\nn
   \eea 
 Or, equivalently,  
   \begin{equation}
\frac{ \wt  \tau( \th_{1}^{n+1})}{ \wt  \tau( \th_{1}^{n})}\leq \frac{ \vf( \th_{1}^{n+1})}{\vf( \th_{1}^{n })}\(1+o( 1_{n })\).\label{8.37qm19}
   \end{equation}
   It follows from    (\ref{8.33})  when    $V\geq V_{0}$ that we can choose $\th_{1}>1$, sufficiently close to $1$, so that 
  \be
  \frac{ \vf( \th_{1}V)}{\vf(V)}\leq 1+\ep/4 \qquad\forall\,   V\geq V_{0} .
  \ee 
We next choose $n_{1}$ so that $ \th_{1}^{n_{1}}\geq V_{0}$, and large enough so that  for all $n\geq n_{1}$, $( 1+\ep/4)\(1+o( 1_{n })\) \leq 1+\ep$ . (Here $\(1+o( 1_{n })\)$ is the expression in (\ref{8.37qm19})).
This completes the proof of   (\ref{8.37qm1}) and hence of (\ref{8.37qm}). 

The statement in 
 (\ref{8.35}) follows similarly by taking  $\th$ less than $1$.  \qed
   
 \begin{remark} {\rm A proof of (\ref{8.35}), the local modulus of continuity, is essentially given  in \cite[Chapter IV, Theorem 1.3]{JM}. But there are some differences with what is given here. Condition \cite[(2.5.9)]{JM} is with regard to a different metric than $\si(s,t)$ but that doesn't matter since what is used in  \cite[Chapter IV, Theorem 1.3]{JM} is (\ref{8.28}). (There is also the requirement that $\vf(2t)\le 2\vf(t)$ in  \cite[Chapter IV, Theorem 1.3]{JM} .) We don't actually prove (\ref{8.35}) here but only point out that it basically the same as the proof of  (\ref{8.37}). The result in  (\ref{8.37}) is not contained in \cite{JM}. Also contained in  \cite[Chapter IV, Theorem 1.3]{JM} is a uniform modulus of continuity of the type given in Theorem \ref{theo-8.1wm}. We do not use it because it doesn't give the constant 1 on the right-hand side of  (\ref{8.35rwv}).

 }\end{remark}
 
   Examples of the processes   $Y$ that we are  studying are the  processes $X_{\al}^{1/2}$ in  Theorem \ref{theo-1.3}.  We see from (\ref{2.58q}) that (\ref{8.28})	holds.   Therefore, we can use Lemma \ref{theo-8.1}.
    	 
 \medskip		 In Theorem \ref{theo-8.1w}  we consider an important class of processes for which we can  lower	the upper bound  in (\ref{8.35}) so that it is best possible.   It uses the next lemma which is an immediate consequence of   (\ref{2.58q}).

	 \bl\label{lem-8.5a} Let  $X_{\al}=\{X_{\al}(t ), t\in [0,1]\}$  be an  $\al$-permanental process with  bounded kernel $u(s,t)$ and sigma   function $\si (s,t)$.  	 Then for any sequence   $\{s_{n}, t_{n}\}$ in $(0,1]\times (0,1]$,   such that $s_{j}\ne t_{j}$ for all $j\in\mathbb N$,
         \begin{equation}
   \limsup_{ n\to\ff} \frac{|X^{1/2}(t_{n})-  X^{1/2}(s_{n}) |}{ \si (s_{n},t_{n}) (\log n )^{1/2} }\le 1\qquad a.s.\label{8.37eea}
   \end{equation}	       
   \el

     	   \bt \label{theo-8.1w}  	   Let  $X_{\al}=\{X_{\al}(t ), t\in [0,1]\}$  be an  $\al$-permanental process with   kernel $u(s,t)$ and sigma   function $\si (s,t)$ for which   (\ref{8.14d})  and (\ref{3.21})  hold for some function  $\vf(t)$ that  is regularly varying at zero with positive index.
   Then,  
	     \begin{equation}
   \limsup_{ h\to 0}\sup_{0\le t\le h}\frac{|X_{\al}^{1/2}(t)-X_{\al}^{1/2}(0)|}{   \vf( h) (\log\log 1/h)^{1/2}}\le 1\qquad a.s.\label{8.35qqq} 
   \end{equation}

 If   
 $u(0,0)=0$,  and if  in addition to the conditions on $\vf$ above,   $\vf^{ 2}(h)= O(u (h,h))$, then  
  \begin{equation}
   \limsup_{ t\to 0} \frac{ X_{\al}^{1/2}  (t) }{(u(t,t)\log\log 1/t))^{1/2} }\le  1  \qquad a.s.\label{8.3jj}
   \end{equation}   
   and
   \begin{equation}
   \limsup_{ h\to 0}\sup_{0\le t\le h}\frac{ X_{\al}^{1/2}(t) }{   (u^{*}(h,h)\log\log 1/h)^{1/2}}\le 1\qquad a.s.\label{8.35kk0},
   \end{equation}
  where $u^{*}$  is defined in   (\ref{ustar}). 
    
    \et

 	   \Proof    
	    We first prove (\ref{8.35qqq}).  Let $\th<1$ and consider Lemma \ref{lem-8.2a} applied to $Y_{\al}(t)=X_{\al}^{1/2}( \th^{n+1}+t)- X_{\al}^{1/2}(\th^{n+1})$,  with  $  S_{n}= \th^{n}-\th^{n+1}$. We have
     \bea
 &&   \sup_{t\in [0,\th^{n}-\th^{n+1 }]}|  Y_{\al}^{1/2}(t)- Y_{\al}^{1/2}(0)| \\
 & &\hspace{1 in}  =\sup_{t\in [0,\th^{n}-\th^{n+1 }]}|   X_{\al}^{1/2}( \th^{n+1}+t)- X_{\al}^{1/2}(\th^{n+1})|\nn\\
  & &\hspace{1 in} =\sup_{t\in [\th^{n+1},\th^{n} ]}|   X_{\al}^{1/2}(  t)- X_{\al}^{1/2}(\th^{n+1})|.\nn
   \eea
  Since it follows from  (\ref{2.58q}) that (\ref{8.28})	holds we see by Lemma \ref{lem-8.2a} that 
   \bea
 &&P\(   \sup_{t\in [\th^{n+1},\th^{n }]}|  X_{\al}^{1/2}(t)- X_{\al}^{1/2}(\th^{n+1})|> a  \vf (  S_{n})+\sum_{p=1}^{\ff }\ka (p)\vf (  S_{n}/ n(p))\)\nn\quad \\
 &&\qquad
   \le n^{2}F(a)+\sum_{p=1}^{\ff }n^{2}(p) F(\ka (p))  .\label{8.15ww}
   \eea
  If we take $a= ((3+\ep) \log n )^{1/2}$ and $\ka(p)=(3 \log n(p) )^{1/2}$, as in  the proof of Lemma \ref{theo-8.1}, we see that for all $n$ sufficiently large, the second  line of (\ref{8.15ww}) is a term of a converging sequence. Consequently, by the calculations in 
 (\ref{3.21q})--(\ref{3.35}),  the event   
  \be    \sup_ {t\in [ \th^{n+1},\th^{n }]}|  X_{\al}^{1/2}(t)-  X_{\al}^{1/2}( \th^{n+1 })|\label{9.6w} >  \(  a  \vf (  S_{n})+\frac{\sqrt3}{\log 2}  \int_{0}^{1/n }  \frac{\vf(S_{n} u)}{u(\log 1/u)^{1/2}}\,du \)   \ee 
 infinitely often,  has probability zero.  
 
 Note that  by the  condition that   $\vf(t)$  is regularly varying at zero with positive index 
 \begin{equation}
   \frac{\vf(S_{n}u)}{\vf(S_{n} )}\le C u^{\bb}\qquad \forall\, u\in[0,1/2]
   \end{equation}
 for some constant $C$ and all  $n$ sufficiently large. Therefore,  
  \begin{equation}
     \int_{0}^{1/n }  \frac{ \vf(S_{n}u)}{u(\log 1/u)^{1/2}}\,du   \le  C\vf(S_{n} )   \int_{0}^{1/2 }  \frac{ u^{\bb-1}}{ (\log 1/u)^{1/2}}\,du \le  C'\vf(S_{n} )  
     \end{equation}
 for some constant $C'$ and all  $n$ sufficiently large. This shows that the integral in (\ref{9.6w}) is $o((\log n)^{1/2}  \vf (  S_{n}))$ as $n\to\ff$.
 
  Using the regular variation hypothesis again we see that  
 \begin{equation}
  \vf (  S_{n})\le 2(\th^{-1}-1 )^{\bb} \vf ( \th^{n+1 }),  
   \end{equation}
   for some $\bb>0$ and all $n$ sufficiently large. 
 Therefore, for any $\ep>0$, by taking $1-\th$ sufficiently small,  the right-hand side of (\ref{9.6w})  
\begin{equation}
 <\ep \vf(\th^{n+1}) (\log n)^{1/2} \label{3.47s},
   \end{equation}
 for  all  $n$ sufficiently large.   It follows from this that for all $\ep >0$ we can choose $\th<1$ such that 
 the probability that  
 \be    \sup_ {t\in [ \th^{n+1},\th^{n }]}|  X_{\al}^{1/2}(t)-  X_{\al}^{1/2}( \th^{n+1})|\label{9.6ww0} >  2\ep \vf(\th^{n+1}) ( \log (n+1))^{1/2} \ee 
  infinitely often, is zero, for all $\ep>0$.
  
  By  Lemma \ref{lem-8.5a}, since  $\si (s,t)\leq \vf \(|t-s|\)$,       \begin{equation}
   \limsup_{ n\to\ff} \frac{|X_{\al}^{1/2}(\th^{n+1})-X_{\al}^{1/2}(0)|}{  \vf(\th^{n+1}) ( \log (n+1))^{1/2} }\le 1\qquad a.s.\label{8.37ez}
   \end{equation}	
   Combining this  with the statement in the sentence containing (\ref{9.6ww0}),    we get  
       \begin{equation}
   \limsup _{n\to\ff}   \frac{\sup_{t\in [\th^{n+1}, \th^{n}]}|X_{\al}^{1/2}(t )-X_{\al}^{1/2}(0)|}{\vf( \th^{n+1}  ) (\log \log \th^{-(n+1)} )^{1/2}}\le 1+2\ep\qquad a.s.\label{9.10w2} 
   \end{equation}
Since   $\vf( t)$ is  asymptotic to a monotonically increasing function near zero we  get,
      \begin{equation}
   \limsup _{n\to\ff}  \sup_{t\in [\th^{n+1}, \th^{n}]} \frac{|X_{\al}^{1/2}(t )-X_{\al}^{1/2}(0)|}{\vf( t  ) (\log \log 1/t )^{1/2}}\le 1+2\ep\qquad a.s.\label{9.10w3} 
   \end{equation}
 and since this holds for any $\ep>0$   we get (\ref{8.35qqq}).  
 
  By (\ref{9.6ww0})  and the fact   that  $\vf(\th^{n+1})\leq C u^{1/2}(\th^{n+1},\th^{n+1})$, 
 \be      \limsup_{ n\to\ff}  \sup_ {t\in [ \th^{n+1},\th^{n }]}  {|X_{\al}^{1/2}(t)-  X_{\al}^{1/2}( \th^{n+1})| \over (u^{ *}(\th^{n+1},\th^{n+1})  \log n)^{1/2}}\leq  2\ep C, \hspace{.2 in}    a.s. \label{9.6ww0j}\ee 
  for all $\ep>0$. By Lemma \ref{lem-8.5}  
       \begin{equation}
   \limsup_{ n\to\ff} \frac{ X_{\al}^{1/2}(\th^{n+1}) }{ u(\th^{n+1 } ,\th^{n+1 } )^{1/2}  (\log n)^{1/2}  }\le 1\qquad a.s.\label{8.37ezdj}
   \end{equation}
   Hence 
        \begin{equation}
   \limsup_{ n\to\ff} \sup_ {t\in [ \th^{n+1},\th^{n }]}\frac{ X_{\al}^{1/2}(t) }{ (u^{ *}(\th^{n+1},\th^{n+1})  \log n)^{1/2} }\le 1 +2\ep C\qquad a.s.\label{8.37ezdj2},
   \end{equation}  
  which gives   (\ref{8.35kk0}).    
 
 \medskip  The proof of  (\ref{8.3jj}) is more subtle. We use the following lemma.

  \begin{lemma}\label{lem-9.1} Let $X=\{X(t ), t\in \TT\}$  be an  $\al$-permanental process with kernel $u(s,t)$ and  sigma  function $\si (s,t)$ and such that  $u\(   t,t\)>0$ for all $t\in \TT$. 
 Let    
 \begin{equation}
   \wt X(t)=\frac{X (t)}{u(t,t)}.\label{3.48mm}
   \end{equation}
 Then  $\wt X=\{\wt X(t), t\in \TT\}$  is an  $\al$-permanental process with kernel 
\begin{equation}
\wt u(s,t)=\frac{u(s,t)}{(u(s,s)  u(t,t))^{1/2}}.\label{tri1}
\end{equation} 
Let $\wt \si (s,t)$ be the  sigma  function of   $(\wt X_{s}, \wt X_{t})$.
Then
\begin{equation}
   \wt \si^{2}(s,t)\le \frac{ \si^{2}(s,t)}{(u(s,s)  u(t,t))^{1/2}}\label{9.2}.
   \end{equation}
 \end{lemma}
 
 \Proof   To verify (\ref{tri1}) note that for any $ t_{1},\ldots,  t_{n}$, 
\begin{equation}
   E\(e^{-\sum_{i=1}^{n}\la_{i}X_{t_{i}}}\) 
 = \frac{1}{ |I+U\La|^{ \al}},   \label{int.1b} 
 \end{equation}
where $U$ is the $n\times n$ matrix  with entries $u(t_{i},t_{j})$ and  $\La$ is a diagonal matrix with entries $\la_{i}$, $1\le i\le n$. Let  $D$ be is the diagonal matrix with entries $u(t_{i},t_{i})$, $1\le i\le n$. Then
\bea
   E\(e^{-\sum_{i=1}^{n}\la_{i}\wt X_{t_{i}}}\)&= &
  E\(e^{-\sum_{i=1}^{n}\la_{i}X_{t_{i}}/u(t_{i},t_{i})}\)  \label{int.1bb} \\
&= & \frac{1}{ |I+U\La D^{-1}|^{ \al}}=\frac{1}{ |I+D^{-1/2}UD^{-1/2}\La |^{ \al}},  \nn
 \eea
which gives (\ref{tri1}).

 Note that $\wt u(s,s)=\wt u(t,t)=1$. By the inequality between geometric and arithmetic mean,    
 \bea
 \wt \si^{2}(s,t)&=&2 -\,  \frac{ 2(u(s,t)  u(t,s))^{1/2}} { (u(s,s)  u(t,t))^{1/2}  }\\
 &=&    \,  \frac{2 (u(s,s)  u(t,t))^{1/2}-  2(u(s,t)  u(t,s))^{1/2}} { (u(s,s)  u(t,t))^{1/2}  }   \nn\\
 &\le&    \frac{ u(s,s) + u(t,t)-  2 (u(s,t)  u(t,s))^{1/2}} {  (u(s,s)  u(t,t))^{1/2} },   \nn
 \eea
 which is (\ref{9.2}).\qed

 	\noindent{\bf Proof of  Theorem \ref{theo-8.1w} continued }  
  To prove    (\ref{8.3jj})  we consider $\wt X_{\al}(t)$ and $\wt\si(s,t)$ as defined in (\ref{3.48mm}) and (\ref{9.2})  for $t\in(0,1]$.  
Let $\th<1$. 
  Then by Lemma \ref{lem-9.1}, for all $\th^{n+1}\le s,t\le \th^{n}$,  
     \begin{equation}
   \wt \si^{2}(s,t)\le  \frac{  \vf^{2}(|t-s|) }{ u^{*}(\th^{n+1 }, \th^{n+1 } )  }   := \wt\vf_{n}^{2}(|t-s|)\label{9.2qw},
   \end{equation}
since $u^{*}(\th^{n+1 },\th^{n+1 }  )=\inf_{h\in [\th^{n+1},\th^{n}]}u^{*}(h,h)$.
 Since $\vf $ is regularly varying with index, say $\bb>0$  we see that for all    $\ep>0$, there exists an $n'_{0}$ such that for all $n\ge n'_{0}$   
 \begin{equation}
 \vf   (\th^{n }-\th^{n+1} )\le (1 +\ep)\( \frac{1-\th }\th \)^{\bb}  \vf ( \th^{n+1  })  \label{9.2qq}.
   \end{equation}
  Also since $\vf$   regularly varying it is asymptotic to an increasing  regularly   varying function with  index $\bb$. To simplify the proof let us simply take $\vf$ to be increasing. By hypothesis $  \vf^{ 2}(h)\le C u^{*}(h,h)$ for all $h$ sufficiently small. Therefore for all $h\in [\th^{n+1},\th^{n}]$, 
  \begin{equation}
   \vf^{2}(\th^{n+1})\le \vf^{ 2}(h)\le C u^{*}(h,h),
   \end{equation}
 which implies that $  \vf^{2}(\th^{n+1})\le C  u^{*}(\th^{n+1 } ,\th^{n+1} )$. Using (\ref{9.2qq}) we see that for all $ \ep'>0$ and $\th$ sufficiently close to 1, and $\th^{n+1}\le s,t\le \th^{n}$, 
 \begin{equation}
   \wt\vf_{n}^{2}(|t-s|)\le  \ep'
   \end{equation}
 It follows from this that the right--hand side of (\ref{9.6w})   applied to $\wt X_{\al}$  with dominating function $\wt \vf_{n}$, is less than $\ep'   (\log n)^{1/2}$  for all $n$ sufficiently large.
  Consequently,   
 the probability that  
 \be    \sup_ {t\in [ \th^{n+1},\th^{n }]}| \wt  X_{\al}^{1/2}(t)-  \wt X_{\al}^{1/2}( \th^{n+1})|\label{9.6ww} >   \ep'   ( \log n)^{1/2} \ee 
  infinitely often, is zero, for all $ \ep'>0$.
  
    It   follows from  Lemma \ref{lem-8.5}  that  
     \begin{equation}
   \limsup_{ n\to\ff} \frac{ \wt X_{\al}^{1/2}(\th^{n+1}) }{   (\log n)^{1/2}  }\le 1\qquad a.s.\label{8.37ezd}
   \end{equation}	
    Combining (\ref{9.6ww}) and (\ref{8.37ezd}) we get  
     \be  
     \sup_ {t\in [ \th^{n+1},\th^{n }]}  \wt   X_{\al}^{1/2}(t) \label{9.6wwpj} \le  ( 1+  \ep')  ( \log n)^{1/2}, 
     \ee
   almost surely as $n\to\ff$.  
   	This gives   
       \begin{equation}
   \limsup_{ t\to 0} \frac{ X_{\al}^{1/2}  (t) }{(u(t,t)\log\log 1/t))^{1/2} }\le  1  \qquad a.s.\label{8.35x}
   \end{equation}  
   which   is    (\ref{8.3jj}). \qed

\medskip	    The next lemma gives relationships between the sigma function of a kernel and its majoring function $\vf.$  
   
 \begin{lemma}\label{lem-3.5qq} Let $X$ be an $\al$-permanental process   with  kernel $u(s,t)$ and suppose that $u(0,0)=0$. Then 
\begin{equation}
   \si^{2}(0,t)=    \si^{2}(t,0)=u(t,t).\label{3.47}
   \end{equation}
If, in addition,  $\si^{2}(s,t)$ is a function of $|t-s|$  then, necessarily, 
\begin{equation}
   \si^{2}(s,t) =  u(|t-s|,|t-s|) .\label{3.49w}
   \end{equation}
    In general, if $\vf$ satisfies (\ref{8.14d})
  \begin{equation}
   \vf^{2}(|t-s|)\ge (u^{1/2}(t,t)-u^{1/2}(s,s))^{2}.\label{3.49q}
   \end{equation} 
   Therefore, if $\lim_{t\to\ff}u (t,t)=\ff$, necessarily  
   \begin{equation}
 \lim_{t\to\ff}  \frac{\vf^{2}(t)}{ u (t,t)}\geq 1.
   \end{equation}

   \el
   
    \Proof    Note that  for
 any $\al$-permanental process $X$  
 \begin{equation}
   E[(X(s)-E(X(s))(X(t)-E(X(t))]=\al u(s,t)u(t,s),\label{3.50ak}
   \end{equation}
   \cite[p. 135]{VJ}. 
Clearly, this is equal to 0 when $s=0$   since $u(0,0)=0$ implies that $X(0)=0$. This     shows that  $u(t,0)u (0,t)=0$.  This and $u(0,0)=0$ gives (\ref{3.47}).

  If, in addition,  $\si^{2}(s,t)$ is a function of $|t-s|$,   then by (\ref{3.47}), when $t>s$, 
\begin{equation}
   \si^{2}(s,t)=\si^{2}(0,t-s)= u(t-s,t-s),
   \end{equation}
which gives (\ref{3.49w}).  

  The statement in (\ref{3.49q}) follows because 
 \be 
 u(s,t)u(t,s)\le u(s,s)u(t,t),\label{3.50xx}
 \ee
which, itself, follows from (\ref{3.50ak})
and the Schwartz Inequality.\qed

  The next theorem is an analogue of Theorem \ref{theo-8.1w} for the uniform modulua of continuity.

  \bt \label{theo-8.1wm}  	   Let  $X_{\al}=\{X_{\al}(t ), t\in [0,1]\}$  be an $\al$-permanental process with  kernel $u(s,t)$ and sigma   function $\si (s,t)$ for which   (\ref{8.14d})  and (\ref{3.21}) hold.  Assume furthermore, that  
	   $\vf(t)$  is regularly varying at zero with positive index.        Then  
         \begin{equation}
   \limsup_{ h\to 0}\sup_{\stackrel{|s-t|\le h}{s,t\in [0,1]}}\frac{|X^{1/2}_{\al} (s)-X^{1/2}_{\al} (t)|}{    \vf(h)( \log  1/h)^{1/2}}\le  1\qquad a.s.\label{8.35rwv}
   \end{equation}  
   \et

     \Proof      Theorem \ref{theo-3.2} asserts that $X_{\al}$ has subgaussian increments.   Then (\ref{8.35rwv}) follows from \cite[Theorem 4]{M}. In this reference the theorem is proved for Gaussian processes with stationary increments but it only uses the estimate in (\ref{1.9d}), along with the fact that $\vf$ is greater than the $L^{2}$ estimate for the increments of the Gaussian process. The difference of the factor 2 is explained by the observation in the line prior to (\ref{1.9f}).  
     
     (Of course we don't have  (\ref{1.9d}) for $X_{\al}$ but only (\ref{1.9fa}). It is easy to see that the additional powers of $\la$ do not  affect  the estimates used in  \cite[Theorem 4]{M}.)

     There is another item in the proof of  \cite[Theorem 4]{M} that needs explanation. In \cite[Theorem 4]{M} which deals with the Gaussian process $\{G(t),t\in [0,1]\}$, with stationary increments
     \begin{equation}
   \si(h)=\(E\(G(t+h)-G(t)\)^{2}\)^{1/2}\label{inc}
   \end{equation}
  is written as $\si(h)=\exp(-g(\log 1/h))$ and it is required that $1/g'(\log 1/h)=o(\log 1/h)$.  Note that when $\si(h)=Ch^{\al}$, for some constant $C$,
    $1/g'(\log 1/h)=1/\al$, much weaker than what is allowed. However it isn't necessary to require that $\si$ is differentiable. All the estimates in the proof of \cite[Theorem 4]{M} that use the condition $1/g'(\log 1/h)=o(\log 1/h)$ follow easily form the condition that  $\vf $  is regularly varying at zero with positive index.    For example when   $\vf $  is regularly varying at zero with index $\bb$, instead of  \cite[(2.16)]{M}, we have
    \begin{equation}
   \frac{1}{\vf(t_{k})}\int_{0}^{t_{k}}\frac{\vf(u)}{u}\,du=\int_{0}^{1}   \frac{\vf(t_{k}u)}{\vf(t_{k})}\frac{\,du}{u}\le (1+\ep)\int_{0}^{1}   \frac{u^{\bb}}{u} \,du =\frac{1+\ep}{\bb},
   \end{equation}
 as $t_{k}\to 0$,     for all $\ep>0$. 
 
 It is also required that $\si(h)$ is concave, but wherever this is used it is easy to get the same estimates when $\vf $  is regularly varying at zero with index $\bb\le 1$, which is always the case.   
  \qed

 	 Our interest in $\{X_{\al}^{1/2}(t)-X_{\al}^{1/2}(s); s,t\in R^{+}\}$ is primarily to use the results obtained to study the behavior of $\{X_{\al} (t)-X_{\al} (s); s,t\in R^{+}\}$. The next lemma does this.
  
 \begin{lemma} 
\label{theo-3.0}  Assume  that 
\begin{equation}
   \limsup_{ h\to 0}\sup_{0\le t\le h}\frac{|X^{1/2}(t)-X^{1/2}(0)|}{\om(h)}\le 1\qquad a.s.\label{8.35a}
   \end{equation}
   for some function $\om$ with $\lim_{h\to 0}\om(h)=0$.
 Then  
   \begin{equation}
   \limsup_{ h\to 0}\sup_{0\le t\le h}\frac{|X(t)-X(0)|}{\om(h)}\le 2 X^{1/2}(0)\qquad a.s.\label{8.35q}
   \end{equation}
   and if $X (0)=0$, 
     \begin{equation}
   \limsup_{ h\to 0}\sup_{0\le t\le h}\frac{ X(t) }{\om^{2}(h)}\le 1\qquad a.s.\label{8.35t}
   \end{equation}
Furthermore, if 
  \begin{equation}
   \limsup_{ h\to 0}\sup_{\stackrel{|s-t|\le h}{s,t\in [0,1]}}\frac{|X^{1/2} _{\al} (s)-X^{1/2} _{\al} (t)|}{  \rho(h)}\le  1\qquad a.s.\label{3.72xx},
   \end{equation}
 then
 \begin{equation}
   \limsup_{ h\to 0}\sup_{\stackrel{|s-t|\le h}{s,t\in [0,1]}}\frac{|X _{\al} (s)-X_{\al} (t)|}{  \rho(h)} \le 2 \sup_{t\in[0,1]}X^{1/2}(t)\qquad a.s.\label{3.73xx},
   \end{equation}
 
 \end{lemma}

 \Proof  The statement in (\ref{8.35t}) is trivial. For (\ref{8.35q}) we note that \begin{equation}
   |X(t  )-X(0)|\le |X^{1/2}(t  )-X^{1/2}(0)|\(|X^{1/2}(t  )-X^{1/2}(0)|+2X^{1/2}(0)\).
   \end{equation}
 Therefore, the left hand side of (\ref{8.35q}) is bounded by
 \begin{equation}
    \limsup_{ h\to 0}\sup_{0\le t\le h}\frac{|X^{1/2}(t)-X^{1/2}(0)|}{\om(h)}\(|X^{1/2}(t  )-X^{1/2}(0)|+2X^{1/2}(0)\).
   \end{equation}
 However, (\ref{8.35a}) implies that $X(t)$ is continuous at $t=0$, almost surely.   Therefore
 \begin{equation}
       \limsup_{ h\to 0}\sup_{0\le t\le h} |X^{1/2}(t)-X^{1/2}(0)|=0 \qquad a.s.   \end{equation}
 Using this and (\ref{8.35a}) we get (\ref{8.35q}).
 
 The result for the uniform modulus follows similarly since
 \begin{equation}
   |X(t  )-X(s)|\le |X^{1/2}(t  )-X^{1/2}(s)| |X^{1/2}(t  )+X^{1/2}(s)|.
   \end{equation}
 \qed 
 
   \noindent {\bf Proofs of  Theorems  \ref{theo-1.4} and \ref{theo-1.4unif}  } The proofs   follow immediately from Theorems \ref{theo-8.1w} and  \ref{theo-8.1wm}  and Lemma \ref{theo-3.0}.

 \section{   Upper bounds for the local moduli of continuity and rate of growth of  permanental processes, II}\label{newest}

 The next corollary exhibits  different upper bounds for the local modulus of continuity that are larger than the ones that hold under the hypotheses of Theorem \ref{theo-1.4} . 
  
 \begin{corollary}\label{cor-3.1}   Let  $X=\{X(t ), t\in [0,1]\}$  be an   $\al$-permanental process with   kernel $u(s,t)$   and with  sigma function $\si(s,t)$ for which   (\ref{8.14d})  and (\ref{3.21}) hold and for which   (\ref{8.33}) also holds  uniformly in  $V\leq V_{0}$ for some $V_{0}<\ff$. Then  
\begin{equation}
   \limsup_{ h\to 0}\sup_{0\le t\le h}\frac{|X(t)-X(0)|}{\Theta(h) }\le \frac{4\sqrt2}{\log 2} X^{1/2}(0)\qquad a.s.\label{8.35qq}
   \end{equation}
 where 
 \begin{equation}
   \Theta(h)=\int_{0}^{h^{2} }  \frac{\vf( u) }{u(\log 1/u )^{1/2}}\,du+  \vf( h)  (\log 1/h )^{1/2}.
   \end{equation}
 \end{corollary}
 The proof follows from  Lemma \ref{theo-8.1}  applied to $Y(t)=X^{1/2}(t)$,  Lemma \ref{theo-3.0}   and  the next lemma.

 \begin{lemma}\label{lem-3.7}  
 \be  \int_{0}^{1/2}  \frac{\vf(hu)}{u(\log 1/u)^{1/2}}\,du\le 2\Theta(h)+o(\Theta(h)).\label{3.64}
\ee
 \end{lemma} 
 
\Proof  We have 
 \begin{equation}
   \int_{0}^{1/2}  \frac{\vf(hu)}{u(\log 1/u)^{1/2}}\,du =   \int_{0}^{h }  \frac{\vf( u/2)} {u(\log 1/u-\log 1/(2h))^{1/2}}\,du      \end{equation}
   Note that
   \be 
 \int_{0}^{h^{2} }  \frac{\vf( u/2)} {u(\log 1/u-\log 1/(2h))^{1/2}}\,du \le  \sqrt 2\int_{0}^{h^{2} }  \frac{\vf( u/2)} {u(\log 1/u )^{1/2}}\,du,  
   \ee 
since $\log 1/(2h)\le (1/2) \log 1/u $ for $u\in [0,h^{2}]$. Furthermore,
\be  
      \int_{h^{2}}^{h }  \frac{\vf( u/2)} {u(\log 1/u-\log 1/(2h))^{1/2}}\,du    \le  \vf( h) \int_{h^{2}}^{h }  \frac{1} {u(\log 1/u-\log 1/(2h))^{1/2}}\,du.  \label{3.50}
 \ee 
The integrand of the last integral above is the derivative of $-2 (\log 1/u	-\log 1/(2h))^{1/2}$. Consequently,   \begin{equation}
   \int_{h^{2}}^{h }  \frac{1} {u(\log 1/u-\log 1/(2h))^{1/2}}\,du = 2(\log 1/h)^{1/2}+o((\log 1/h)^{1/2})\quad\mbox{as }h\to0.\label{3.51q}
   \end{equation}
  Combining these relationships we get (\ref{3.64}).\qed

   It is interesting to note that for certain functions $\vf$, except for a multiplicative factor, we can reverse the inequality in  (\ref{3.64}).
 
  \begin{lemma}   Suppose that
  \begin{equation}
   \limsup_{h\to 0}\frac{\vf(h^{2} )}{\vf(h)}= C.\label{3.67mm}
   \end{equation}
   Then    
 \be  \int_{0}^{1/2}  \frac{\vf(hu)}{u(\log 1/u)^{1/2}}\,du\ge  C \Theta(h)+o(\Theta(h)).\label{3.64q}
\ee
 \end{lemma} 
 
\Proof    We have
\bea
    \int_{0}^{1/2}  \frac{\vf(hu)}{u(\log 1/u)^{1/2}}\,du&\ge&    \int_{0}^{h/2}  \frac{\vf( u)}{u(\log 1/u)^{1/2}}\,du\\
    &\ge&    \int_{0}^{h^{2}}  \frac{\vf( u)}{u(\log 1/u)^{1/2}}\,du+ \int_{h^{2}}^{h/2}  \frac{\vf( u)}{u(\log 1/u)^{1/2}}\,du\nn
   \eea
and
\bea
   \int_{h^{2}}^{h/2}  \frac{\vf( u)}{u(\log 1/u)^{1/2}}\,du&\ge& \vf(h^{2})   \int_{h^{2}}^{h/2}  \frac{1}{u(\log 1/u)^{1/2}}\,du\\
   &=&\vf(h^{2}) ((\log 1/h)^{1/2}+o(h))\nn\\
     &\sim&C\vf(h ) ((\log 1/h)^{1/2}+o(h))\nn.
   \eea
Consequently, we get (\ref{3.64q}).\qed

  \begin{example} {\rm Let $\ga >1/2$. Then if  
    \begin{equation}
     \vf(h) =\frac{1}{(\log 1/h)^{\ga }},\quad \Theta(h)=\frac{2^{ (3/2)-\ga}+2\ga-1}{2\ga -1}\vf( h)  (\log 1/h )^{1/2}.
   \end{equation}
 Let $\bb> -1$. Then if
  \begin{equation}
     \vf(h) =\frac{(\log\log 1/h)^{\bb}}{(\log 1/h)^{1/2}},\quad \!\Theta(h)= \(\frac{1}{1+\bb}+o(h)\)\vf( h)  ( \log 1/h )^{1/2} \log\log 1/h.
   \end{equation}
  Note that in both these cases $C=1$ in (\ref{3.67mm}).

 }\end{example}

  For more examples see \cite[Lemma 7.6.5 and Example 7.6.6]{book}.

   \medskip		We now examine the relationship between $\si(s,t)$ and the $L_{2}$ metric for Gaussian processes.    Let  
\begin{equation}
   \rho ^{2}(s,t):=   u(s,s)+u(t,t)-(u(s,t)+u(t,s)) .
   \end{equation}
   Although we don't require that $u(s,t)$ is symmetric, when $u(s,t)$ is symmetric $\rho (s,t)=\si(s,t)$. In general we get the next lemma.
   
    \begin{lemma} \label{lem-3.8} Let  $\si(s,t)$  be as defined in (\ref{1.2}). Then    

\begin{equation}
    \si^{2} (s,t) = \rho^{2}(s,t)  + (u^{1/2}(s,t)-u^{1/2}(t,s) )^{2}.\label{3.74a}
   \end{equation}
    In addition, when 
$ u (s,t)\vee u(t,s)\le u (s,s)\wedge u (t,t)$ for all $s,t\in T$,
\be
  \si (s,t)\le \sqrt 2\rho (s,t) \label{3.75a}.
 \end{equation}
 \end{lemma}

 \Proof We have 
 \bea    \si^{2}(s,t)&=&\rho^{2}(s,t)  + u(s,t)+u(t,s)-2(u(s,t)u(t,s))^{1/2} \label{3.76a}\\
 &=&\rho^{2}(s,t)  + (u^{1/2}(s,t)-u^{1/2}(t,s) )^{2}.\nn
   \eea 
 The inequality in (\ref{3.75a}) is given in \cite[Lemma 5.5]{MRsuf}. 
\qed

  \begin{corollary} \label{cor-4.2}

 If 
\begin{equation}
   u(s,t)=v(s,t) +h(t)\label{canon}
   \end{equation}
 and $v$  is symmetric  \begin{equation}
    \si^{2} (s,t)\le v(s,s)+v(t,t)-2v(s,t)   +|h (s )-h (t ) |.\label{3.74as}
   \end{equation}
  If, in addition $\inf _{s,t \in I}u(s,t)\ge \de$, for some interval $I$, then for all $s,t\in I$,
  \begin{equation}
    \si^{2} (s,t)\le v(s,s)+v(t,t)-2v(s,t)   +\frac{|h (s )-h (t ) |^{2}}{4\de } .\label{3.74ff}
   \end{equation}
   
 \end{corollary}

\Proof  The inequality in (\ref{3.74as}) follows immediately from (\ref{3.74a}).  To obtain (\ref{3.74ff}) note that for $a<b$
\begin{equation}
   b^{1/2}-a^{1/2}=\int_{a}^{b}\frac{1}{2u^{1/2}}\,du\le \frac{b-a}{2a^{1/2}}.\label{dog}
   \end{equation}
Consider $u(s,t)$ in (\ref{canon}) and suppose that $h(t)>h(s)$.  Then by (\ref{dog}),  for $s,t\in I$ 
\begin{equation}
   u^{1/2}(s,t)-u^{1/2}(t,s) \le\frac{u(s,t)-u(t,s)}{2\de^{1/2}} \le \frac{h(t)-h(s)}{2\de^{1/2}}.\label{4.20}
   \end{equation}
  Using this and (\ref{3.74a}) and the fact that $v$ is symmetric, we get (\ref{3.74ff}).\qed

 	  \begin{remark}\label{rem-4.1} {\rm The inequality in   (\ref{3.74ff}) may be  smaller than the one in (\ref{3.74as}) even when $u(s,t)$ has the form of (\ref{canon}). For example, suppose
  \begin{equation}
   u(s,t)=e^{-\la|t-s|} +e^{-r t}\label{can}\qquad r,\la>0.
   \end{equation}
  Obviously 
 \begin{equation}
   |e^{-r t}-e^{-r s}|\sim r |t-s|\qquad \mbox{ as $s,t\to 0$}.
   \end{equation}
 Consequently, it follows from (\ref{3.74ff}) that
 $
   \si^{2}(s,t)$ is bounded by $\la |t-s|$ as $s,t\to 0$,
   whereas (\ref{3.74as})  only gives   that  is bounded by $(\la+r) |t-s|$ as $s,t\to 0$.
    }\end{remark}

\noindent{\bf Proof of  upper bounds in Theorem \ref{theo-1.8new} }    We show in \cite[Section 5]{KMR} that 
 \begin{equation}
   u_{T_{0};\ga,\bb}(x,y)=R(x,y)_{\ga,\bb}+H_{\ga,\bb}(x,y)  \label{1.22w}
   \end{equation}
where $R_{\ga,\bb}$ is symmetric and $H_{\ga,\bb}$ is antisymmetric. Explicitly,
 \begin{equation}
 R_{\ga,\bb}(x,y)=
C_{\ga,\bb}  \(|x|^{\ga}+|y|^{\ga}-
|x-y|^{\ga}\),     \label{9.6}
\ee
       \be {H_{\ga,\bb}(x,y) }\label{9.6a}  
  =\bb C_{\ga,\bb}  \( { sign }(x)|x|^{\ga }- { sign }( y)|y|^{\ga }-
 { sign }(x-y)|x-y|^{\ga}\), 
\ee
and
\begin{equation}
 C_{\ga,\bb}={-\sin \left( (\ga+1)\frac{\pi}{2} \right) 
\Gamma( -\ga) \over \pi (1+\bb^{2}  \,\tan^{2} ( (\ga+1)\pi/2))}>0.
 \end{equation}
By Lemma \ref{lem-3.8} the  sigma function for $Y_{\al:\ga,\bb}$, which we denote by $\si_{T_{0}:\ga\bb}$, satisfies
\bea
     \si^{2}_{T_{0};\ga,\bb}(x,y)&\le&  u_{T_{0};\ga,\bb}(x,x)+   u_{T_{0};\ga,\bb}(y,y)-(   u_{T_{0};\ga,\bb}(x,y)  +  u_{T_{0};\ga,\bb}(x,y))\nn\\
     &&\qquad+|   u_{T_{0};\ga,\bb}(x,y)  -  u_{T_{0};\ga,\bb}(y,x)|\label{4.24}\\
     &=&R_{\ga,\bb}(x,x)+R_{\ga,\bb}(y,y)-2R_{\ga,\bb}(x,y)+2|H_{\ga,\bb}(x,y)|,\nn
   \eea
  where, for the last equation we use the facts that $R_{\ga,\bb}$ is symmetric, $H_{\ga,\bb}$ is antisymmetric and $H_{\ga,\bb}(x,x)\equiv 0$.
      It is easy to see that 
   \begin{equation}
   |H(x,y)|\le |\bb|C_{\ga,\bb}|x-y|^{\ga}.\label{3.85mm}
   \end{equation}
Using this and (\ref{9.6}), we get 
 \begin{equation}
   \si_{T_{0};\ga,\bb}(x,y)\le (2  (1+|\bb|) C_{\ga,\bb})^{1/2}|x-y|^{\ga/2}:=\vf(|x-y|).\label{3.87}
   \ee  
  Since $ \vf$    is regularly varying at zero with positive index  we see that the upper bounds in Theorem \ref{theo-1.8new} follows from Theorems \ref{theo-1.4},  \ref{theo-1.4unif} and \ref{theo-3.2q}.
\qed

 We refer to the permanental process  with kernel  $u_{T_{0},\ga,\bb}$ as   $\mbox{FBMQ}^{\ga ,\bb}$.
We use this notation because when $\bb=0$, $u_{T_{0},\ga,\bb}$ is the covariance of   fractional Brownian motion of index   $\ga $, (i.e. FBM). We add the $Q$, for quadratic, to denote the square of this process, as one does in the designation of the squared Bessel processes, (BESQ).  \label{page 1}

 	\medskip	\noindent{\bf Proof of    Theorem \ref{theo-1.8exp}   }  
It follows from \cite[Lemma 5.2]{KMR} that 
\bea
   u_{\rho;\ga,\bb}(x,y)&=&\frac{1}{\pi}\int_{0}^{\ff}\frac{\cos\la(x-y)\RR e(\rho+\psi_{\ga,\bb}(\la))}{|\rho+\psi_{\ga,\bb}(\la)|^{2}}\,d\la\\
   &&\qquad +\frac{1}{\pi}\int_{0}^{\ff}\frac{\sin\la(x-y)\II m \,  \psi_{\ga,\bb}(\la)  }{|\rho+\psi_{\ga,\bb}(\la)|^{2}}\,d\la\nn.
   \eea
Consequently, as in (\ref{4.24})
\bea
   \si^{2}_{\rho;\ga,\bb}(x,y)&\le&\frac{2}{\pi}\int_{0}^{\ff}\frac{(1-\cos\la(x-y))\,\RR e(\rho+\psi_{\ga,\bb}(\la))}{|\rho+\psi_{\ga,\bb}(\la)|^{2}}\,d\la\label{4.28q} \\
   &&\qquad +\frac{2}{\pi}\left |\int_{0}^{\ff}\frac{\sin\la(x-y)\II m \, \psi_{\ga,\bb}(\la) }{|\rho+\psi_{\ga,\bb}(\la)|^{2}}\,d\la\right |\nn.
   \eea
We write
\bea
 &&  \int_{0}^{\ff}\frac{ (1-\cos\la(x-y))\,\RR e(\rho+\psi_{\ga,\bb}(\la))}{|\rho+\psi_{\ga,\bb}(\la)|^{2}}\,d\la\\
 & &\qquad\le \int_{0}^{\ff}\frac{ (1-\cos\la(x-y))\,\RR e\, \psi_{\ga,\bb}(\la) }{| \psi_{\ga,\bb}(\la)|^{2}}\,d\la+\nn  2\rho \int_{0}^{\ff}\frac{\sin^{2} \frac{\la(x-y)}{ 2}\, }{|\rho+\psi_{\ga,\bb}(\la)|^{2}}\,d\la  .
   \eea
It follows from \cite[(5.39)]{KMR} that the first integral to the right of the inequality sign is equal to $\pi C_{\ga,\bb}|x-y|^{\ga}$. The second integral to the right of the inequality sign is bounded by
\begin{equation}
 \frac{2}{\rho} \int_{0}^{1} \sin^{2} \frac{\la\th }{ 2} \,d\la +2\rho \int_{1}^{1/\th}\frac{\sin^{2} \frac{\la\th}{ 2}\, }{| \psi_{\ga,\bb}(\la)|^{2}}\,d\la+2\rho \int_{1/\th}^{\ff}\frac{1}{| \psi_{\ga,\bb}(\la)|^{2}}\,d\la,
   \end{equation}
where $\th=|x-y|$. The first of these integrals is $O(\th^{2})$ as $\th\to 0$.
The third is $O(\th^{2\ga+1})$ as $\th\to 0$. The second is $O(\th^{2})$ if $\ga>1/2$, $O(\th^{2})\log 1/\th$ if $\ga=1/2$ and $O(\th^{1+2\ga})\log 1/\th$ if $\ga<1/2$, all as $\th\to 0$. Thus we see that the first integral in (\ref{4.28q}) is bounded by $2 C_{\ga,\bb}|x-y|^{\ga}$.

 	We now consider the last integral  in (\ref{4.28q})
\bea
 && \int_{0}^{\ff}\frac{\sin\la(x-y)\II m \, \psi_{\ga,\bb}(\la) }{|\rho+\psi_{\ga,\bb}(\la)|^{2}}\,d\la\label{4.31} \\
 &&\qquad =   \int_{0}^{\ff}\frac{\sin\la(x-y)\II m \, \psi_{\ga,\bb}(\la) }{| \psi_{\ga,\bb}(\la)|^{2}}\,d\la+\int_{0}^{\ff}\frac{\sin\la(x-y)\II m \, \psi_{\ga,\bb}(\la) }{| \Psi_{\ga,\bb}(\la)|^{2}}\,d\la\nn
 \eea
where 
\bea
   \frac{1}{| \Psi_{\ga,\bb}(\la)|^{2}}&=&\frac{|\rho+\psi_{\ga,\bb}(\la)|^{2} -|\psi_{\ga,\bb}(\la)|^{2}}{|\psi_{\ga,\bb}(\la)|^{2}\,|\rho+\psi_{\ga,\bb}(\la)|^{2}}\label{4.32}\\
 &=&\frac{(\rho+\RR e\psi_{\ga,\bb}(\la))^{2} - ( \RR e\psi_{\ga,\bb}(\la))^{2}}{|\psi_{\ga,\bb}(\la)|^{2}\,|\rho+\psi_{\ga,\bb}(\la)|^{2}}\nn  \\
 &\le &\frac{2\rho }{|\psi_{\ga,\bb}(\la)| \,|\rho+\psi_{\ga,\bb}(\la)|^{2}}.\nn  
   \eea
It follows from \cite[(5.40)]{KMR} that the first integral to the right of the  equal  sign in (\ref{4.31}) is equal to $\pi\bb \,\mbox{sign}\,( x-y )  C_{\ga,\bb}|x-y|^{\ga}$. We now show that the second integral is little o of this. Using (\ref{4.32}) we see that the last integral in (\ref{4.31}) is bounded by $2\rho$ times 
\begin{equation}
   \int_{0}^{\ff}\frac{|\sin\la\th|   }{|\rho+ \Psi_{\ga,\bb}(\la)|^{2}}\,d\la\le \frac{\th}{2\rho^{2}}+ \int_{1}^{1/\th}\frac{ \la\th   }{|  \Psi_{\ga,\bb}(\la)|^{2}}\,d\la + \int_{1/\th}^{\ff}\frac{ 1  }{|  \Psi_{\ga,\bb}(\la)|^{2}}\,d\la\label{4.33}
   \end{equation}
It is easy to see that the first integral to the right of the  equal  sign in (\ref{4.33}) is  $O(\th )$ as $\th\to 0$ and 
the second  is   $O(\th ^{1+2\ga})$ as $\th\to 0$. Therefore,  the absolute value of the second  integral in (\ref{4.28q}) is bounded by $2|\bb |C_{\ga,\bb}|x-y|^{\ga}$.

Using the bounds for the last two integrals in (\ref{4.28q}) we see that 
\begin{equation}
   \si_{\rho;\ga,\bb}(x,y)\le (2  (1+|\bb|) C_{\ga,\bb})^{1/2}|x-y|^{\ga/2}:=\vf(|x-y|),\label{3.87ss}
   \ee  
   the same as in (\ref{3.87}).
  Since $ \vf$    is regularly varying at zero with positive index  we see that the upper bounds in Theorem \ref{theo-1.8new} follows from Theorems \ref{theo-1.4},  \ref{theo-1.4unif} and \ref{theo-3.2q}. \qed

\noindent{\bf Proof of   Theorem \ref{theo-1.9mm} }  
The $\al$-permanental process $\ov X_{\al, f} $ has kernel $\wh u_{f}  (s,t)$. We see from (\ref{3.74ff}) in  Corollary  \ref{cor-4.2} that   for $s,t\in [0,\de/\la]$, the sigma function of $\ov X_{\al,f} $ satisfies 
\be
  \si^{2}_{f}(s,t)\le \(2(1- e^{-\la|s-t|}\) +\frac{|f(t)-f(s)|^{2} }{4(1-\de)}  
  \ee 
  for all $\de$ sufficiently small.   Using this and   the fact that $f\in C^{ 2}$  implies that for $t>s$, $|f(t)-f(s)|\le f'(t)|t-s|$, we see that
\begin{equation}
   \si^{2}_{f}(s,t)\sim 2\la| t-s|,\qquad\mbox{as $s,t\to 0$}.
   \end{equation}
   Therefore, Theorem \ref{theo-1.9mm}
 follows from  Theorems \ref{theo-1.4} and \ref{theo-1.4unif}.  \qed

\section{Rate of growth of permanental processes at infinity}	\label{sec-4}

\medskip	   We begin by considering another  important class of processes for which we can  decrease	the upper bounds that can be obtained by  (\ref{8.37}). First we need some preliminary results.
 
 \medskip	 Let $X=\{X(t ), t\in R^{+}\}$  be an $\al$-permanental  process with kernel $u(s,t)$ and sigma function $\si(s,t)$. Then  by 
 \cite[Lemma 3.1]{MRnec}, for $\la>2(\al-1)\vee 0$,
   \begin{equation}
   P\(X(t)\ge u(t,t)\la\)\le \frac{2\la^{\al-1}e^{-\la}}{\Ga(\al)} \label{8.5q},
   \end{equation}
   and for $\la\ge 2$
   \begin{equation}
  \frac{2\la^{\al-1}e^{-\la}}{3\Ga(\al)} \le    P\(X(t)\ge u(t,t)\la\)   \label{8.5qq}.
   \end{equation}
    (It is interesting to note that   since $E(X(t))=\al u(t,t)$, if we were to consider  the rate of growth of $ X(t)/E(X(t))$, it would depend on $\al$. The results we give for  the rate of growth of   $ X(t)/u(t,t)$ do  not depend on $\al$.) 
 
 \medskip	The next observation is elementary.
 
 \bl\label{lem-8.5} Let $X=\{X(t_{n}), n\in \mathbf{N}\}$  be an $\al$-permanental sequence with kernel   $u(t_{i},t_{j})$. Then
     \begin{equation}
   \limsup_{ n\to\ff} \frac{X(t_{n})}{  u(t_{n},t_{n})\log n  }\le 1\qquad a.s.,\label{8.37e}
   \end{equation}	      
 	    or, equivalently,    
  
         \begin{equation}
   \limsup_{ n\to\ff} \frac{X^{1/2}(t_{n})}{  (u(t_{n},t_{n})\log n )^{1/2} }\le 1\qquad a.s.\label{8.37ee}
   \end{equation}	   
   \el

   \Proof  The  statement in (\ref{8.37e}) follows from (\ref{8.5q}) and the Borel-Cantelli   \qed	
	
 \medskip	\noindent	 {\bf Proof of Theorem \ref{theo-3.1w} }
 The proof follows from   Lemma \ref{lem-8.2a} similarly to how the results about the local modulus of continuity are obtained in  Section \ref{sec-8}. \njc We can assume that $u^{\ast}(T_{0},T_{0})>0 $ for some $T_{0}$. {\bf I don't see why. Delete to star based on   comment prior to Theorem 1.1.  }
Redefine  
   \begin{equation}
  F (a)= \sup_{s, t \geq T_{0} }P\(\frac{X^{1/2}(s)- X^{1/2}(t)}{\si (s,t)}\geq a\).\label{8.15zx}  
     \end{equation} 
   It follows from   (\ref{2.58q})  
    that (\ref{8.28}) holds. Let $\de>0$.   
 Then, by 
 Lemma \ref{lem-8.2a},  for all  $n$ with $n\de \geq T_{0}$,
   \bea
 &&P\(   \sup_{t\in [n\de,(n+1)\de]}|  X^{1/2}(t)- X^{1/2}(n\de)|> a  \vf (  \de)+\sum_{p=1}^{\ff }\th (p)\vf (  \de/ n(p))\)\nn\quad \\
 &&\qquad
   \le n^{2}F(a)+\sum_{p=1}^{\ff }n^{2}(p) F(\th(p))  .\label{8.15wee}
   \eea
Taking  $a= ((3+\ep) \log n )^{1/2}$ and $\th(p)=(3 \log n(p) )^{1/2}$, as in (\ref{8.39}), we see that for all $n$ sufficiently large, the second  line of (\ref{8.15wee}) is a term of a converging sequence. Consequently, as in the beginning of the proof of Lemma \ref{theo-8.1}, the event  
 \be   
  \sup_{t\in [n\de,(n+1)\de]}|  X^{1/2}(t)-  X^{1/2}( n\de)|\label{9.6wq} >  \(  a  \vf (  \de)+\frac{\sqrt3}{\log 2}  \int_{0}^{1/n }  \frac{\vf (\de u)}{u(\log 1/u)^{1/2}}\,du  \)  
   \ee 
 infinitely often, is zero.

 For any $\ep>0$ we can find a $\de>0$ so that 
 \begin{equation}
    a  \vf (  \de)  \le \frac{\ep} 2 ( \log n)^{1/2}
   \end{equation}
and   
  \begin{equation}
 \int_{0}^{1/n }  \frac{\vf (\de u)}{u(\log 1/u)^{1/2}}\,du     
     \end{equation}
  is bounded uniformly in $\de\leq 1$. 
 Therefore, for any $\ep>0$ we can find a $\de$ so that   the right hand side of (\ref{9.6wq}) is $\le\ep (  \log n)^{1/2}$ for all $n$ sufficiently large.
It follows from this  that 
 the probability that 
 \be     \sup_{t\in [n\de,(n+1)\de]}|  X^{1/2}(t)-  X^{1/2}( n\de)|\label{9.6tt} >   \ep(  \log n)^{1/2} \ee 
  infinitely often, is zero, for all $\ep>0$.
 Note that   
 \bea
   &&     \sup_{t\in [n\de,(n+1)\de]}  \frac{ X^{1/2}(t) }{(u^*(\de n,\de n) \log n)^{1/2}}\label{9.6as}\\
   &&\qquad \le \frac{ X^{1/2}(\de n) }{(u^*(\de n,\de n) \log n)^{1/2}}\nn+ \sup_{t\in [n\de,(n+1)\de]} \frac{|  X^{1/2}(t)-  X^{1/2}( n\de)|}{(u^*(\de n,\de n) \log n)^{1/2}}.
   \eea 
 It follows  that 
\begin{equation}
 \limsup_{t\to\ff} \frac{ X^{1/2} (t) }{   (u^{*}(t,t) \log  t)^{1/2} }=  \limsup_{n\to\ff}\sup_{t\in [n\de,(n+1)\de]}  \frac{ X^{1/2}(t) }{(u^*(\de n,\de n) \log n)^{1/2}}.
   \end{equation}
 Therefore, using (\ref{9.6as})   we see that
 \bea
          \limsup_{t\to\ff} \frac{ X^{1/2} (t) }{   (u^{*}(t,t) \log  t)^{1/2} }&\le& \limsup_{n\to\ff}\frac{ X^{1/2}(\de n) }{(u^*(\de n,\de n) \log n)^{1/2}}\label{5.14op}\\
          && + \limsup_{n\to\ff} \sup_{t\in [n\de,(n+1)\de]} \frac{|  X^{1/2}(t)-  X^{1/2}( n\de)|}{(u^*(\de n,\de n) \log n)^{1/2}}\nn.       
   \eea
  Writing $\log n=\log n\de+\log 1/\de$ we see from (\ref{8.37ee}) that the first term to the right of the inequality in (\ref{5.14op}) is less than or equal to 1 almost surely.     By (\ref{9.6tt}),  the second term to the right of the inequality in (\ref{5.14op})   is bounded by  $\ep/(u^{\ast}(T_{0},T_{0}))^{1/2}$  almost surely. Since this is true for all $\ep>0$ we get (\ref{9.10wr}).\qed

 \medskip	 \noindent	{\bf Proof of Theorem \ref{theo-3.2q} }     Let  $u(t,t)$ be regularly varying at infinity  with  index $\bb>0$. Let $ t_{n}=\th^{n}$, where $\th>1$  so that    $\wt S_{n}=\th^{n+1}-\th^{n}$. Since   
  \be
   \vf^{2}(t)\le O( u(t,t))\quad\mbox{as}\quad t\to\ff,\label{4.30}
  \ee
   we see that 
  \begin{equation}
   \vf^{2}(\wt S_{n})\le  Cu( \wt S_{n},\wt S_{n})\le  C(\th-1)^{\bb} u( \th^{n}, \th^{n})\quad\mbox{as $n\to\ff$},\label{4.35}
   \end{equation}  
 for some constant $C$. 
 
  Let $a_{n}= ((3+\ep) \log n )^{1/2}$. As in (\ref{9.6wq}), the probability that
 \be   
  \sup_{t\in [\th^{n},\th^{n+1}]}|  X^{1/2}(t)-  X^{1/2}(\th^{n})|\label{9.6wq8} >  \(  a_{n}  \vf(\wt S_{n})+\frac{\sqrt3}{\log 2}  \int_{0}^{1/n }  \frac{\vf (\wt S_{n} u)}{u(\log 1/u)^{1/2}}\,du  \)  
   \ee 
 infinitely often, is zero.   Note that by (\ref{4.35})
 \begin{equation}
 a_{n}  \vf(\wt S_{n})\leq  C(\th-1)^{\bb/2} (u( \th^{n}, \th^{n}) \log n )^{1/2}.\label{9.6wq9}
 \end{equation}
 
  We now show that this dominates the integral in (\ref{9.6wq8}). 
For all $n$ sufficiently large, the integral in  (\ref{9.6wq8}) is equal to
  \bea
\lefteqn{ \int_{0}^{\wt S_{n}/n }  \frac{\vf(  u)}{u(\log 1/u+\log \wt S_{n})^{1/2}}\,du\le\int_{0}^{1/2}  \frac{\vf( u)}{u(\log 1/u )^{1/2}}\,du}\label{3.104}\\\nn
&& +\int_ {1/2}^{1}  \frac{\vf( u)}{u(\log \wt S_{n} )^{1/2}}\,du + \int_{1}^{\wt S_{n}/n }  \frac{\vf(  u)}{u( \log \wt S_{n}-\log u)^{1/2}}\,du. \eea
   Using (\ref{3.21w}) we see that the first two integrals on the right-hand side of the inequality sign in (\ref{3.104}) are finite.  
In addition,  
\begin{equation}
    \int_{1}^{\wt S_{n}/n }  \frac{\vf(  u)}{u( \log \wt S_{n}-\log u)^{1/2}}\,du\le   \frac{ 1 }{(\log n)^{1/2}}\int_{1}^{\wt S_{n} }  \frac{\vf(  u)}{u }\,du.\label{3.104f}
   \end{equation}
Using (\ref{4.30}) and   (\ref{4.35}) and the regular variation of $u(t,t)$, we see that  (\ref{3.104f}) is  
 \bea
 &&   \le \frac{ C'   }{(\log n)^{1/2}}\int_{1}^{\wt S_{n}}  \frac{(u (x,x))^{1/2}}{ x  }\,dx\le \frac{ C''   }{(\log n)^{1/2}}( u(\wt S_{n},\wt S_{n}))^{1/2}\nn\\& & \le    \frac{\wt C''  }{(\log n)^{1/2}} ( \th-1)^{\bb/2}(u (\th^{n}, \th^{n}))^{1/2}.
  \eea
 Thus we see that the right-hand side of (\ref{9.6wq8})  is asymptotic to  (\ref{9.6wq9}) as $n\to\ff$. Consequently,
  \begin{equation}
    \limsup _{n\to\ff}  \sup _{t\in [\th^{n},  \th^{n+1}]}{|  X^{1/2}(t)-  X^{1/2}(\th^{n})| \over (u( \th^{n}, \th^{n}) \log n )^{1/2}} \leq C(\th-1)^{\bb/2}\qquad a.s. \label{9.6wq10}
  \end{equation}
    It follows from this and (\ref{8.37ee}),  with $t_{n}=\th^{n}$,   that 
  \begin{equation}
    \limsup _{n\to\ff}  \sup _{t\in [\th^{n},  \th^{n+1}]}{   X^{1/2}(t)  \over (u( \th^{n}, \th^{n}) \log n )^{1/2}} \leq 1+C(\th-1)^{\bb/2}\qquad a.s. \label{9.6wq11}
  \end{equation}
  Since $u(t,t)$ is regularly varying at infinity it is asymptotic to a monotonic function at infinity. Therefore,  
  \begin{equation}
    \limsup _{t\to\ff}  {   X^{1/2}(t)  \over (u( t, t) \log n )^{1/2}} \leq 1+C(\th-1)^{\bb/2}\qquad a.s. \label{9.6wq12}
  \end{equation}
Since this holds for all $\th>1$ we get (\ref{3.67}).\qed

  	 \noindent {\bf Proof of  upper bound in Theorem \ref{theo-1.9inf} }  
 Let  
\begin{equation}
\wh V_{\al ,f}(t)=\frac{\wh X_{\al ,f}(t)}{\wh u_{f}(t,t)},\qquad t\ge 0.\label{cup.1w}
\end{equation}
By Lemma \ref{lem-9.1},    $\wh V_{\al ,f}=\{\wh V_{\al ,f}(t), t\geq 0\}$  is    an  $\al$-permanental process with 
 sigma function  
\be 
   \wh \si_{f}^{2}(s,t)=2 -\,  \frac{ 2(\wh u_{f}(s,t) \wh  u_{f}(t,s))^{1/2}} { (\wh u_{f}(s,s)  \wh u_{f}(t,t))^{1/2}  }. \label{cup.3f}
   \ee 
   Note that 
   \bea 
   \frac{   \wh u_{f}(s,t) \wh u_{f}(t,s) } {  \wh u_{f}(s,s)  \wh u_{f}(t,t)  }  &=&   \frac{ e^{-\la|t-s|}+f(s)}{1+f(s)} \,\,\,  \frac{e^{-\la|t-s|}+f(t)}{1+f(t)} \nn\\
   &\ge&e^{-2\la|t-s|}. \eea
  Consequently
\begin{equation}
    \wh \si_{f}^{2}(s,t)\le 2(1-e^{-\la|t-s|})\le (\la|t-s|)\wedge 1:=\vf^{2} (|t-s|).
   \end{equation} 
 Therefore, it follows from Theorem \ref{theo-3.1w}  that 
      \begin{equation}
   \limsup_{ t\to\ff} \frac{ \wh V_{\al ,f}(t) }{\log t }\le 1\qquad a.s.\label{5.40}
   \end{equation}   
 This is (\ref{1.43rw}).
 
  For the last remark in this theorem, suppose that    $f$ is a potential for $\ov  B$, with  $h\in L_{+}^{1}\(0,\ff\)$. Then we have  
\be
f(t)=\int_{0}^{\ff} e^{-\la|t-s|}h(s)\,ds. 
\ee
For any $\ep>0$,  choose $s_{0}$ so that $\int _  {s_{0}}^{\ff} h(s)\,ds\leq \ep$. For $t\geq s_{0}$,
\bea
\label{pots.23}
f(t)&=& \int_{0}^{s_{0}} e^{-\la|t-s|}h(s)\,ds + \int _  {s_{0}}^{\ff} h(s)\,ds\\
&\leq& e^{-\la|t-s_{0}|} \,\|h\|_{1}+  \ep \nn.
\eea
Therefore, $\lim_{t\to\ff} f(t)\le \ep$ 
    for all $\ep>0$.  
 \qed

We now give some background material that may be needed to understand    Example \ref{ex-1.3}.
  A function  $f$ is excessive for $\wt  B$  if  and only if  $f$ is    positive and concave on $D$, which implies  that $f$ is increasing. This follows from the fact that
  $f$ is excessive for $\wt  B$  if  and only if $f$ is a positive
  superharmonic function on $D$,  \cite[Section 4.5, Theorem 3]{CW}.    That is,  $f$ is finite,  lower semi-continuous and midpoint concave, which implies that $f$ is concave,  \cite[Chapter I, Section 4.4, Corollary 1]{Bbk}.  
  It follows from this and Theorem \ref{theo-borelN} that for any positive  concave function $f$ on $D$ 
   \begin{equation}
 \wt u_{f}(s,t)=s  \wedge t+f(t),\qquad s,t>0,\label{1.39mmj}
  \end{equation}
is the kernel of an $\al$-permanental processes, for all $\al>0$.   In Theorem \ref{theo-1.10mm} we   denote  this process  by  $  \wt Z_{\al,f}  = \{\wt Z_{\al ,f}(t),t>0\} $. 

  Since $f$ is positive and increasing, we can define  $f(0)=\lim_{t\downarrow 0}f(t)$. It is easy to check that this extended function $f$ on $[  0,\ff)$ is positive and concave. 
Any positive  concave function $f$ on $[  0,\ff)$ can be written in the form  
\begin{equation}
f(t)=\wt f(t)+C_{0}t\label{1.39mmw},
\end{equation}
where $C_{0}\ge 0$ is a constant and $\wt f(t)$ is a positive  concave function that is $ o(t)$  at   infinity. 
To see this, note that $f_{r}'\(   t\)$, the right hand derivative of $f$, is decreasing in $t$. Let 
$C_{0}=\lim_{t\to \ff}f_{r}'\(   t\)$. We must have $C_{0}\ge 0$ since otherwise $f$ could not remain positive. Let $\wt f(   t)=f(t)-C_{0}t$.  This function $\wt f $ is concave and $\wt f_{r}'\(   t\)\geq 0$,  which implies  that  $\wt f$ is increasing. Since $\wt f\( 0\)=f(0)\geq 0$, we see that $\wt f$ is positive. 
In addition,  since  $ \lim_{t\to \ff}\wt f_{r}'\(   t\)=0$,  $\wt f$ is $ o(t)$  at   infinity.

\medskip	\noindent{\bf Proof of   upper bounds in Theorem \ref{theo-1.10mm} }
To prove (\ref{1.25ss})  it suffices to work with the  $\al$-permanental process   $\{\wt Z_{\al,f}(t),\,t\geq 1\}$. This process has  kernel $  \wt u_{f}(s,t)=s\wedge t+f(t)$, $s,t\geq 1$. Note that
$  \wt u_{f}(t,t)=t+f(t)$ is increasing, and by (\ref{1.39mmw}),   is regularly varying at infinity with positive index. We see from (\ref{3.74as}) in Lemma \ref{lem-3.8} that  when $f$ is  concave   the sigma function of $\wt Z_{\al,f} $ satisfies
\bea
  \si^{2}_{f}(s,t)&\le& |s-t| +|f(s)-f(t)|  \label{4.28}\\
&\le &\(1 +f_{l}'(s\wedge t)\)  |t-s|   \nn\\
&\le &\(1 +f_{l}'(1)\)  |t-s|  :=\vf^{2}_{f}( |t-s|)\nn,
  \eea
where $f_{l}'\(   x\)$ denotes the  left-hand derivative of $f$ at $x$. 

Clearly, $\vf^{2}_{f}( t)=O( \wt u_{f}(t,t))$ as $t\to\ff$. Therefore   the upper bound in  (\ref{1.25ss})  follows from Theorem \ref{theo-3.2q} and (\ref{1.39mmw}). 

  To obtain (\ref{1.25ssj}) note that when  
\be
f(t)=\int_{0}^{\ff}  (s\wedge t)h(s)\,ds\qquad   \mbox{  and }\, h\in L_{+}^{1},  
\ee then for all $\ep>0$,  
\bea
\label{pot.23}
\frac{f(t)}{t}&=&\frac{1}{t}\int_{0}^{\ep t}  (s\wedge t)h(s)\,ds +\frac{1}{t}\int _  {\ep t}^{\ff}  (s\wedge t)h(s)\,ds\\
&\leq& \ep  \,\|h\|_{1}+ \int_ {\ep t}^{\ff}  h(s)\,ds\nn.
\eea
Therefore,
\begin{equation}
 \lim_{t\to\ff} \frac{f(t)}{t}\le \ep  \,\|h\|_{1},
   \end{equation}
    for all $\ep>0$. This gives (\ref{1.25ssj}).\qed

 \noindent{\bf  Proof of  Theorem \ref{theo-cup}  }
Let $  \wt Z_{\al,f}  = \{\wt Z_{\al ,f}(t),t>0\} $ be the $\al$-permanental processes  with kernel $  \wt u_{f}(s,t)=s\wedge t+f(t)$, $s,t>0 $.  It follows from Remark \ref{rem-6.1},   with the isolated point   $\ast$ replaced by $0$, that there also exists an $\al$-permanental processes  that extends $ \wt Z_{\al,f}$ to $   \{\wt Z_{\al ,f}(t),t\ge 0\} $, with kernel  
\begin{eqnarray}
&&  \wt u_{f}(s,t)=s\wedge t+f(t), \hspace{.2 in}s,t>0
\label{rp.2q}\\
&&  \wt u_{f}(0,t)=f(t),  \quad t>0\hspace{.2 in}\mbox{ and }\hspace{.2 in}  \wt u_{f}(s,0)= \wt u_{f}(0,0)=1.
\nonumber
\end{eqnarray}   

Let  
\begin{equation}
\wh Z_{\al ,f}(t)=\frac{\wt Z_{\al ,f}(t)}{\wt u_{f}(t,t)},\qquad t\ge 0.\label{cup.1}
\end{equation}

Then by Lemma \ref{lem-9.1},    $\wh Z_{\al ,f}=\{\wh Z_{\al ,f}(t), t\geq 0\}$  is an $\al$-permanental process with 
sigma function
\begin{equation}
   \wh \si_{f}^{2}(s,t)=2 -\,  \frac{ 2(\wt u_{f}(s,t) \wt  u_{f}(t,s))^{1/2}} { (\wt u_{f}(s,s)  \wt u_{f}(t,t))^{1/2}  }=2 \(1-\,  \(\frac{ \wt u_{f}(s,t) } {    \wt u_{f}(t,t)  }\)^{1/2} \)\label{cup.3}
   \end{equation}
 when $0\leq s\leq t$. 
 
  Let $\vf_{f}(u)=u/f(u)$.  Note that since $\vf_{f}(u)/u=1/f(u)$, (\ref{ddd}) implies that   for some $\de>0$,
\begin{equation}
 \int_{0}^{\de }  \frac{ \vf_{f}( u) }{u(\log 1/u )^{1/2}}\,du<\ff.\label{5.49}
   \end{equation} 
  Changing variables, $u=e^{-s^{2}}$,  
 we write the integral in (\ref{5.49}) as  
\begin{equation}
  2\int_{(\log 1/\de)^{1/2}}^{\ff}\vf_{f}(e^{-s^{2}})\,ds.\label{5.50m}
  \end{equation}  
 Thus, since this integral is finite and $\vf_{f}(e^{-s^{2}})$ is decreasing, 	for any  $s> s_{0}$  
\begin{equation}
s  \vf_{f}(e^{-s^{2}}) \le   s_{0}  \vf_{f}(e^{-s^{2}}) +\int_{s_{0}}^{\ff}  \vf_{f}(e^{-u^{2}})\,du
  \end{equation} 
  Therefore,  
  \begin{equation}
\lim_{s\to\ff}s \vf_{f}(e^{-s^{2}}) \le    \int_{s_{0}}^{\ff}  \vf_{f}(e^{-u^{2}})\,du.
\ee
Since this holds for all $s_{0}$ we get  
\be
\lim_{s\to\ff}s  \vf_{f}\(e^{-s^{2}}\)=0,\label{5.53}
\ee
or equivalently
\begin{equation}
    \lim_{t\to 0}  \vf_{f}(t) (\log 1/t)^{1/2}=0.\label{5.50}
   \end{equation}
 
In particular,  (\ref{5.50}) implies that $\lim_{t\to 0}t/f(t)=0$. Therefore,   when  $0\leq s\leq t$, for all $\ep>0$, 
\begin{equation}
   \( \frac{ \wt u_{f}(s,t)  }{   \wt u_{f}(t,t)}\)^{1/2}= \( \frac{   s+f( t)  }{   t+f( t)  }\)^{1/2} \geq 1- (1+\ep)\frac{t-s}{f( t)} \qquad\mbox{as $\quad\,t\to 0$}.
   \end{equation}
  Therefore, 
\begin{equation}
 \wh \si_{f}^{2}(s,t)\le(1+\ep) \frac{|t-s|}{f( t)}\le (1+\ep)\frac{|t-s|}{f( t-s)}:=\wh \vf_{f}( t-s) \qquad\mbox{as $\quad\,s,t\to 0$}.\label{5.41q}
   \end{equation}

 \medskip	In preparation for using Lemma \ref{theo-8.1}  we first note that by concavity, $f(t)/t\leq f(s)/s$ for $t\geq s$ so that $\wh \vf_{f}( t)$ is increasing. We now show that (\ref{8.33}) holds.   If $\ga>1$, then, since $f$ is concave, we have  \begin{equation}
{   f(\ga V)-f(V) \over (\ga-1)V} \le  {f(  V) \over V}.
   \end{equation}
  Consequently  
  \be
 \frac{f(\ga V)}{f(  V)}-1\le \ga-1,
  \ee
  which gives (\ref{8.33}) when $\ga>1$.
If $\ga<1$,      
   \begin{equation}
{   f(V)-f(\ga V) \over (1-\ga)V} \le  {f( \ga  V) \over \ga V},
   \end{equation}
which implies that 
  \be
 \frac{f( V)}{f(\ga  V)}-1 \le \frac{1-\ga}{\ga},
  \ee
  which gives (\ref{8.33}) when $\ga<1$.

We can now use Lemma \ref{theo-8.1} to see  
	     \begin{equation}
   \limsup_{ t\to 0} \frac{|\wh Z^{1/2}_{\al ,f}(t)-\wh Z^{1/2}_{\al ,f}(0)|}{  \tau ( t)  }\le \sqrt 3\qquad a.s.\label{cup.4x} ,
   \end{equation} 
 where
  \begin{equation}
   \tau(t)=  \wh \vf_{f}(t)(\log\log 1/t)^{1/2}+\frac{1}{\log 2}\int_{0}^{1/2}\frac{\wh \vf_{f}(tu)}{u(\log 1/u)^{1/2}}\,du .\label{3.23q}
   \end{equation}  
   By Lemma \ref{lem-3.7}
    \be  \int_{0}^{1/2}  \frac{\wh \vf_{f}(tu)}{u(\log 1/u)^{1/2}}\,du\le 2\Theta(t)+o(\Theta(t)).\label{3.64z}
\ee
  where 
   \begin{equation}
   \Theta(t)=\int_{0}^{t^{2} }  \frac{\wh \vf_{f}( u) }{u(\log 1/u )^{1/2}}\,du+\wh \vf_{f}( t) (\log 1/t )^{1/2}.
   \end{equation}
By (\ref{5.49}) and (\ref{5.50})
   \begin{equation}
   \lim_{t\to 0} \Th(t)=0,\label{5.46}
   \end{equation}
  which gives 
\be
     \lim_{t\to 0}\tau ( t)  =0.\label{5.47q}
\ee  
  Since  $\wh Z_{\al ,f}(t)=\wt Z_{\al ,f}(t)/\(   t+f( t)\)$,   it follows from (\ref{cup.4x}) that\begin{equation}
\label{cup.5}
\lim_{t\to 0}\frac{\wt Z_{\al ,f}(t)}{f( t)}=\wh Z_{\al ,f}(0)\hspace{.2 in}a.s.\end{equation}
The theorem now follows from the fact that an $\al$-permanental random variable with kernel $1$, such as $\wh Z_{\al ,f}(0)$, is   a   gamma random variable $\xi_{\al,1}$ with shape $\al$ and scale $1$; see   \cite[(1.3)]{MRnec}.

  \qed

 \noindent{\bf  Proof of Theorem \ref{theo-cupx}  } If $f$ is regularly varying at 0 with index $0<\ga<1$,  $\wh \vf_{f}$ is regularly varying at 0 with index   $1-\ga $.  
By Theorem \ref{theo-8.1w} we get 
	     \begin{equation}
   \limsup_{ t\to 0} \frac{|\wh Z^{1/2}_{\al ,f}(t)-\wh Z^{1/2}_{\al ,f}(0)|}{  \wh \vf_{f}( t) (\log\log 1/t)^{1/2}}\le 1\qquad a.s.\label{cup.4w} 
   \end{equation} 
  Using the fact that  $\wh Z_{\al ,f}(t)=\wt Z_{\al ,f}(t)/\(   t+f( t)\)$, we get (\ref{cup.4}).\qed

\noindent{\bf  Proof of Theorem \ref{lem-3.9a} }   Using    Theorem \ref{theo-3.2} on the pairs  $\{t_{n},0\}$,  we see that   
 \begin{equation}
      \limsup_{n\to\ff}\frac{|X^{1/2}(t_{n} )-X^{1/2}(0)|}{ \si(t_{n} ,0)(\log n)^{1/2}}\le 1\qquad a.s.\label{3.73}
   \end{equation}
    
     When  $ \bb=\ff$, since $u(0,0)<\ff$ implies that $X(0)$ is finite almost surely, we immediately get (\ref{3.73q}).
   
 To consider the case when     $0<\bb<\ff$, we write 
 \begin{equation}
   |X(t_{n} )-X(0)|\le |X^{1/2}(t_{n} )-X^{1/2}(0)|\(|X^{1/2}(t_{n} )-X^{1/2}(0)|+2X^{1/2}(0)\)\label{prod.7}
   \end{equation}
   
  Using  (\ref{3.73}) we see that that
 \begin{equation}
   \limsup_{n\to\ff}\(|X^{1/2}(t_{n} )-X^{1/2}(0)|+2X^{1/2}(0)\)\le \bb+2X^{1/2}\label{3.75}(0)\qquad a.s.
   \end{equation}
   Combining this in (\ref{prod.7}) with (\ref{3.73}) we get (\ref{3.71}).
   
   When $\bb=0$, (\ref{3.75}) still holds and using (\ref{3.73}) we get (\ref{3.72}).\qed
   
   \begin{remark} {\rm Let $X=\{X(t ),t\in [0,1]\}$ be an $\al$-permanental sequence with kernel $u(s,t )$ and  sigma  function $\si(s,t )$. Consider the permanental sequence $Y=\{Y(n ),n\in \mathbb N \}=\{X(1/n), n\in \mathbb N\}$.  Obviously, the results in Theorem \ref{lem-3.9a} hold with $X(t_{n})$ replaced by $Y( {n})$ and $\si(t_{n},0)$ replaced by $\si(1/n, 0)$. 
   
  The advantage of Theorem  \ref{lem-3.9a} is that it doesn't require that (\ref{3.21w}) holds. This is   significant because   (\ref{3.21w}) requires that
  \begin{equation}
   \lim_{u\to 0} \si(u,0)(\log 1/u)^{1/2}=0.
   \end{equation}
(Under the additional assumptions that $\si=\vf$ and $\si(u,0)=\si(u)$.) Therefore, the results in  (\ref{8.35q}) would not give the results in Theorem \ref{lem-3.9a} when $\bb>0$. Additionally there are also some cases when $\bb=0$ but (\ref{3.21w}) does not hold.
}\end{remark}

   \begin{example}\label{ex-4.2} {\rm Let $ \{X_{\al}(0), X_{\al}(1), \ldots,X_{\al}( n),  \ldots\,n \in \mathbb N\}$ be an $\al$-permanental sequence determined by the kernel
   \be 
  u(0,0) =2,\qquad u( j,0)=1+f_{j},\qquad  u(0, k)=1+g_{k}, \qquad\mbox{$j,k= 1,\ldots$}\label{4.50}
   \ee 
   and
   \begin{equation}\label{4.51}
   u( j,k)= \la_{j}\de_{j,k}+1+f_{j}g_{k},\qquad\mbox{$j,k= 1,\ldots$}
   \end{equation}
 where    $\la_{j}\to 0$,  
 \begin{equation}
   0\le  f_{j}=1-p_{j}\le 1,\quad0\le  g_{j}=1-q_{j}\le 1,\quad   p_{j}=o(\la^{1/2}_{j}),\quad   q_{j}=o(\la^{1/2}_{j}) . 
   \end{equation} 
 It is easy to show that  the inverse of $\{u(j,k)\}_{j,k=1}^{m}$ is an $M$-matrix with positive row sums, which implies  that $u$ is the kernel of an $\al$-permanental sequences.  (See \cite{MRnonsym} for details.)

  We have  
      \bea
   \si^{2}(j,0)&=&u(j,j)+u(0,0)-2\(u(j,0)u(0,j)\)^{1/2}\\
  & =&\la_{j}+1+f_{j}g_{j}+2-2\((1+f_{j} )(1+ g_{j})\)^{1/2}\nn\\& =&\la_{j}+4-p_{j}-q_{j}+p_{j}q_{j} -2\((2-p_{j})(2-g_{j}) \)^{1/2}\nn
   \eea
  which gives   
   \begin{equation}
   \si^{2}(j,0)=\la_{j}+o(\la_{j}),\qquad\mbox{as}\quad j\to\ff.
   \end{equation}
 Using this,  it follows from   Theorem      \ref{lem-3.9a} that
   when  $\bb=\ff$,
  \begin{equation}
      \limsup_{n\to\ff}\frac{ X_{\al} (n ) }{ \la_{n}\log n }\le 1\qquad a.s.\label{3.73qm}
   \end{equation}
 when $0<\bb<\ff$,
 \begin{equation}
   \limsup_{n\to\ff} {| X_{\al} (n ) - X_{\al} (0 ) |} \le \bb^{2}+2\bb X_{\al}^{1/2}(0)\qquad a.s.\label{3.71m}
   \end{equation}
and when $ \bb=0$,
 \begin{equation}
   \limsup_{n\to\ff}\frac {| X_{\al} (n ) - X_{\al} (0 ) |}{(\la_{n}\log n)^{1/2}}\le  2  X^{1/2}\label{3.72m}(0)\qquad a.s.
   \end{equation}
    We show in \cite{MRnonsym} that the $\limsup$ in (\ref{3.73qm}) is actually equal to 1. 
   }\end{example}

    \section{Partial rebirthing of transient Borel right processes }\label{sec-Borel}

      Let $S$ a be locally compact set with a countable base. 
  Let     $X\!=\!
(\Om,  \FF_{t}, X_t,\th_{t},\newline P^x
)$ be a transient Borel right process with state space $S$,  and  continuous strictly positive  potential   densities $u(x,y)$ with respect to some $\si$-finite measure $m$ on $S$.   

Let $\ze=\inf\{t\,|\,X_t=\De \}$,  where $\De$ is the cemetery state for $X$,   and assume that $\ze<\ff$ a.s. Let $\mu$ be a finite measure on $S$. 
We call the   function    
\begin{equation}
f(y)= \int_{S} u(x,y)\,d\mu(x)\label{rp.1}
\end{equation}
 a left potential for $X$. Since $u(x,y) $ is continuous in $y$ uniformly in $x$ and $\mu$ is a finite measure
we see that $f(y)  $ is continuous. See \cite[Section 2]{FR}

The next theorem,   which is interesting on its own, is also  used in the proof of Theorem \ref{theo-borelN}. Note that   it does not require that $u$ is symmetric.  In this theorem we add a point   $*$ to the state space $S$ of $ X$ and modify $X$ so that instead of going to $\De$ it goes to $*$. We then   allow the process to return to $S$ from $*$ with a probability $p<1$, or to go to $\De$ with probability $1-p$. Let $\wt X$ denote the modified process on the enlarged space. We see that   when  $X$ ``dies'', $\wt X$ has a chance to be  reborn, after which  it continues to evolve that way $X$ did.

\bt\label{theo-borel}
Let $X\!=\!
(\Om,  \FF_{t}, X_t,\th_{t},P^x
)$ be a transient Borel right process with state space $S$,   as above. Then for any left potential $f$  for $X$, there exists a transient Borel right process
 $\wt X\!=\!
(\Om,  \FF_{t},\wt  X_t,\th_{t},\wt P^x
)$ with state space $\wt S=S\cup \{\ast\}$, where $\ast$ is  an isolated point, such that $\wt X$ has potential densities 
\begin{eqnarray}
&&  \wt    u(x,y)= u(x,y)+f(y), \hspace{.2 in}x,y\in S
\label{rp.2}\\
&&   \wt    u(\ast,y)=  f(y), \hspace{.2 in}\mbox{ and }\hspace{.2 in} \wt    u(x, \ast)=\wt    u(\ast, \ast) =1, 
\nonumber
\end{eqnarray}
with respect to the measure $\wt m$ on $\wt S$ which is equal to $m$ on $S$ and assigns 
a unit mass to $\ast$.
\et

   \Proof    We construct $\wt X$ as described prior to the statement of this theorem.  Let $\rho$ be the total mass of $\mu$. If  $\wt X$  starts in $S$
it proceeds just  like   $X$ until time   $\ze$, at which time   it goes  to $\ast$. It stays there for an independent exponential time with parameter $1+\rho$, $\rho>0$,  after which it   returns to  $S$ with   initial law ${\mu/(1+\rho)}$. (This is what we mean by partial rebirthing.)   

 Once in $S$, $\wt X$ continues as we just described for $\wt X$ starting in $S$. Since the measure ${\mu/(1+\rho)}$ has total mass ${\rho /( 1+\rho)}$, after each visit to $\ast$,  $\wt X$   only has probability ${\rho /( 1+\rho)}$ to be reborn. With probability $  {1 /(1+\rho)}$ the process enters a cemetery state   $ \De$.   

We now calculate the   potential densities for  $\wt X$. Let $g$ be a function on $\wt S$ with $g(\ast)=0$.
Then for any $x\in S$

\begin{eqnarray}
&& \hspace{-.1in}E^{x}\(\int_{0}^{\ff}g\(\wt X_{t}\)\,dt\)\label{rp.3}\\
&&\quad=E^{x}\(\int_{0}^{\ze}g\(X_{t}\)\,dt\)+\sum_{n=1}^{\ff}\({\rho \over 1+\rho}\)^{n-1}\int 
{\,d\mu(z) \over 1+\rho}E^{z}\(\int_{0}^{\ze}g\(X_{t}\)\,dt\).
\nn
\eea
This is equal to
\bea
&&\hspace{-.2in}  \int u(x,y)g(y)\,dm(y) +{1\over 1+\rho} \sum_{n=1}^{\ff}\({\rho \over 1+\rho}\)^{n-1}
\int \,d\mu(z)\int u(z,y)g(y)\,dm(y)  \nonumber\\
&& =\int u(x,y)g(y)\,dm(y) + 
 \int f(y)g(y)\,dm(y),  \nonumber
\end{eqnarray}
which gives the first line of (\ref{rp.2}). The first half of the second line of (\ref{rp.2}) follows from a similar  computation, where now we no longer have the first term in the second line of (\ref{rp.3}). Finally, since at each visit to $\ast$ the process waits there an independent exponential time with parameter $1+\rho$, and then returns to $\ast$ with probability ${\rho /1+\rho)}$,  we have, for some sequence of functions $h^{ ( n)}\to h $.
\be E^{x}\(\int_{0}^{\ff}1_{\{\wt X_{t}=\ast\}} \,dt\)
\label{rp.4} {1\over 1+\rho} +\sum_{n=1}^{\ff}\({\rho \over 1+\rho}\)^{n} {1\over 1+\rho}=1. \ee
The same computation holds if we start at $\ast$.\qed

  The next lemma is used in the proof of Theorem \ref{theo-borelN}.

 \begin{lemma} \label{lem-6.1} Assume that   for each $n\in \mathbb N$,   $u^{ ( n)}(s,t),\,s,t\in S $, is the kernel of an $\al$-permanental process.  If
  $  u^{ ( n)}(s,t)\to   u(s,t)$ for all $s,t\in S$, then $  u(s,t) $ is the kernel of an $\al$-permanental process.
  \end{lemma}
 
\Proof  By the hypothesis,  for all $k$ and  $x_{1},\ldots, x_{k}\in S$,   there exists  an $\al$-permanental vector 
    $\(X^{\(n\)}_{\al}(x_{1}),\ldots, X^{\(n\)}_{\al}(x_{k})\)$ with kernel $K^{\(n\)}_{i,j}=u^{\(n\)}(x_{i},x_{j}) $, $1\leq i,j\leq k$.   Therefore, by definition, for all $s_{1},\ldots, s_{k}\geq 0$,
  \begin{equation}
   E\(e^{-\sum_{i=1}^{k}s_{i}X^{\(n\)}_{\al}(x_{i})}\) 
 = \frac{1}{ |I+K^{\(n\)}S|^{ \al}}.   \label{int.1exy} 
 \end{equation}  
In addition, since  $u^{\(n\)}(x_{i},x_{j}) \to u (x_{i},x_{j}) $, we have     $ |I+K^{\(n\)}S|\to  |I+KS|$, where $K_{i,j}=u(x_{i},x_{j}) $, $1\leq i,j\leq k$. 
It follows from the extended continuity theorem  \cite[Theorem 5.22]{K}, that there exists a random vector $(X_{\al}(x_{1}),\ldots, \newline	X_{\al}(x_{k}))$ with 
\begin{equation}
   E\(e^{-\sum_{i=1}^{k}s_{i}X_{\al}(x_{i})}\) 
 = \frac{1}{ |I+KS|^{ \al}}.   \label{int.1ex2y} 
 \end{equation}
Since this is true for all $k$ and all  $x_{1},\ldots, x_{k}\in S$,  it follows by the Kolmogorov extension theorem
that    $ \{u(s,t) ,s,t\in S\}$ is the kernel of an $\al$-permanental process. \qed

 	 \noindent{\bf  Proof of Theorem \ref{theo-borelN} } We   apply Lemma \ref{lem-6.1} twice to prove the theorem. Consider  a general excessive function $f$. It follows from \cite[II, ( 2.19)]{BG} that  there exists a sequence of functions $g_{n}\geq 0$, with both  $g_{n}$ and  
\begin{equation}
Ug_{n}(x)= \int_{S} u(x,y)g_{n}(y)\,dy\label{rp.11}
\end{equation}
  bounded such that $f(x)$
is the increasing limit of  $Ug_{n}(x)$.   If the $g_{n}$ are in $L^{ 1}$ then,  since $u$ is symmetric,  $Ug_{n}$ is a left potential  as in (\ref{rp.1}). Therefore, by Theorem \ref{theo-borel}, $\{ u(s,t)+U  g_{n}  ( t),s,t\in S\}$ are kernels of  $\al$-permanental   processes.  Consequently,  by Lemma \ref{lem-6.1}, $\{ u(s,t)+f( t),s,t\in S\}$  is the kernel  of  $\al$-permanental process.

If 
 $g_{n}$ is not   in $L^{ 1}$ we proceed as follows: Let $C_{m}$ be an increasing sequence of compact sets whose union is $S$. Then $g_{n}1_{C_{m}}\in L^{ 1}$, so that    by Theorem \ref{theo-borel} $\{ u(s,t)+U\( g_{n}1_{C_{m}}  \)( t),s,t\in S\}$ is the kernel of an $\al$-permanental process. Taking the limit as $m\to \ff$, it follows from  Lemma \ref{lem-6.1} that  $ \{u(s,t)+U  g_{n} ( t),s,t\in S\}$ is the kernel of an $\al$-permanental process. Since $U  g_{n}\to f$ we can use Lemma \ref{lem-6.1} again to see that  $\{ u(s,t)+f( t),s,t\in S\}$   is the kernel  of  $\al$-permanental process.
 
 \qed

\begin{remark}\label{rem-6.1}
   Theorem \ref{theo-borelN}   shows that there exists an $\al$-permanental process $Z_{\al}( t), t\in S$ with the kernel given in (\ref{1.10}).  The same proof  also shows that there exists an $\al$-permanental process $\{   Z_{\al}( t), t\in S\}\cup Z_{\al}(\ast)$ with the kernel given in (\ref{rp.2}) for any function $f$ which is excessive for $X$.     
   \end{remark}

  		  \section{Lower bounds}  \label{sec-7}

We use results from \cite{MRnec} to obtain  lower bounds for the rate of growth of permanental process or for their behavior at 0. There are several different situations that can arise	depending on the kernels of the permanental processes. We give several criteria that can be used on kernels that behave differently. 
 	\begin{lemma} \label{lem-5.1}   Let $X_{\al}=\{X_{\al}(t ), t\in R^{+}\}$  be an $\al$-permanental process with  kernel $u(s,t)$ such that $u(t,t)>0$ for all $t\in R^{+}$. Set
\begin{equation}
 \wt u(s,t)  = \frac{u(s,t)}{(u(s,s)u(t,t))^{1/2}}\qquad s,t\in R^{+}.
   \end{equation}
  Let $\{t_{j}\}_{j=1}^{\ff}$ be a sequence in $R^{+}$. Set 
\begin{equation}
   \phi^{2}(i,j)= 2-   ( \wt u(t_{i},t_{j})+ \wt u(t_{j},t_{i}) )  \quad\mbox{and}\quad (\phi_{n}^{*})^{2}= \inf_ {\stackrel{1\le i,j\le n} {i\ne j}  }\phi^{2}(i,j).
   \end{equation}
If
\begin{equation}
   \sup_ {\stackrel{1\le i,j\le \ff} {i\ne j}  }   \wt u(t_{i},t_{j})+ \wt u(t_{j},t_{i})   \le \ep_{1}\label{10.3q}
   \end{equation}
  and  
  \begin{equation}
   \sup _ {\stackrel{1\le i,j\le n} {i\ne j}  }  | \wt u(t_{i},t_{j})- \wt u(t_{j},t_{i})  | \le \ep_{2}   (\phi_{n}^{*})^{2}\label{10.5},
   \end{equation}
  for $\ep_{1}, \ep_{2}$ sufficiently small, then 
  \begin{equation}
   \limsup_{i\to\ff} \frac{X_{\al}(t_{i})}{u(t_{i},t_{i})\log i}\ge 1  -3(\ep_{1}+\ep_{2})\qquad  {a.s.}\label{5.5}
   \end{equation}	
    \end{lemma}	
  
  \Proof   Let $ \{\wt X(t )=X(t )/u(t,t), t\in R^{+}\}$.
  We show in Lemma \ref{lem-9.1} that   $\wt X= \{\wt X(t ), t\in R^{+}\}$ is an $\al$-permanental process with kernel $ \wt u(s,t) $.    
Now consider the matrix $K_{n}=\{ \wt u(t_{i},t_{j})\}_{i,j =1}^{n}$. This is the kernel of the $\al$-permanental vector $(\wt X(t_{1}),\ldots,\wt X(t_{n}) )$. Let  $\{a_{i,n}\}_{i=1} ^{n}$ denote the diagonal elements of $K_{n}^{-1}$.  By (\ref{10.3q})
   \begin{equation}
    (\phi_{n}^{*})^{2}\ge 2-\ep_{1}.
   \end{equation} 
By (\ref{10.5}) we can take $C=\ep_{2}$ in \cite[(5.5)]{MRnec} and, since $ \wt u(t,t)\equiv 1$, use \cite[Lemma 5.2]{MRnec}  to get   
\begin{equation} 
   a_{i,n}\le \frac{2}{(1-\ep_{2})(\phi_{n}^{*})^{2}}\le \frac{2}{(1-\ep_{2})(2-\ep_{1})}\le 1+2(\ep_{1}+\ep_{2}),
   \end{equation}
for all $1\le i\le n$ and all $\ep_{1},\ep_{2}$ sufficiently small.

\medskip	 To complete the proof we use the next two lemmas. 
 
\begin{lemma} \label{lem-4.2}Let $\{\xi^{(i)}_{u,v}\}_{i=l}^{n}$ be independent   copies of $\xi_{u,v}$.   (See (\ref{1.2q})). Then  for all   $0<\ep<1$,  and  $l\ge l_{0}=l_{0}(\ep)$ with  $ (2l_{0}^{\ep}/(3  \,\Gamma(u)  \log l_{0}) )\ge 1$,
 \begin{equation}
   P\(\max_{l\le i\le n}{ \xi^{(i)}_{u,v  } \over \log i}\ge  \frac{ (1-\vep)}{v }\)\ge 1- {l+1 \over n+1}.\label{4.5}
   \end{equation}
 \end{lemma}
 
  \Proof  
We have  
  \bea
   P\(\max_{l\le i\le n} { \xi^{(i)}_{u,v  } \over \log i}> \frac{(1-\ep) }{v  }\)\   &=& 1-P\(\max_{l\le i\le n} { \xi^{(i)}_{u,v  } \over \log i} \leq  \frac{(1-\ep) }{v }\)\label{4.8}\\& =&1-\prod_{i=l}^{n }\(1- P\({ \xi^{(i)}_{u,v  } \over \log i}> \frac{(1-\ep) }{v  }\) \)\nn.
   \eea
  For any $i\geq l_{0}$,
\begin{equation}
   P\( { \xi^{(i)}_{u,v  } \over \log i}> \frac{(1-\ep) }{v  }\)\ge 
  \frac{2e^{-(1-\ep)\log i}} {3\Gamma(u)(1-\ep)\log i}\ge \frac{1}{i}.\label{7.11}
   \end{equation}
  (See e.g.     \cite[(3.2)]{MRnec}.)   Using   (\ref{4.8}) and (\ref{7.11}), we see that   
    \begin{equation}
P\(\max_{l\le i\le n}{ \xi^{(i)}_{u,v  } \over \log i}> \frac{(1-\ep) }{v  }\)\ge  1- \prod_{i=l}^{n }\(1-  \frac{1}{i}  \) >1-e^{-\sum_{i=l}^{n }1/i} \nn.
   \end{equation}
   \qed

 	  \begin{lemma}  \label{lem-4.5}Let $X=(X_{1}, X_{2}, \ldots)$ be an $\al$-permanental sequence, and for each $n$, let  $K(n)$ be the kernel of $X=(X_{1}, X_{2}, \ldots, X_{n})$. If   $(K(n))^{-1}$ has diagonal elements $a_{n, i}\leq a_{i}$,  $i_{0}\leq i\leq n-i_{0}$, for some $i_{0}$ and  all  $n$ sufficiently large, then
            \be
 \limsup_{  i\to\ff }  { X_{ i}\over a^{-1} _{i }\log i} \ge 1, \hspace{.2 in}a.s.\label{4.10new2}   
   \ee   
   \end{lemma}

  \Proof   
Using   \cite[(1.7)]{MRnec}  and then (\ref{4.5}) we see that for any $\ep>0$ and $l\geq l_{0}\vee i_{0}(\ep)$  
  \bea
    P\(\max_{l\le i\le n-i_{0}}  { X_{ i}\over a^{-1} _{i }\log i} \ge 1-\vep \)&\geq &  P\(\max_{l\le i\le n-i_{0}}  { X_{ i}\over a^{-1} _{n,i }\log i} \ge 1-\vep \)\nn\\
  &\ge &   P\(\max_{l\le i\le n -i_{0}}  { \xi^{(i)}_{\al,1 } \over  \log i} \ge 1-\vep\)\\&\ge & 1-{l+1 \over n-i_{0}+1}.\label{4.10qj}   \nn
  \eea
It follows from this that for any $\ep>0$ and $l\geq l_{0}(\ep)\vee i_{0}$,
         \be
   P\(\sup_{l\le i }  { X_{ i}\over a^{-1} _{i }\log i} \ge 1-\vep \)= 1.\label{4.10qm}   
   \ee
We take the limit as  $l\to\ff$ and use monotone convergence to get
               \be
   P\(\limsup_{  i\to\ff }  { X_{ i}\over a^{-1} _{i }\log i} \ge 1-\ep \)=1.\label{4.10qp}   
   \ee
Since this holds for all $\ep>0$ we obtain (\ref{4.10new2}).
  \qed 
  
   \noindent{\bf Proof of Lemma \ref{lem-5.1} continued:  } It follows from Lemma \ref{lem-4.5} that

 \begin{equation}
   \limsup_{i\to\ff} \frac{\wt X(t_{i})}{\log i}\ge  (1-3(\ep_{1}+\ep_{2})) \qquad  {a.s.}
   \end{equation}
   This gives (\ref{5.5}).  \qed

\noindent{\bf Proof of  lower bounds in Theorem \ref{theo-1.8new} }  We obtain the lower bounds in   (\ref{8.35sw5}) and (\ref{3.67ss}).   Following Lemma \ref{lem-5.1} set  
 \begin{equation}
\wt u_{T_{0},\ga,\bb}(s,t)  = \frac{u_{T_{0},\ga,\bb}(s,t)}{(u_{T_{0},\ga,\bb}(s,s)u_{T_{0},\ga,\bb}(t,t))^{1/2}}\qquad s,t\in R^{+}.
   \end{equation}
  Let   $ t_{j}=\th^{j} $. Then  
\begin{equation}
\wt u_{T_{0},\ga,\bb}(\th^{i},\th^{j})+\wt u_{T_{0},\ga,\bb} (\th^{j},\th^{i})  =  \frac{  R_{\ga,\bb}(\th^{i},\th^{j})}{\(R_{\ga,\bb}(\th^{i},\th^{i})R_{\ga,\bb}(\th^{j},\th^{j})\)^{1/2}}\label{10.3qq}.
   \end{equation}
   Using (\ref{9.6}) it is easy to check that for $i\ne j$, this is 
   \begin{equation}
   \le  C \th^{-\ga/2}\quad\mbox{for $\th\gg1$}\qquad\mbox{and}\qquad \le C'\th^{ \ga/2}\quad\mbox{for $\th\ll1$},
   \end{equation}
  for constants $C$ and $C'$. Similarly \begin{equation}
\wt u_{T_{0},\ga,\bb}(\th^{i},\th^{j})-\wt u_{T_{0},\ga,\bb}(\th^{j},\th^{i})  =  \frac{2|H_{\ga,\bb}(\th^{i},\th^{j})|}{\(R_{\ga,\bb}(\th^{i},\th^{i})R_{\ga,\bb}(\th^{j},\th^{j})\)^{1/2}}.\label{10.3ww}
   \end{equation}Using (\ref{9.6a}) we see that this is 
    \begin{equation}
   \le |\bb|\th^{-\ga/2}\quad\mbox{for $\th\gg1$}\qquad\mbox{and}\qquad \le  |\bb|\th^{ \ga/2}\quad\mbox{for $\th\ll1$}.
   \end{equation} 
    Therefore,    (\ref{5.5}) holds    for $\th\gg1$ or $\th\ll1$,  and since $\wt u_{T_{0},\ga,\bb}(\th^{n},\th^{n})=2C_{\al,\bb}\th^{\al n}$ we get 
   \begin{equation}
   \limsup_{n\to\ff} \frac{X(\th^{n})}{  \th^{\al n} \log n}\ge 2C_{\al,\bb}( 1-\ep)  \qquad  \label{10.11q}{a.s.}
   \end{equation}
   where $\ep$ depends on $\th$ and goes to 0 as $\th$ goes to 0 or $\ff$, depending on whether $\th\ll 1$ or $\th\gg 1$.  Using the facts that for $\th\gg 1$, $\lim_{n\to\ff}\log n /\log\log \th^{n}=1$ and for $\th\ll 1$, $\lim_{n\to\ff}\log n /\log\log \th^{-n}=1$
  we get the lower bounds in  (\ref{8.35sw5}) and (\ref{3.67ss}).    \qed
  
   \noindent {\bf Proof of  lower bound in Theorem \ref{theo-1.9inf} } This is an immediate application of Lemma \ref{lem-5.1}. 
 Consider  $  \{\wh X_{\al,f}(n j ), j\in  \mathbb N\}$. It is easy to see that    \begin{equation}
   \sup_{ {1\le j,k\le n} \atop{j\ne k}} \wt u_{f}(nj,nk)=    \sup_{ {1\le j,k\le n} \atop{j\ne k}} f(nk)+e^{-\la n|k-j|} 
   \end{equation}
  Therefore, since $\lim_{t\to\ff}f(t)=0$, for all $\ep>0$ we can choose $n$ such that (\ref{10.3q}) and (\ref{10.5}) hold with $\ep_{1}$ and $\ep_{2}$ less that $\ep$. Consequently, (\ref{1.43ro}) follows from (\ref{5.5}).\qed

 \medskip	\noindent{\bf Proof of  lower bound in Theorem \ref{theo-1.10mm} }  
 Let   $U $ be a non-singular $n\times n$ matrix. We use  $U^{-1}$ to denote the inverse, and $U^{j,k}$ to denote the elements of $U^{-1}$.

Let  $U_{f}$ be the $(n+1)\times (n+1)$ matrix   
   \bea \hspace{-.3 in} U_{f}&=&\left (
\begin{array}{ cccc } 1 &f(1 )&\ldots&f(n )  \\
1   &U_{1 ,1 }+f(1 )&\ldots&U_{1 ,n }+f(n )  \\
\vdots& \vdots &\ddots &\vdots  \\
1   &U_{n ,1 }+f(1 )&\ldots&U_{n ,n+ }+f(n )  
\end{array}\right ).\label{19.39}
    \eea 
One can check that 
\be  U_{f}^{-1} =\left (
\begin{array}{ cccc  }  1+\rho &-\sum_{i=1}^{n}f(i)U ^{i ,1 }   &\dots&-\sum_{i=1}^{n}f(i)U ^{i ,n }  \\
- \sum_{j=1}^{n}U ^{1 ,j } & U ^{1 ,1 } & \dots &  U ^{ 1 ,n }    \\
\vdots&\vdots& \ddots&\vdots  \\
 - \sum_{j=1}^{n}U ^{n ,j }& U ^{n ,1 }&  \dots & U ^{n ,n }   \end{array}\right ),\label{19.40}
    \ee 
    where
    \begin{equation}
\rho =:\sum_{j=1}^{n}\sum_{i=1}^{n}f(i)U ^{i ,j} .\label{19.10v}
\end{equation}
We now apply this with $U$ replaced by  $  W(n)=\{ s_{j}\wedge s_{k}\}_{j,k=1}^{n} $, where $t_{i}>0 $, $i=1,\ldots,n$, and 
\begin{equation}
   s_{j}=t_{1}+\cdots+t_{j},\qquad j=1,\ldots,n.
   \end{equation}
That is,
  \be W(n) = \left (
\begin{array}{ cccccc cc}  
s_{1} &s_{1}&s_{1}&\dots &s_{1} &s_{1}  \\
s_{1} &s_{2}&s_{2}&\dots &s_{2}&s_{2}  \\
\vdots&\vdots&\vdots&\ddots&\vdots &\vdots \\
s_{1} &s_{2}&s_{3}&\dots &s_{n-1} &s_{n-1}   \\
s_{1} &s_{2}&s_{3}&\dots &s_{n-1}&s_{n} \end{array}\right ).\label{harr2.211q}
     \ee 
It is easy to check that  
  \be W(n)^{-1}= \left (
\begin{array}{ cccccc cc}  
\frac{1}{t_{1}}+\frac{1}{t_{2}}&-\frac{1}{t_{2}}&0&\dots &0 &0  \\
-\frac{1}{t_{2}}&\frac{1}{t_{2}}+\frac{1}{t_{3}}&-\frac{1}{t_{3}}&\dots &0&0  \\
\vdots&\vdots&\vdots&\ddots&\vdots &\vdots \\
0&0&0&\dots &\frac{1}{t_{n-1}}+\frac{1}{t_{n}}&-\frac{1}{t_{n}}    \\
0&0&0&\dots &-\frac{1}{t_{n}}&\frac{1}{t_{n}} \end{array}\right ).\label{harr2.211}
     \ee

 Now, let 
\begin{equation}
   U_{f}(n+1)=\{ s_{j}\wedge s_{k}+f(k)\}_{j,k=0}^{n} 
   \end{equation}  
   where $s_{0}=0$  and $f(0)=1$. 
It follows from (\ref{19.40}) that all the diagonal entries of $U_{f}(n+1)^{-1}$, except for the first one, are equal to diagonal entries of $W(n)^{-1}$, that is they are equal to $ {1}/{t_{j}}+{1}/{t_{j+1}}$, $j=1,\ldots, n-1$ and $ {1}/{t_{n}} $.

       Let $s_{j}=\th^{j}$, for $\th\gg 1$. Then $t_{j}=\th^{j}-\th^{j-1}$ and  
\begin{equation}
 U_{f}(n+1)^{j,j}=  \frac{1}{\th^{j}}\(\frac{\th+1}{\th-1}\),\qquad j=2,\ldots,n-1,\label{6.24mm}
   \end{equation}
   and 
   \begin{equation}
 U_{f}(n+1)^{n,n } = \frac{1}{\th^{n}}\({\th \over  \th -1}\).\label{6.24m}
   \end{equation}
  It now follows from Lemma \ref{lem-4.5}  that
    \be   \limsup_{n\to \ff}\frac{\wt Z_{\al,f } (\th^{n}) }{\th^{n}\,  \log n}\ge \frac{\th-1}{\th+1} \qquad a.s.\label{6.26mm}
   \end{equation}
  Taking $\th$ arbitrarily large gives the lower bound in (\ref{1.25ss}).\qed

 \section{Appendix}\label{Appendix}

  In this section  we simply write $d$ for the metric  $ d_ {C ,\si}$ in (\ref{1.53}).  
  
  \medskip	Let $(\TT,d)$   be aÊ separable metric or pseudo-metric space. Let $B_{d}(t,u)$  denote a closed ball of radius $u$ in $(\TT,d)$ and   $\mu$ a probability measure on $\TT$ we define
\begin{equation}
  J_{\TT, d,\mu}( a) =\sup_{t\in \TT}\int_0^a\(\log\frac1{\mu(B_{d}(t,u))}\)^{1/2} \,du.\label{tau}
   \end{equation}

The next theorem follows  from \cite[Theorem 3.1]{KMR}. The proof of \cite[Theorem 3.1]{KMR} is a consequence of the fact that a 1/2-permanental process is subgaussian. Using Theorem \ref{theo-1.3} it extends it as follows:

 \bt\label{theo-1.1a}Ê Ê  Let $X_{\al}=\{X_{\al}(t),t\in \TT\}$ be an $\al$-permanental process with kernel $u(s,t)$.  
 Assume   that $\TT$ is separable for $  d $ with finite diameter $D$ and   that there exists
a probability measure $\mu$Ê Ê on $\TT$ such that  
\be
J_{\TT,  d,\mu }(D)<\ff\label{1.1v}.
\ee 
Then there exists a version $X'_{\al}=\{X'_{\al}(t),t\in \TT\}$ of $X_{\al}$ which is bounded almost surely.
 
 If
\begin{equation}
   \lim_{\de\to 0}J_{\TT,  d,\mu }(\de)=0,\label{1.8jv}
   \end{equation}
there exists a version $X'_{\al}=\{X'_{\al}(t),t\in \TT\}$ of $X_{\al}$ such that
 \be     \lim_{\de\to 0}\sup_{\stackrel{s,t\in \TT}{ d(s,t)\le \de}} |X '_{\al}( s)-X' _{\al}( t)|\label{3.6.90hhv} =0,\qquad a.s.
\ee 

If (\ref{1.8jv}) holds and 
\begin{equation} 
 \lim_{\de\to 0}{J_{  \TT, d,\mu }(\de)\over \de}=\ff
  \label{2.18v},
   \end{equation}
   then  
   \begin{equation}
   \lim_{\de\to0}   \sup_{\stackrel{s,t\in \TT}{  d (s,t)\le \de}}\frac{ |X '_{\al}( s)-X' _{\al}( t)|}{J_ {\TT,  d,\mu } (  d(s,t) /2) }\le 30\(\sup_{t\in \TT}X'_{\al}(t)\)^{1/2}\quad a.s.\label{2.1sv}
   \end{equation}
   \et
   
   The next theorem follows  from \cite[Theorems  4.2]{MRsuf} and  Theorem \ref{theo-1.3}.  
  
\begin{theorem}\label{theo-1.1au} Under the hypotheses of Theorem \ref{theo-1.1a} assume that 
(\ref{1.8jv}) holds. For any $t_{0}\in \TT$   and $\de>0$, let $\TT_{\de}:=\{s:  d(s,t_{0})\le \de/2\}$.  Suppose   $0<\de \le \de_{0}<D$ which implies that $\TT_{\de}\le \TT_{D}$. Assume that for some $\bb<1$
\begin{equation}
   \limsup_{k\to\ff}\frac{\mu(\TT_{\bb^{k}})}{\mu(\TT_{\bb^{k+1}} )}\le C ,\label{2.27}
   \end{equation}
for some constant $C$, and consider the probability measures $\mu_{\de}(\cd):=\mu(\cd\cap \TT_{\de})/\mu(\TT_{\de})$, $0<\de\le \de_{0}$. 
 Then if $X_{\al}(t_{0})\ne 0$     there exists a version $X'_{\al}=\{X'_{\al}(t),t\in \TT\}$ of $X_{\al}$  such that  
\begin{equation}
   \lim_{\de\to0}   \sup_{   d(s,t_{0})\le \de/2 }\frac{ |X'_{\al}(s)-X'_{\al}(t_{0})|}{\ov H_{\TT_{\de},d,\mu_{\de}} (\de/4)}\le C\, X'^{1/2}_{\al}(t_{0})   \qquad a.s.\label{2.1ww}
   \end{equation}
   where  
     \be 
 { \ov H_{\TT_{\de},d,\mu_{\de}}(\de/4)\label{2.20aa}} :=\de(\log \log 1/\de)^{1/2}+J_{\TT_{\de},d,\mu_{\de}} (\de/4).
   \ee 
 
   If $X_{\al}(t_{0})= 0$     there exists a version $X'_{\al}=\{X'_{\al}(t),t\in \TT\}$ of $X_{\al}$  such that  

   \be  
      \lim_{\de\to0}   \sup_{ \hat d(s,t_{0})\le \de/2 }\frac{  X'_{\al}(s) }{\(\ov H_{\TT_{\de},d,\mu_{\de}} (\de/4) \)^{2}}\le C' \qquad a.s.,\label{2.1qqq}
   \end{equation}
 for some constant $C'$.
    \end{theorem}

 \begin{remark} {\rm We have pointed out on page \pageref{p3}, that when $\{u(s,t);s,t\in \TT\}$ is the potential density of a transient Markov process,  $\{d_{C,\si}(s,t);s,t\in \TT\}$ defined in (\ref{1.2}) and (\ref{1.53}), is a metric on $\TT$. In general, if we only assume that $\{u(s,t);s,t\in \TT\}$ is a kernel of $\al$-permanental processes, we don't know whether $d_{C,\si}$ is a metric. Actually Theorems \ref{theo-1.1a} and \ref{theo-1.1au} still hold if $\{d_{C,\si}(s,t);s,t\in \TT\}$ is not a metric.   We continue to define 
 \begin{equation}
   B_{d_{C,\si}}(t,u)=\{s; d_{C,\si}(s,t)\le u\}
   \end{equation}
and everything goes through. (This is the approach we took in \cite{MRsuf} which we wrote before we knew that  when $\{u(s,t);s,t\in \TT\}$ is the potential density of a transient Markov process,  $\{d_{C,\si}(s,t);s,t\in \TT\}$   is a metric on $\TT$. See, in particular, \cite[Theorems 1.1 and 1.2]{MRsuf} and the paragraph that preceeds \cite[Theorems 1.1 ]{MRsuf}.)
 }\end{remark}

\def\noopsort#1{} \def\printfirst#1#2{#1}
\def\singleletter#1{#1}
            \def\switchargs#1#2{#2#1}

\def\bibsameauth{\leavevmode\vrule height .1ex
            depth 0pt width 2.3em\relax\,}
\makeatletter
\renewcommand{\@biblabel}[1]{\hfill#1.}\makeatother
\newcommand{\bysame}{\leavevmode\hbox to3em{\hrulefill}\,}

\noindent
\begin{tabular}{lll} & \hskip20pt Michael B.  Marcus
     & \hskip20pt  Jay Rosen\\  & \hskip20pt 253 West 73rd. St., Apt. 2E
   & \hskip20pt Department of Mathematics \\    &\hskip20pt  New York, NY 10023, USA
 \hskip20pt 
     & \hskip20pt College of Staten Island, CUNY \\    & \hskip20pt mbmarcus@optonline.net
    & \hskip20pt  Staten Island, NY
10314, USA \\    & \hskip20pt      & \hskip20pt   jrosen30@optimum.net \\   & \hskip20pt 

     & \hskip20pt 
\end{tabular}

\end{document}